\documentclass[12pt]{article}
\usepackage{amsthm,amsfonts,amssymb,amscd}

\begin{document}

\begin{center}
\large \bf Birationally rigid Fano-Mori fibre spaces
\end{center}\vspace{0.5cm}

\centerline{A.V.Pukhlikov}\vspace{0.5cm}

\parshape=1
3cm 10cm \noindent {\small \quad\quad\quad \quad\quad\quad\quad
\quad\quad\quad {\bf }\newline In this paper we prove the
birational rigidity of Fano-Mori fibre spaces $\pi\colon V\to S$,
every fibre of which is a Fano complete intersection of index 1
and codimension $k\geqslant 3$ in the projective space ${\mathbb
P}^{M+k}$ for $M$ sufficiently high, satisfying certain natural
conditions of general position, in the assumption that the fibre
space $V\slash S$ is sufficiently twisted over the base. The
dimension of the base $S$ is bounded from above by a constant,
depending only on the dimension $M$ of the fibre (as the dimension
of the fibre $M$ grows, this constant grows as $\frac12 M^2$).

Bibliography: 28 items.}


\vspace{0.3cm}

\section*{Introduction}

{\bf 0.1. Fano complete intersections.} In the present paper we
study the birational geometry of algebraic varieties, fibred into
Fano complete intersections of codimension $k\geqslant 3$
(fibrations into Fano hypersurfaces were studied in
\cite{Pukh15a}, into Fano complete intersections of codimension 2
in \cite{Pukh2022a}). We start with a description of fibres of
these fibre spaces. Let us fix an integer $k\geqslant 3$ and set
$$
\varepsilon(k)=\mathop{\rm min} \left\{a\in {\mathbb Z}\,\left|\,
a\geqslant 1, \left(1+\frac{1}{k}\right.\right)^a\geqslant
2\right\}.
$$
Now let us fix $M\in{\mathbb Z}$, satisfying the inequality
\begin{equation}\label{14.11.22.1}
M\geqslant 10 k^2+8k+2\varepsilon(k)+3.
\end{equation}
The right hand side of that inequality denote by the symbol
$\rho(k)$. Let
$$
\underline{d}=(d_1,\dots, d_k)
$$
be an ordered tuple of integers,
$$
2\leqslant d_1\leqslant d_2\leqslant\dots\leqslant d_k,
$$
satisfying the equality
$$
d_1+\dots +d_k=M+k.
$$
Fano varieties, considered in this paper, are complete
intersections of type $\underline{d}$ in the complex projective
space ${\mathbb P}^{M+k}$. More precisely, let the symbol ${\cal
P}_{a,N}$ stand for the space of homogeneous polynomials of degree
$a\in {\mathbb Z}_+$ in $N\geqslant 1$ variables. Set
$$
{\cal P}=\prod^k_{i=1} {\cal P}_{d_i,M+k+1}
$$
to be the space of all tuples
$$
\underline{f}=(f_1,\dots, f_k)
$$
of homogeneous polynomials of degree $d_1$, \dots, $d_k$ on
${\mathbb P}^{M+k}$. If for $\underline{f}\in {\cal P}$ the scheme
of common zeros of the polynomials $f_1,\dots, f_k$ is an
irreducible reduced factorial variety $F=F(\underline{f})$ of
dimension $M$ with terminal singularities, then $F$ is a primitive
Fano variety:
$$
\mathop{\rm Pic} F={\mathbb Z} H_F,\quad K_F=-H_F,
$$
where $H_F$ is the class of a hyperplane section of the variety
$F$ (the Lefschetz theorem). Assuming that this is the case, let
us give the following definition.

{\bf Definition 0.1.} The variety $F$ is {\it divisorially
canonical}, if for every effective divisor $D\sim nH_F$ the pair
$(F,\frac{1}{n}D)$ is canonical, that is, for every exceptional
divisor $E$ over $F$ the inequality
$$
\mathop{\rm ord}\nolimits_E D\leqslant n\cdot a(E)
$$
holds, where $a(E)$ is the discrepancy of $E$ with respect to $F$.

Below is the first main result of the present paper.

{\bf Theorem 0.1.} {\it There exist a Zariski open subset ${\cal
F}\subset{\cal P}$, such that for every tuple
$\underline{f}\in{\cal F}$ the scheme of common zeros of the tuple
$\underline{f}$ is an irreducible reduced factorial divisorially
canonical variety $F(\underline{f})$ of dimension $M$ with terminal
singularities, and the codimension of the complement ${\cal
P}\setminus{\cal F}$ satisfies the inequality}
$$
\mathop{\rm codim}(({\cal P}\setminus{\cal F})\subset{\cal P})
\geqslant M-k+5+{M-\rho(k)+2 \choose 2}.
$$

(Thus for a fixed $k$ and growing $M$ the codimension of the
complement ${\cal P}\setminus{\cal F}$ grows as $\frac12 M^2$.)

It is convenient to express the property of divisorial canonicity
in terms of the {\it global canonical threshold} of the variety
$F$.

Recall that for a Fano variety $X$ with the Picard number 1 and
terminal ${\mathbb Q}$-factorial singularities its global
canonical threshold $\mathop{\rm ct}(X)$ is the supremum of
$\lambda\in{\mathbb Q}_+$ such that for every effective divisor
$D\sim -nK_X$ (here $n\in{\mathbb Q}_+$) the pair
$\left(X,\frac{\lambda}{n}D\right)$ is canonical. Therefore,
Theorem 0.1 claims that for every $\underline{f}\in{\cal F}$ the
inequality $\mathop{\rm ct} (F(\underline{f}))\geqslant 1$ holds.

If in the definition of the global canonical threshold instead of
``for every effective divisor $D\sim -nK_X$'' we put ``for a
general divisor $D$ in any linear system $\Sigma\subset|-nK_X|$
with no fixed components'', we get the definition of the {\it
mobile canonical threshold} $\mathop{\rm mct}(X)$; obviously,
$\mathop{\rm mct}(X)\geqslant\mathop{\rm ct}(X)$. The inequality
$\mathop{\rm mct}(X)\geqslant 1$ is equivalent to the birational
superrigidity of the Fano variety $X$, see \cite{Ch05c}. If in the
definition of the global canonical threshold the property of the
pair $(X,\frac{\lambda}{n}D)$ to be canonical we replace by the
log canonicity of that pair, we get the definition of the {\it
global log canonical threshold} $\mathop{\rm lct}(X)$; again,
$\mathop{\rm lct}(X)\geqslant\mathop{\rm ct}(X)$.

For simplicity we write $F\in{\cal F}$ instead of
$F=F(\underline{f})$ for $\underline{f}\in{\cal F}$.\vspace{0.3cm}


{\bf 0.2. Fano-Mori fibre spaces.} By a {\it Fano-Mori fibre
space} we mean a surjective morphism of projective varieties
$$
\pi\colon V\to S,
$$
where $\dim V\geqslant 3 + \dim S$, the base $S$ is non-singular
and rationally connected, and the following conditions are
satisfied:

(FM1) every scheme fibre $F_s=\pi^{-1}(s)$, $s\in S$, is an
irreducible reduced factorial Fano variety with terminal
singularities and the Picard group $\mathop{\rm Pic} F_s\cong
{\mathbb Z}$,

(FM2) the variety $V$ itself is factorial and has at most terminal
singularities,

(FM3) the equality
$$
\mathop{\rm Pic} V={\mathbb Z} K_V\oplus\pi^* \mathop{\rm Pic}S
$$
holds.

So Fano-Mori fibre spaces are Mori fibre spaces with additional
very good properties.

{\bf Definition 0.2.} A Fano-Mori fibre space $\pi\colon V\to S$
is {\it stable with respect to fibre-wise birational
modifications}, if for every birational morphism $\sigma_S\colon
S^+\to S$, where $S^+$ is a non-singular projective variety, the
morphism
$$
\pi_+\colon V^+=V\mathop{\times}\nolimits_S S^+\to S^+
$$
is a Fano-Mori fibre space.

We will consider birational maps $\chi\colon V\dashrightarrow V'$,
where $V$ is the total space of a Fano-Mori fibre space and $V'$
is the total space of a fibre space $\pi'\colon V'\to S'$ which
belongs to one of the two classes:

(1) {\it rationally connected fibre spaces}, that is, $V'$ and
$S'$ are non-singular and the base $S'$ and a fibre of general
position $(\pi')^{-1}(s')$ are rationally connected,

(2) {\it Mori fibre spaces}, where $V'$ and $S'$ are projective
and the variety $V'$ has ${\mathbb Q}$-factorial terminal
singularities.

For a birational map $\chi\colon V\dashrightarrow V'$, where
$V'/S'$ is a rationally connected fibre space, we want to answer
the question: is it fibre-wise, that is, is there a rational
dominant map $\beta\colon S\dashrightarrow S'$, making the diagram
\begin{equation}\label{15.11.22.1}
\begin{array}{rcccl}
   & V & \stackrel{\chi}{\dashrightarrow} & V' & \\
\pi\!\!\!\!\! & \downarrow &   &   \downarrow & \!\!\!\!\!\pi' \\
   & S & \stackrel{\beta}{\dashrightarrow} & S'
\end{array}
\end{equation}
a commutative one, that is, $\pi'\circ\chi=\beta\circ\pi$?

For a birational map $\chi\colon V\dashrightarrow V'$, where
$V'/S'$ is a Mori fibre space with the additional properties (2)
(only such Mori fibre spaces are considered in this paper), we
want to answer the question: is there a {\it birational} map
$\beta\colon S\dashrightarrow S'$, for which the diagram
(\ref{15.11.22.1}) is commutative? If the answer to this question
is always affirmative (that is, it is affirmative for every fibre
space from the class (2)), then the fibre space $V/S$ is {\it
birationally rigid}.

Now let us state the second main result of the present paper.

{\bf Theorem 0.2.} {\it Assume that a Fano-Mori fibre space
$\pi\colon V\to S$ is stable with respect to fibre-wise birational
modifications, and moreover,

{\rm (i)} for every point $s\in S$ the fibre $F_s$ satisfies the
inequalities $\mathop{\rm lct} (F_s)\geqslant 1$ and $\mathop{\rm
mct} (F_s)\geqslant 1$,

{\rm (ii)} (the $K$-condition) every mobile (that is, with no
fixed components) linear system on $V$ is a subsystem of a
complete linear system $|-nK_V+\pi^* Y|$, where $Y$ is a
pseudoeffective class on $S$,

{\rm (iii)} for every family ${\overline{\cal C}}$ of irreducible
curves on $S$, sweeping out a dense subset of the base $S$, and
$\overline{C}\in{\overline{\cal C}}$, no positive multiple of the
class
$$
-(K_V\cdot \pi^{-1}(\overline{C}))-F\in A^{\dim S} V,
$$
where $A^iV$ is the numerical Chow group of classes of cycles of
codimension $i$ on $V$ and $F$ --- the class of a fibre of the
projection $\pi$, is represented by an effective cycle on $V$.

Then for every rationally connected fibre space $V'\slash S'$
every birational map $\chi\colon V\dashrightarrow V'$ (if such
maps exist) is fibre-wise, and the fibre space $V\slash S$ itself
is birationally rigid.}

By what was said in Subsection 0.1, the assumption (i) can be
replaced by the single inequality $\mathop{\rm ct}(F_s)\geqslant
1$ for every $s\in S$, that is, it is sufficient to assume that
every fibre of the fibre space $V\slash S$ is a divisorially
canonical variety.

As we will see from the proof of Theorem 0.2, instead of the
conditions (ii) and (iii) it is sufficient to require that for
every family $\overline{\cal C}$ of irreducible curves on $S$,
sweeping out a dense subset, and $\overline{C}\in \overline{\cal
C}$ the class
$$
-N(K_V\cdot \pi^{-1}(\overline{C}))-F\in A^{\dim S} V
$$
is not represented by an effective cycle on $V$ for any
$N\geqslant 1$. The last condition is especially easy to verify:
it is enough to have a numerically effective $\pi$-ample class
$H_V$ on $V$, satisfying the inequality
\begin{equation}\label{18.11.22.1}
\left(K_V\cdot \pi^{-1}(\overline{C})\cdot H_V^{\dim V-\dim
S}\right)\geqslant 0
\end{equation}
for every dense family $\overline{\cal
C}\ni\overline{C}$.\vspace{0.3cm}


{\bf 0.3. An explicit construction of a fibre space.} Now let us
construct a large class of Fano-Mori fibre spaces, satisfying the
conditions of Theorem 0.2. Let $S$ be a non-singular projective
rationally connected positive-dimensional variety and $\pi_X\colon
X\to S$ a locally trivial fibration with the fibre ${\mathbb
P}^{M+k}$, where $k$ and $M$ are the same as in Subsection 0.1. We
say that the subvariety $V\subset X$ of codimension $k$ is a {\it
fibration into complete intersections of type} $\underline{d}$, if
the base $S$ can be covered by Zariski open subsets $U$, over
which the fibration $\pi_X$ is trivial, $\pi^{-1}_X(U)\cong
U\times{\mathbb P}^{M+k}$, and for every $U$ there is a regular
map
$$
\Phi_{U}\colon U\to{\cal P},
$$
such that $V\cap\,\pi^{-1}_X(U)$ in the sense of the
above-mentioned trivialization is the scheme of common zeros of a
tuple
$$
\underline{f}(s)=\Phi_{U}(s)=(f_1(x_*,s),\dots,f_k(x_*,s)),
$$
where $x_*$ are homogeneous coordinates on ${\mathbb P}^{M+k}$ and
$s$ runs through $U$.

Below (in \S 1) it will be clear, that the open subset ${\cal F}$
from Theorem 0.1 is invariant under the action of the group
$\mathop{\rm Aut}{\mathbb P}^{M+k}$. For that reason, the
following definition makes sense.

{\bf Definition 0.3.} A fibration $V\subset X$ into complete
intersections of type $\underline{d}$ is a ${\cal F}$-{\it
fibration}, if for any trivialization of the bundle $\pi_X$ over
an open set $U\subset S$ we have $\Phi_{U}(U)\subset{\cal F}$.

Obviously, if the inequality
\begin{equation}\label{18.11.22.2}
\mathop{\rm dim} S\leqslant M-k+4+{M-\rho(k)+2 \choose 2}
\end{equation}
holds, then we may assume that $V$ is a ${\cal F}$-fibration. Set
$\pi=\pi_X|_V$. Now from Theorems 0.1 and 0.2 it is easy to obtain
the third main result of the present paper.

{\bf Theorem 0.3.} {\it Any ${\cal F}$-fibration $\pi\colon V\to
S$ constructed above is a Fano-Mori fibre space. If the conditions
(ii) and (iii) of Theorem 0.2 hold, then for every rationally
connected fibre space $V'/S'$ every birational map $\chi\colon
V\dashrightarrow V'$ is fibre-wise, and the fibre space $V\slash
S$ itself is birationally rigid.}

{\bf Example 0.1.} Let $H_X$ be a numerically effective divisorial
class on $X$, the restriction of which onto the fibre
$\pi_X^{-1}(s)\cong {\mathbb P}^{M+k}$ is the class of a
hyperplane. Let $\Delta_1$,\dots, $\Delta_k$ be very ample classes
on the base $S$. Let us construct a ${\cal F}$-fibration $V/S$ as
a complete intersection of $k$ general divisors
$$
V=G_1\cap \dots \cap G_k,
$$
where $G_i\in|d_i H_X+\pi_X^*\Delta_i|$. Let us find out, when
$V/S$ satisfies the conditions (ii) and (iii) of Theorem 0.2.
Write
$$
K_X=-(M+k+1) H_X+\pi^*_X \Delta_X,
$$
then we get
$$
K_V=\left.\left(-H_X+\pi^*_X\left(\Delta_X+\sum^k_{i=1}\Delta_i\right)
\right)\right|_V.
$$
It is easy to check that the inequality (\ref{18.11.22.1}) in this
case takes the form of the estimate
$$
\left(\left(\Delta_X+\sum^k_{i=1}\left(1-\frac{1}{d_i}\right)\Delta_i\right)\cdot
\overline{C}\right)\geqslant \left(H_X^{M+k+1}\cdot
\pi_X^{-1}(\overline{C})\right),
$$
where for the class $H_V$ we took $H_X|_V$. This inequality must
be satisfied for every dense family $\overline{\cal
C}\ni\overline{C}$.

Let us consider a very particular case, when $X={\mathbb
P}^m\times{\mathbb P}^{M+k}$ and $G_i$ are divisors of bi-degree
$(m_i,d_i)$, $i=1,\dots,k$. Taking for $H_X$ the pull back on $X$
of the class of a hyperplane in ${\mathbb P}^{M+k}$, we get that
the last inequality is equivalent to the numerical inequality
\begin{equation}\label{18.11.22.3}
\sum^k_{i=1}\left(1-\frac{1}{d_i}\right)m_i\geqslant m+1.
\end{equation}
If it is satisfied and the dimension $m=\mathop{\rm dim} S$
satisfies the inequality (\ref{18.11.22.2}), then the intersection
$V=G_1\cap\dots\cap G_k$ of general (in the sense of Zariski
topology) divisors of bi-degree $(m_1,d_1),\dots,(m_k,d_k)$,
fibred over $S={\mathbb P}^m$, is a birationally rigid Fano-Mori
fibre space and every birational map of $V$ onto the total space
of a rationally connected fibre space is fibre-wise. The
inequality (\ref{18.11.22.3}) shows that this claim holds for
almost all tuples $(m_1,\dots,m_k)\in{\mathbb Z}^k_+$ (except for
finitely many of them), that is, for almost all families of
Fano-Mori fibre spaces, obtained by means of this construction.
Note that the condition (\ref{18.11.22.3}) is lose to a criterial
one: if
$$
m_1+\dots+m_k\leqslant m,
$$
then the projection of $V$ onto ${\mathbb P}^{M+k}$ defines on $V$
a structure of a Fano-Mori fibre space (and a rationally connected
fibre space), which is ``transversal'' to the original structure
$\pi\colon V\to S$ (and is not fibre-wise), so that in this case
$V/S$ is not birationally rigid.\vspace{0.3cm}


{\bf 0.4. The structure of the paper.} The paper is organized in
the following way. In \S 1 we produce the explicit local
conditions defining the open subset ${\cal F}\subset{\cal P}$. The
proof of divisorial canonicity of a variety $F\in{\cal F}$ (that
is, of the inequality $\mathop{\rm ct}(F)\geqslant 1$) is reduced
in \S 1 to a number of technical facts that will be shown in the
subsequent sections (\S\S 3-7). In \S 2 we show Theorem 0.2.

The proof of Theorem 0.1 consists of several pieces. The fact that
the local conditions for the singularities that a variety
$F\in{\cal F}$ can have ({\it multi-quadratic singularities}, see
Subsection 1.2), guarantee that the variety $F$ is factorial and
its singularities are terminal, is proven in \S 4, where we give a
general definition of multi-quadratic singularities and study
their properties. The estimate for the codimension of the
complement ${\cal P}\setminus {\cal F}$ (which is very important
for constructing families of Fano-Mori fibre spaces, satisfying
the assumptions of Theorem 0.2) is shown in \S 8. However, the
main (and the hardest) part of the proof of Theorem 0.1 is to show
that a variety $F\in{\cal F}$ is divisorially canonical. We assume
that for some effective divisor $D\sim nH_F$ the pair
$(F,\frac{1}{n}D)$ is not canonical, that is, for some exceptional
divisor $E$ over $F$ the inequality
$$
\mathop{\rm ord}\nolimits_E D>n\cdot a(E)
$$
holds. Now we have to show that this assumption leads to a
contradiction. In Subsections 1.3-1.6 it is shown how (using the
inequalities for the multiplicity of subvarieties of the variety
$F$ at a given point, proven in \S 7) to obtain a contradiction in
the case when a point of general position $o\in B$, where $B$ is
the centre of $E$ on $F$, either is non-singular on $F$, or is a
quadratic singularity. The hardest task is to obtain a
contradiction when the point $o$ is a multi-quadratic singularity
of the variety $F$. A plan of solving this problem is given in
Subsection 1.7, where we introduce the concept of a {\it working
triple} and describe the procedure of construction a sequence of
subvarieties of the variety $F$, in which each subvariety is a
hyperplane section of the previous one and the last subvariety
delivers the desired contradiction.

This program is realized in \S 3, where we study the properties of
working triples, however a number of key technical facts is only
stated there --- their proof is put off for a greater clarity of
exposition. These key facts are shown in \S\S 5,6 (and the proof
makes use of the facts on linear subspaces on complete
intersections of quadrics, proven in Subsection 4.5).

Finally, in \S 7 we prove the estimates for the multiplicities of
certain subvarieties of the variety $F$ at given points in terms
of the degrees of these subvarieties in ${\mathbb P}^{M+k}$. Here
we use the well known technique of hypertangent divisors. For the
purposes of our proof of Theorem 0.1 we have to somewhat modify
this technique.\vspace{0.3cm}


{\bf 0.5. General remarks.} The birational rigidity of Fano-Mori
fibre spaces over a positive-dimensional base was one of the most
important topics in birational geometry in the past 40 years. For
its history and place in the context of the modern birational
geometry of rationally connected varieties, see \cite[Subsection
0.4]{Pukh2022a}. Here we just mention a few recent papers in the
areas that are close to the direction, to which the present paper
belongs.

These areas are: the birational rigidity, explicit birational
geometry of Mori fibre spaces (including the studies of their
groups of birational automorphisms and, wider, Sarkisov links),
the rationality problem, computing and estimating the global
canonical thresholds and, related to these problems, the theory of
$K$-stability.

In the papers \cite{KrylovOkadaetal22,AbbanKrylov22,Krylov18}
important results on the birational rigidity and rigidity-type
results for fibrations over ${\mathbb P}^1$ were obtained. The
paper \cite{Stibitz21} links the Sarkisov program with the problem
of estimating the canonical threshold of certain divisors on Fano
varieties. The papers \cite{AbbanOkada18,KrylovOkada20} prove the
stable non-rationality of very general conic bundles and
fibrations into del Pezzo surfaces, respectively, over a
higher-dimensional base. The problem of stable rationality for
hypersurfaces of various bi-degrees in the products of projective
spaces (see Example 0.1 above) is considered in
\cite{NicaiseOttem22}. The theory of $K$-stability, which is on
the border of birational geometry, is investigated in many papers
(especially in the recent past), in particular, see
\cite{StibitzZhuang19,Zhuang21,CheltsovPark22,CheltsovDenisovaetal22};
we mentioned the papers that are the closest to the birational
rigidity-type problems. Finally, there was a lot of development
recently in the direction of applying the theory of Sarkisov links
and relations between them to the study of the groups of
birational automorphisms of such varieties that have a very large
this group, see, for instance,
\cite{BlancYasinsky20,BlancLamyZ21}.

Getting back to the topic of this paper, we note that its
immediate predecessor is \cite{Pukh2022a}, however, that paper
investigates the non canonical singularities, the centre of which
is contained in the set of bi-quadratic points of the variety
(from the technical viewpoint, this is the hardest part of the
proof of divisorial canonicity), using the secant varieties of
subvarieties of codimension 2 on an intersection of two quadrics.
It is not possible to apply this approach to subvarieties of
higher codimension on an intersection of $k\geqslant 3$ quadrics,
and the present paper is based on a completely different
construction (which applies to the bi-quadratic singularities,
considered in \cite{Pukh2022a}, as well).

The author is grateful to the members of Divisions of Algebraic Geometry and
Algebra at Steklov Institute of Mathematics for the interest to
his work, and also to the colleagues in Algebraic Geometry
research group at the University of Liverpool for general support.


\section{Fano complete intersections}

In this section we describe the local conditions defining the open
subset ${\cal F}\subset{\cal P}$ (Subsections 1.2 and 1.4). For a
complete intersection $F\in{\cal F}$ the proof of its divisorial
canonicity is reduced to a number of technical claims, which will
be shown later. A more detailed plan of the proof of Theorem 0.1
is given in Subsection 1.1.\vspace{0.3cm}

{\bf 1.1. A plan of the proof of Theorem 0.1.} In order to prove
Theorem 0.1, one has to give an explicit definition of the open
set ${\cal F}\subset{\cal P}$. This definition consists of two
groups of conditions, which should be satisfied by the polynomials
$f_1,\dots, f_k$ at every point $o\in {\mathbb P}^{M+k}$ at which
they all vanish. The first group of conditions is about the
singularities of the complete intersection $F(\underline{f})$:
they can be quadratic or multi-quadratic of a rank bounded from
below. The corresponding definitions and facts are given in
Subsection 1.2. Assuming that the conditions of the first group
are satisfied, we get that the scheme of common zeros of the
polynomials $f_1,\dots, f_k$ is an irreducible reduced factorial
variety $F=F(\underline{f})\subset {\mathbb P}^{M+k}$ with
terminal singularities, and so $\mathop{\rm Pic} F={\mathbb Z}
H_F$ and $K_F=-H_F$, so that the question, is it divisorially
canonical, makes sense.

Assuming that $F$ is not divisorially canonical, let us fix an
effective divisor $D_F\sim n(D_F)H_F$, where $n(D_F)\geqslant 1$,
such that the pair
$$
\left(F,\frac{1}{n(D_F)}D_F\right)
$$
is not canonical, that is, there is an exceptional divisor $E$
over $F$, satisfying the Noether-Fano inequality:
$$
\mathop{\rm ord}\nolimits_E D_F>n(D_F)\, a(E).
$$
We have to show that the existence of such a divisor leads to a
contradiction. Let $B\subset F$ be the centre of the exceptional
divisor $E$ on $F$. The information about the singularities of the
varieties $F$ makes it possible to easily exclude the option when
$\mathop{\rm codim}(B\subset F)=2$. This is done in Subsection
1.3.

After that in Subsection 1.4 we produce the second group of local
conditions for the tuple of polynomials $\underline{f}\in{\cal
F}$: now they are the regularity conditions. Assuming that they
are satisfied at every point $o\in F$, we exclude the option
$B\not\subset\mathop{\rm Sing}F$ in Subsection 1.5, and in
Subsection 1.6 the option that the point $o\in B$ of general
position is a quadratic singularity of $F$. In Subsection 1.7 we
describe the procedure of excluding the multi-quadratic case, when
the point $o\in B$ of general position is a multi-quadratic
singularity of the type $2^l$, $l\in\{2,\dots,k\}$. This is the
hardest part of the work, which is completed in the subsequent
sections.\vspace{0.3cm}


{\bf 1.2. Multi-quadratic singularities.} Let $o\in{\mathbb
P}^{M+k}$ be a point, at which $f_1,\dots,f_k$ all vanish. Let us
consider a system of affine coordinates $z_*=(z_1,\dots,z_{M+k})$
with the origin at the point $o$ on an affine chart ${\mathbb
A}^{M+k}\subset{\mathbb P}^{M+k}$, containing that point. Write
down
$$
\begin{array}{ccccc}
f_1=f_{1,1}+f_{1,2}+ & \dots & + f_{1,d_1}, & &\\
f_2=f_{2,1}+f_{2,2}+ & \dots &  & + f_{2,d_2}, & \\
&                      \dots & & &\\
f_k=f_{k,1}+f_{k,2}+ & \dots & & & +f_{k,d_k},
\end{array}
$$
where we use the same symbols $f_i$ for the non-homogeneous
polynomials in $z_*$, corresponding to the original polynomials
$f_i$, and $f_{i,a}$ is a homogeneous polynomial of degree $a$ in
$z_*$. Obviously, if the linear forms $f_{1,1},\dots,f_{k,1}$ are
linearly independent, then in a neighborhood of the point $o$ the
scheme of common zeros of the polynomials $f_1,\dots,f_k$ is a
non-singular complete intersection of codimension $k$. In order to
give the definition of a multi-quadratic singularity, we will need
the concept of the rank of a tuple of quadratic forms.

{\bf Definition 1.1.} (\cite{Pukh2022a}) {\it The rank of the
tuple of quadratic forms} $q_1,\dots,q_l$ in $N$ variables is the
number
$$
\mathop{\rm rk}(q_1,\dots,q_l)=\mathop{\rm min}\{\mathop{\rm
rk}(\lambda_1 q_1+\dots+\lambda_l
q_l)\,|\,(\lambda_1,\dots,\lambda_l)\neq (0,\dots,0)\}.
$$
Obviously, $\mathop{\rm rk}(q_1,\dots,q_l)\leqslant N$. For that
reason, in the sequel the inequality $\mathop{\rm
rk}(q_*)\geqslant a$ means implicitly that the forms $q_i$ depend
on a sufficient $(\geqslant a)$ number of variables.

Take $l\in\{1,2,\dots,k\}$.

{\bf Definition 1.2.} The tuple $\underline{f}$ has at the point
$o$ a {\it multi-quadratic singularity of type} $2^l$ of rank $a$,
if the following conditions are satisfied:

\begin{itemize}

\item $\mathop{\rm dim}\langle f_{1,1},\dots,f_{1,k}\rangle=k-l$
(and in order to simplify the notations we assume that the forms
$$
f_{l+1},\dots, f_k
$$
are linearly independent),

\item the rank of the tuple of quadratic forms
$$
f^*_{i,2}=f_{i,2}-\sum^k_{j=l+1}\lambda_{i,j}f_{j,2},
$$
$i=1,\dots,l$, where $\lambda_{i,j}\in{\mathbb C}$ are defined by
the equalities
$$
f_{i,1}=\sum^k_{j=l+1}\lambda_{i,j}f_{j,1},
$$
is equal to the number $a$.
\end{itemize}

Now the first condition, defining the subset ${\cal F}\subset{\cal
P}$, is stated in the following way.

(MQ1) For every point $o\in{\mathbb P}^{M+k}$, such that
$$
f_1(o)=\dots=f_k(o)=0
$$
either the linear forms $f_{1,1},\dots,f_{k,1}$ are linearly
independent, or $\underline{f}$ has at the point $o$ a
multi-quadratic singularity of type $2^l$, where
$l\in\{1,2,\dots,k\}$, of rank
$$
\geqslant 2l+4k+2\varepsilon(k)-1.
$$

{\bf Theorem 1.1.} {\it Assume that $\underline{f}$ satisfies the
condition (MQ1). Then the scheme of common zeros of the
polynomials $f_1,\dots,f_k$ is an irreducible reduced factorial
variety F=F(\underline{f}) --- a complete intersection of
codimension $k$ with terminal singularities, and, moreover,}
$$
\mathop{\rm codim}(\mathop{\rm Sing}F\subset F)\geqslant
4k+2\varepsilon(k).
$$

{\bf Proof} is given in \S 4 (Subsections 4.1-4.3).

Assume that $\underline{f}$ satisfies the condition (MQ1). For a
point $o\in F=F(\underline{f})$ the symbol $T_oF$ stands for the
subspace $\{f_{1,1}=\dots=f_{k,1}=0\}\subset{\mathbb C}^{M+k}$.
For the proof of Theorem 0.1 we will need one more property of the
tuple $\underline{f}$, which we include in the definition of the
subset ${\cal F}$.

(MQ2) For any point $o\in F$, which is a multi-quadratic of type
$2^l$, where $l\geqslant 2$, the rank of the tuple of quadratic
forms
$$
f_{1,2}|_{T_oF},\dots,f_{k,2}|_{T_oF}
$$
is at least $10k^2+8k+2\varepsilon(k)+5$.

The condition (MQ2) for multi-quadratic points of type $2^l$ with
$l\geqslant 2$ implies the condition (MQ1), because the rank of a
quadratic form, restricted to a hyperplane, drops at most by 2,
however for the convenience of references we state the conditions
(MQ1) and (MQ2) independently of each other. These conditions are
used in the proof of Theorem 0.1 in different ways.

So every tuple $\underline{f}\in{\cal F}$ satisfies (MQ1) and
(MQ2).\vspace{0.3cm}


{\bf 1.3. Subvarieties of codimension 2.} Following the plan,
given in Subsection 1.1, let us fix an effective divisor $D_F\sim
n(D_F)H_F$, $n(D_F)\geqslant 1$, such that the pair
$(F,\frac{1}{n(D_F)}D_F)$ is not canonical. By the symbol
$$
\mathop{\rm CS}\left(F,\frac{1}{n(D_F)}D_F\right)
$$
we denote the union of the centres on $F$ of all exceptional
divisors over $F$, satisfying the Noether-Fano inequality (that is
to say, of all non-canonical singularities of that pair). This is
a closed subset of $F$. Let $B$ be an irreducible component of
maximal dimension of that set.

{\bf Proposition 1.1.} {\it The following inequality holds:}
$\mathop{\rm codim}(B\subset F)\geqslant 3$.

{\bf Proof.} Assume the converse: $\mathop{\rm codim}(B\subset
F)=2$. Then $B\not\subset\mathop{\rm Sing}F$. Moreover, let
$P\subset{\mathbb P}^{M+k}$ be a general linear subspace of
codimension $2k+2$. Theorem 1.1 implies that $P\cap\mathop{\rm
Sing}F=\emptyset$, so that $F\cap P$ is a non-singular complete
intersection of type $\underline{d}$ in $P\cong{\mathbb
P}^{2k+2}$. Furthermore, the pair
$$
\left(F\cap P,\frac{1}{n(D_F)}D_F|_{F\cap P}\right)
$$
is not canonical, and the irreducible subvariety $B\cap P$ is an
irreducible component of maximal dimension of the set
$$
\mathop{\rm CS}\left(F\cap P,\frac{1}{n(D_F)}D_F|_{F\cap
P}\right),
$$
so that (as $F\cap P$ is non-singular)
$$
\mathop{\rm mult}\nolimits_{B\cap P}D_F|_{F\cap P}>n(D_F).
$$
However, $D_F|_{F\cap P}\sim n(D_F)H_{F\cap P}$ (where $H_{F\cap
P}$ is the class of a hyperplane section of $F\cap P$), so that by
\cite[Proposition 3.6]{Pukh06b} or \cite{Suzuki15} we get a
contradiction, proving the proposition. Q.E.D.\vspace{0.3cm}


{\bf 1.4. Regularity conditions.} In order to continue the proof
of Theorem 0.1, we need a second group of conditions defining the
set ${\cal F}$. Let $o\in F$ be a point. We use the notations of
Subsection 1.2. By the symbol $T_oF$ we denote the {\it linear}
tangent space
$$
\{f_{1,1}=\dots=f_{k,1}=0\}\subset{\mathbb C}^{M+k},
$$
and by the symbol ${\mathbb P}(T_oF)$ its projectivization. Let
${\cal S}=(h_1,\dots,h_M)$ be the sequence of homogeneous
polynomials
$$
f_{i,j}|_{{\mathbb P}(T_oF)},
$$
where $j\geqslant 2$, placed in the lexicographic order:
$(i_1,j_1)$ precedes $(i_2,j_2)$, if $j_1<j_2$ or $j_1=j_2$, but
$i_1<i_2$. By the symbol ${\cal S}[-m]$ denote the sequence ${\cal
S}$ with the last $m$ members removed. Finally, the symbol ${\cal
S}[-m]|_{\Pi}$ stands for the restriction of that sequence (that
is, the restriction of each its member) onto a linear subspace
$\Pi\subset{\mathbb P}(T_oF)$. The regularity conditions depend on
the type of the singularity $o\in F$.

First, let the point $o\in F$ be non-singular, so that ${\mathbb
P}(T_oF)\cong{\mathbb P}^{M-1}$. In that case the regularity
condition is stated in the following way.

(R1) The sequence
$$
{\cal S}[-(k+\varepsilon(k)+3)]|_{\Pi}
$$
is regular for every subspace $\Pi\subset{\mathbb P}(T_oF)$ of
codimension $k+\varepsilon(k)-1$.

The condition (R1) is assumed for every non-singular point $o\in
F$. It implies the following key fact.

{\bf Theorem 1.2.} {\it Let $P\subset{\mathbb P}^{M+k}$ be an
arbitrary linear subspace of codimension $\leqslant k+
\varepsilon(k)-1$. Then for every non-singular point $o\in F\cap
P$ and every prime divisor $Y\sim n(Y)H_{F\cap P}$ on $F\cap P$
the inequality}
$$
\mathop{\rm mult}\nolimits_oY\leqslant2n(Y)
$$
{\it holds.}

{\bf Proof} is given in \S 7 (Subsections 7.1, 7.2).

Now let $o\in F$ be a quadratic singularity (this case corresponds
to the value $l=1$ in Definition 1.2). Here ${\mathbb
P}(T_oF)\cong{\mathbb P}^M$. In this case the regularity condition
is stated as follows.

(R2) The sequence
$$
{\cal S}[-4]|_{\Pi}
$$
is regular for every hyperplane $\Pi\subset{\mathbb P}(T_oF)$.

The condition (R2) is assumed for every quadratic singular point
$o\in F$ and implies the following key fact.

{\bf Theorem 1.3.} {\it Let $o\in F$ be a quadratic singularity
and $W\ni o$ the section of $F$ by a hyperplane that is not
tangent to $F$ at the point $o$, and $Y\sim n(Y)H_W$ a prime
divisor on $W$. Then the following inequality holds:}
$$
\mathop{\rm mult}\nolimits_oY\leqslant 4n(Y).
$$

{\bf Proof} is given in \S 7 (Subsection 7.3).

(The symbol $H_W$ stands for the class of a hyperplane section of
the variety $W$; the linear form, defining the hyperplane that
cuts out $W$, is not a linear combination of the forms
$f_{1,1},\dots,f_{k,1}$.)

Now let $o\in F$ be a multi-quadratic point of type $2^l$, where
$l\in\{2,\dots,k\}$. Here we will need two regularity conditions.
In the first of them the symbol $T_oF$ means the projective
closure of the linear subspace
$$
\{f_{1,1}=\dots=f_{k,1}=0\}\subset{\mathbb C}^{M+k}
$$
in ${\mathbb P}^{M+k}$.

(R3.1) For every subspace $P\subset T_oF$ of codimension
$\varepsilon(k)$, containing the point $o$, the scheme of common
zeros of the polynomials
$$
f_1|_P,\dots,f_k|_P,\quad f_{i,2}|_P\quad\mbox{for all}\quad i:
d_i\geqslant 3,
$$
is an irreducible reduced subvariety of codimension
$k+k_{\geqslant 3}$ in $P$, where
$$
k_{\geqslant 3}=\sharp\{i=1,\dots,k\,|\, d_i\geqslant 3\}.
$$
Note that in the condition (R3.1) the homogeneous polynomials
$f_{i,2}$ in the {\it affine} coordinates $z_*$ are considered as
quadratic forms in {\it homogeneous} coordinates on ${\mathbb
P}^{M+k}$.

In the second regularity condition for multi-quadratic points the
symbol $T_oF$ means a linear subspace in ${\mathbb C}^{M+k}$.

(R3.2) For every linear subspace $\Pi\subset{\mathbb P}(T_oF)$ of
codimension $\varepsilon(k)$ the sequence
$$
{\cal S}[-m^*]|_{\Pi}
$$
is regular, where $m^*=\mathop{\rm max}\{\varepsilon(k)+4-l,0\}$.

The conditions (R3.1) and (R3.2) are assumed for every
multi-quadratic singular point $o\in F$. They imply the following
key inequality. In Theorem 1.4, stated below, the symbol $T_oF$
stands for the projective closure of the embedded tangent space,
that is, a linear subspace in ${\mathbb P}^{M+k}$, containing the
point $o$.

{\bf Theorem 1.4.} {\it Let $P\subset T_oF$ be an arbitrary linear
subspace of codimension $\leqslant\varepsilon(k)$ and $Y\ni o$ a
prime divisor on $F\cap P$, $Y\sim n(Y)H_{F\cap P}$. Then the
following inequality holds:}
$$
\mathop{\rm mult}\nolimits_o Y\leqslant\frac32\cdot 2^kn(Y).
$$

{\bf Proof} is given in \S 7 (Subsections 7.4, 7.5).

(The symbol $H_{F\cap P}$ stands for the class of a hyperplane
section of the variety $F\cap P$; we will show below, see \S 4,
that $F\cap P$ is an irreducible factorial complete intersection.)

Summing up, let us give a complete definition of the subset ${\cal
F}\subset{\cal P}$: it consists of the tuples $\underline{f}$,
satisfying the conditions (MQ1,2), the condition (R1) at every
non-singular point $o\in F(\underline{f})$, the condition (R2) at
every quadratic point $o\in F(\underline{f})$ and the conditions
(R3.1,2) at every multi-quadratic point $o\in F(\underline{f})$.

The inequality for the codimension of the complement ${\cal
P}\setminus{\cal F}$, given in Theorem 0.1, is shown in \S
8.\vspace{0.3cm}


{\bf 1.5. Exclusion of the non-singular case.} We carry on with
the proof of divisorial canonicity of the variety $F\in{\cal F}$.
In the notations of Subsection 1.3 assume that the point of
general position $o\in B$ is a non-singular point of $F$. We know
(Proposition 1.1), that $\mathop{\rm codim}(B\subset F)\geqslant
3$. Consider a general subspace $P\ni o$ of dimension $k+3$. Then
$F\cap P$ is a non-singular three-dimensional variety and the
point $o$ is a connected component of the set
$$
\mathop{\rm CS}\left(F\cap P,\frac{1}{n(D_F)}D_F|_{F\cap P}\right)
$$
(if $\mathop{\rm codim}(B\subset F)\geqslant 4$, then $\mathop{\rm
CS}$ can be replaced, by inversion of adjunction, by $\mathop{\rm
LCS}$), that is, outside the point $o$ in a neighborhood of that
point the pair
\begin{equation}\label{23.11.22.1}
\left(F\cap P,\frac{1}{n(D_F)}D_F|_{F\cap P}\right)
\end{equation}
is canonical. It is well known (see \cite[Proposition 3]{Pukh05}
or \cite[Chapter 7, Proposition 2.3]{Pukh13a}), it follows from
here that either the inequality
$$
\mathop{\rm mult}\nolimits_oD_F>2n(D_F)
$$
holds, or on the exceptional divisor $E\cong{\mathbb P}^{M-1}$ of
the blow up $F^+\to F$ of the point $o$ there is a hyperplane
$\Theta\subset E$ (uniquely determined by the pair
(\ref{23.11.22.1})), such that the inequality
$$
\mathop{\rm mult}\nolimits_oD_F+\mathop{\rm
mult}\nolimits_{\Theta}D^+_F> 2n(D_F)
$$
holds, where $D^+_F$ is the strict transform of $D_F$ on $F^+$.

The first option is impossible as it contradicts Theorem 1.2. In
the second case denote by the symbol $|H-\Theta|$ the projectively
$k$-dimensional linear system of hyperplane sections of $F$, a
general element of which $W\ni o$ is non-singular at the point $o$
and satisfies the equality
$$
W^+\cap E=\Theta.
$$
The restriction $D_W=(D_F\circ W)$ is an effective divisor on $W$,
and $n(D_W)=n(D_F)$ and the inequality
$$
\mathop{\rm mult}\nolimits_oD_W\geqslant\mathop{\rm
mult}\nolimits_oD_F+\mathop{\rm
mult}\nolimits_{\Theta}D^+_F>2n(D_W)
$$
holds, which again contradicts Theorem 1.2. We have shown the
following fact.

{\bf Proposition 1.2.} {\it The subvariety $B$ is contained in the
singular locus of $F$:} $B\subset\mathop{\rm Sing}
F$.\vspace{0.3cm}


{\bf 1.6. Exclusion of the quadratic case.} Again let $o\in B$ be
a point of general position.

{\bf Proposition 1.3.} {\it The point $o$ is a multi-quadratic
singularity of type} $2^l$, $l\geqslant 2$.

{\bf Proof.} Assume the converse: the point $o$ is a quadratic
singularity of $F$. Let $P\ni o$ be a general $(k+3)$-dimensional
linear subspace in ${\mathbb P}^{M+k}$. By the condition (MQ1) and
Theorem 1.1 the intersection $F\cap P$ is a three-dimensional
variety with the unique singular point $o$, which is a
non-degenerate quadratic singularity. This intersection can be
constructed in two steps: first, we consider the inetrsection
$F\cap P'$ with a general linear subspace $P'\subset{\mathbb
P}^{M+k}$, $P'\ni o$, of dimension
$$
k+\mathop{\rm codim}(\mathop{\rm Sing}F\subset F)
$$
and after that the intersection with a general subspace $P\subset
P'$, $P\ni o$, of dimension $(k+3)$. Now we get: the pair
$$
\left(F\cap P,\frac{1}{n(D_F)}D_F|_{F\cap P}\right)
$$
is not log canonical, but canonical out side the point $o$. Let us
consider the blow up
$$
\varphi_P\colon P^+\to P
$$
of the point $o$ with the exceptional divisor ${\mathbb
E}_P\cong{\mathbb P}^{k+2}$ and let $(F\cap P)^+\subset P^+$ be
the strict transform of $F\cap P$ on $P^+$, so that $(F\cap
P)^+\to F\cap P$ is the blow up of the quadratic singularity $o$
with the exceptional divisor $E_P=(F\cap P)^+\cap{\mathbb E}_P$,
which is a non-singular two-dimensional quadric in the
three-dimensional subspace $\langle E_P\rangle\subset {\mathbb
E}_P$. Obviously, $a(E_P, F\cap P)=1$, so that, writing down
$$
D_P=D_F|_{F\cap P}\sim n(D_F)H_{F\cap P}
$$
and $D_P^+\sim n(D_F)H_{F\cap P}-\nu E_P$ (the strict transform of
$D_P$ on $(F\cap P)^+$), we obtain two options:
\begin{itemize}

\item either $\nu>2n(D_F)$, so that $E_P$ is a non log canonical
singularity of the pair $(F\cap P,\frac{1}{n(D_F)} D_P)$,

\item or $n(D_F)<\nu\leqslant 2n(D_F)$, and then the closed set
$$
\mathop{\rm LCS}\left(\left(F\cap P,\frac{1}{n(D_F)} D_P\right),
(F\cap P)^+\right)
$$
--- the union of the centres of all non log canonical singularities
of the original pair $(F\cap P,\frac{1}{n(D_F)} D_P)$ on $(F\cap
P)^+$ is a connected closed subset of the non-singular quadric
$E_P$, which can be either a (possibly reducible) connected curve
$C_P\subset E_P$, or a point $x_P\in E_P$.
\end{itemize}

(It is well known, see, for instance, \cite[Chapter 2, Proposition
3.7]{Pukh13a}, that the inequality $\nu\leqslant n(D_F$) is
impossible.) In the case $\nu>2n(D_F)$ we get
$$
\mathop{\rm mult}\nolimits_o D_P=\mathop{\rm mult}\nolimits_o D_F>
4n(D_F),
$$
which contradicts Theorem 1.3, so that this case is impossible.
Coming back to the original variety $F$, let us consider the blow
ups $\varphi_{\mathbb P}\colon({\mathbb P}^{M+k})^+\to{\mathbb
P}^{M+k}$ and $\varphi\colon F^+\to F$ of the point $o$, where
$F^+$ is identified with the strict transform of $F$ on $({\mathbb
P}^{M+k})^+$, with the exceptional divisors ${\mathbb E}$ and $E$,
respectively, so that $E=F^+\cap{\mathbb E}$ is a quadratic
hypersurface $E$ in the subspace $\langle E\rangle\subset{\mathbb
E}$ of codimension $(k+1)$. By the condition for the rank (MQ1)
the case of a point $x_P\in E_P$ is impossible: in that case the
quadric $E$ would contain a linear subspace of codimension 2 (with
respect to $E$), which can not happen. Now, arguing word for word
as in \cite[Subsection 3.2]{Pukh2022a} and using \cite[Theorem
3.1]{Pukh2022a}, we get that on the quadric $E$ there is a
hyperplane section $\Lambda\subset E$, such that
$$
\nu+\mathop{\rm mult}\nolimits_{\Lambda}D^+_F>2n(D_F).
$$
Taking the linear system $|H_F-\Lambda|$ (of the projective
dimension $(k-1)$) of hyperplane sections of the variety $F$, a
general divisor $W\in|H_F-\Lambda|$ in which contains the point
$o$ and its strict transform $W^+$ cuts out $\Lambda$ on $E$ (that
is, $W^+\cap E=\Lambda$), we set $D_W=(D_F\circ W)$ and obtain the
inequality
$$
\mathop{\rm mult}\nolimits_oD_W=2(\nu+\mathop{\rm
mult}\nolimits_{\Lambda}D^+_F)>4n(D_F)=4n(D_W),
$$
which contradicts Theorem 1.3. This completes the proof of
Proposition 1.3.\vspace{0.3cm}


{\bf 1.7. Exclusion of the multi-quadratic case.} This is the
hardest and the longest part of our work. Fix a point $o\in B$ of
general position, which by what was proven is a multi-quadratic
singularity of type $2^l$, satisfying the conditions (MQ1,2). The
pair $(F,\frac{1}{n(D_F)}D_F)$ has a non-canonical singularity,
the centre $B$ of which is a component of the maximal dimension of
the set $\mathop{\rm CS}(F,\frac{1}{n(D_F)}D_F)$, so that in a
neighborhood of the point $o$ this pair is canonical outside $B$.
We will show that this is impossible. This will be done in a few
steps, and now we describe the scheme of the proof and state the
key intermediate claims.

{\bf Definition 1.3.} A pair $[X,o]$, where
$$
X\subset{\mathbb P}(X)={\mathbb P}^{N(X)}
$$
is an irreducible reduced factorial complete intersection of type
$\underline{d}$ in the projective space ${\mathbb P}(X)$,
$\mathop{\rm dim}X=N(X)-k\geqslant 3$, and $o\in X$ is a point, is
called a {\it complete intersection with a marked point} or, for
brevity, a {\it marked complete intersection of level}
$(l_X,c_X)$), where $l_X$, $c_X$ are positive integers, satisfying
the inequalities
$$
2\leqslant l_X\leqslant k\quad\mbox{and} \quad c_X\geqslant l_X+4,
$$
if the following conditions are satisfied:

(MC1) the inequality
$$
\mathop{\rm codim}(\mathop{\rm Sing}X\subset X)\geqslant c_X
$$
holds,

(MC2) the point $o\in X$ is a multi-quadratic singularity of type
$2^{l_X}$, the rank of which satisfies the inequality
$$
\mathop{\rm rk}(o\in X)\geqslant 2l_X+c_X-1,
$$

(MC3) the non-singular part $X\setminus\mathop{\rm Sing}X$ of the
variety $X$ satisfies the condition of divisorial canonicity,
$$
\mathop{\rm ct}(X\setminus\mathop{\rm Sing}X)\geqslant 1,
$$
that is, for every effective divisor $A\sim aH_X$ we have
$\mathop{\rm CS}(X,\frac{1}{a}A)\subset \mathop{\rm Sing}X$.

The non-singular set of integers
$$
I_X=[k+l_X+3,k+c_X-1]\cap{\mathbb Z}
$$
is called the {\it admissible set} of the marked complete
intersection $[X,o]$.

{\bf Remark 1.1.} (i) Since $X\subset{\mathbb P}(X)$ is a complete
intersection, the factoriality of the variety $X$ follows from
Grothendieck's theorem \cite{CL} by the condition (MC1). For that
reason $\mathop{\rm Pic}X={\mathbb Z}H_X$, where $H_X$ is the
class of a hyperplane section.

(ii) By (MC1) for every $m\leqslant k+c_X-1$ and a general
subspace $P\ni o$ of dimension $m$ in ${\mathbb P}(X)$ the point
$o$ is the only singularity of the variety $X\cap P$.

(iii) Let ${\mathbb P}(X)^+\to{\mathbb P}(X)$ be the blow up of
the point $o$ with the exceptional divisor ${\mathbb
E}_X\cong{\mathbb P}^{N(X)-1}$. The strict transform
$X^+\subset{\mathbb P}(X)^+$ is the result of blowing up the point
$o$ on $X$ with the exceptional divisor $E_X=X^+\cap{\mathbb
E}_X$. Obviously, $E_X$ is an irreducible reduced non-degenerate
complete intersection of $l_X$ quadrics in a linear subspace of
codimension $(k-l_X)$ in ${\mathbb E}_X$ (this follows from (MC2),
see Proposition 1.4).

{\bf Proposition 1.4.} {\it The following inequality holds:}
$$
\mathop{\rm codim}(\mathop{\rm Sing}E_X\subset E_X)\geqslant c_X.
$$

{\bf Proof} see in \S 4 (Subsection 4.2; by the condition (MQ2)
the claim of the proposition follows from Proposition 4.2, (ii)).

{\bf Remark 1.2.} Proposition 1.4 implies the estimate
$$
\mathop{\rm codim}(\mathop{\rm Sing}E_X\subset{\mathbb
E}_X)\geqslant k+c_X.
$$
Therefore for every $m\leqslant k+c_X$ and a general subspace
$P\ni o$ of dimension $m$ in ${\mathbb P}(X)$ the strict transform
$P^+\subset{\mathbb P}(X)^+$ does not meet the set $\mathop{\rm
Sing}E_X$, since $P^+\cap{\mathbb E}_X$ is a general linear
subspace of dimension $m-1\leqslant k+c_X-1$ in ${\mathbb E}_X$.
Therefore, for $m=k+c_X$ an isolated, and for $m\leqslant k+c_X-1$
the unique singularity $o$ of the variety $X\cap P$ is resolved by
the blow up of that point, and moreover the exceptional divisor
$$
E_{X\cap P}=P^+\cap E_X
$$
of that blow up is a non-singular complete intersection of $l_X$
quadrics in the linear subspace of codimension $(k-l_X)$ in
${\mathbb E}_{X\cap P}=P^+\cap{\mathbb E}_X$. The discrepancy of
that exceptional divisor is
$$
a(E_{X\cap P})=a(E_{X\cap P}, X\cap P)=m-1-k-l_X,
$$
so that for $m=k+l_X+3$ we have: $a(E_{X\cap P})=2$. The meaning
of the lower end of the admissible set is in that equality.

In the following definition we use the notations of Remarks 1.1
and 1.2. We continue to consider a marked complete intersection
$[X,o]$ of level $(l_X,c_X)$.

{\bf Definition 1.4.} A triple $(X,D,o)$, where $D\sim n(D)H_X$ is
an effective divisor on $X$, $n(D)\geqslant 1$, is called a {\it
working triple}, if for a general subspace $P\ni o$ of dimension
$k+c_X-1$ in ${\mathbb P}(X)$ the pair
\begin{equation}\label{28.11.22.1}
\left(X\cap P,\frac{1}{n(D)}D|_{X\cap P}\right)
\end{equation}
is not log canonical at the point $o$.

{\bf Remark 1.3.} Since the point $o$ is the unique singularity of
the variety $X\cap P$, and by (MC3) the pair (\ref{28.11.22.1}) is
canonical outside the point $o$, there is a non log canonical
singularity of that pair, the centre of which on $X\cap P$ is
precisely the point $o$. By inversion of adjunction, the same is
true for a general subspace $P\ni o$ of dimension $m\leqslant
k+c_X-2$.

Let us introduce one more notation. For the strict transform $D^+$
of the divisor $D$ on $X^+$ write
$$
D^+\sim n(D)H_X-\nu(D)E_X
$$
(in order to simplify the notations, the pull back of the
divisorial class $H_X$ on $X^+$ is denoted by the same symbol
$H_X$). Respectively, for a general subspace $P\ni o$ in ${\mathbb
P}(X)$ of dimension $m\leqslant k+c_X-1$ we have
$$
D_P=D|_{X\cap P}\sim n(D)H_{X\cap P}
$$
and
$$
D^+_P\sim n(D)H_{X\cap P}-\nu(D)E_{X\cap P},
$$
where $H_{X\cap P}=H_X|_{X\cap P}$ is the class of a hyperplane
section of the variety $X\cap P\subset P\cong{\mathbb P}^m$.

{\bf Proposition 1.5.} {\it Assume that $c_X\geqslant 2l_X+4$.
Then the inequality} $\nu(D)>n(D)$ {\it holds.}

{\bf Proof} is given in \S 3 (Subsection 3.2).

Let us come back to the task of excluding the multi-quadratic
case. Recall that $F\in{\cal F}$, so that we can use the
conditions (MQ1,2) and the statement of Theorem 1.4. We fix a
point of general position $o\in B$, where $B$ is an irreducible
component of the maximal dimension of the closed set $\mathop{\rm
CS}(F,\frac{1}{n(D_F)}D_F)$.

{\bf Proposition 1.6.} {\it The pair $[F,o]$ is a marked complete
intersection of level $(l,c_F)$, where $c_F=4k+2\varepsilon(k)$,
and $(F,D_F,o)$ is a working triple.}

{\bf Proof} is given in \S 3 (Subsection 3.1).

Assume now that $l\leqslant k-1$. The symbol $T_oF$ stands again
for a subspace of codimension $(k-l)$ of the projective space
${\mathbb P}^{M+k}$. Set
$$
T=F\cap T_oF.
$$
This is subvariety of codimension $(k-l)$ in $F$ and a complete
intersection of type $\underline{d}$ in ${\mathbb P}(T)=T_oF$.

{\bf Remark 1.4.} Let us state here two well known facts which we
will use many times in the sequel: when a quadratic form is
restricted to a hyperplane, its rank either remains the same or
drops by 1 or 2; when a complete intersection in the projective
space is intersected with a hyperplane, the codimension of its
singular locus either remains the same or drops by 1 or 2 (for a
proof of the second claim, see \cite{IP} or \cite{Pukh00a}).

If $l=k$, then for uniformity of notations we set $T=F$.

{\bf Proposition 1.7.} {\it The pair $[T,o]$ is a marked complete
intersection of level $(k,c_T)$, where $c_T=2k+2\varepsilon(k)+4$.
There is an effective divisor $D_T\sim n(D_T)H_T$ on $T$, such
that $(T,D_T,o)$ is a working triple.}

{\bf Proof} is given in \S 3 (Subsection 3.4).

Proposition 1.5 (taking into account Remark 1.4) implies that
$\nu(D_T)>n(D_T)$. Now the main stage in the exclusion of the
multi-quadratic case (and thus in the proof of Theorem 0.1) is
given by the following claim.

{\bf Proposition 1.8.} {\it There is a sequence of marked complete
intersections
$$
[R_0=T,o],\quad [R_1,o],\quad\dots,\quad [R_a,o],
$$
where $a\leqslant\varepsilon(k)$ and ${\mathbb P}(R_{i+1})$ is a
hyperplane in ${\mathbb P}(R_{i})$, containing the point $o$, and
of effective divisors $D_i\sim n(D_i)H_{R_i}$ on $R_i$,
$n(D_i)\geqslant 1$, such that $D_0=D_T$ and
$$
(R_0,D_0,o),\quad (R_1,D_1,o),\quad\dots,\quad (R_a,D_a,o)
$$
are working triples, and moreover for every $i=0,\dots,a-1$ the
inequality
$$
2-\frac{\nu(D_{i+1})}{n(D_{i+1})}<\frac{1}{1+\frac{1}{k}}\left(
2-\frac{\nu(D_i)}{n(D_i)}\right)
$$
holds and} $\nu(D_a)>\frac32 n(D_a)$.

{\bf Proof} is given in \S 3 (Subsection 3.5) and \S 5.

Now let us complete the exclusion of the multi-quadratic case. The
variety $R_a$ is a section of $T=F\cap T_oF$ by a subspace of
codimension $\leqslant\varepsilon(k)$, containing the point $o$,
and $D_a$ is an effective divisor on $R_a$, satisfying the
inequality
$$
\mathop{\rm mult}\nolimits_o D_a=2^k\nu(D_a)>\frac32\cdot
2^kn(D_a).
$$
This contradicts Theorem 1.4.

The contradiction completes the proof of divisorial canonicity of
the variety $F\in{\cal F}$.


\section{Fano-Mori fibre spaces}

In this section we prove Theorem 0.2. In Subsection 2.1 we
associate with a birational map $\chi\colon V\dashrightarrow V'$ a
mobile linear system $\Sigma$ on $V$ and state the key Theorem 2.1
about this system. In Subsection 2.2 we construct a fibre-wise
birational modification of the fibre space $V/S$ for the system
$\Sigma$. In Subsection 2.3 we consider a mobile algebraic family
of irreducible curves ${\cal C}$ on $V$, and use it to prove (in
Subsection 2.4) Theorem 2.1, which implies the first claim of
Theorem 0.2 (that $\chi$ is fibre-wise). In Subsection 2.5 we
prove the birational rigidity of the fibre space
$V/S$.\vspace{0.3cm}


{\bf 2.1. The mobile linear system $\Sigma$.} Assume that the
Fano-Mori fibre space $\pi\colon V\to S$ satisfies all conditions
of Theorem 0.2. Fix a fibre space $\pi'\colon V'\to S'$ that
belongs to one of the two classes: either the class of rationally
connected fibre spaces (and then we say that the rationally
connected case is being considered), or the class of Mori fibre
spaces in the sense of Subsection 0.2 (and then we say that the
case of a Mori fibre space is being considered). We will study
both cases simultaneously.

In the rationally connected case let $Y'\ni\mathop{\rm Pic}S'$ be
a very ample class. Set
$$
\Sigma'=|(\pi')^*Y'|=|-mK_V'+(\pi')^*Y'|,
$$
where $m=0$. This is a mobile complete linear system on $V'$ (it
defines the morphism $\pi'$).

In the case of a Mori fibre space let
$$
\Sigma'=|-m K_V'+(\pi')^*Y'|
$$
be a complete linear system on $V'$, where $m\geqslant 0$ and $Y'$
is a very ample divisorial class on $S'$, and moreover, for
$m\geqslant 1$ the system $\Sigma'$ is very ample.

In both cases set
$$
\Sigma=(\chi^{-1})_*\Sigma'\subset|-n K_V+{\pi}^*Y|
$$
to be the strict transform of $\Sigma'$ on $V$ with respect to the
birational map $\chi\colon V\dashrightarrow V'$. Note that if
$m=0$ and $n=0$, then by construction of these linear systems the
map $\chi$ is fibre-wise.

{\bf Theorem 2.1.} {\it The following inequality holds:}
$n\leqslant m$.

{\bf Proof.} Assume the converse: $n>m$.  In particular, if $m=0$,
then $\chi$ is not fibre-wise. Let us show that this assumption
leads to a contradiction.\vspace{0.3cm}


{\bf 2.2. A fibre-wise birational modification of the fibre space
$V/S$.} Let $\sigma_S\colon S^+\to S$ be a composition of blow ups
with non-singular centres,
$$
S^+=S_N\stackrel{\sigma_{S,N}}{\to}S_{N-1}
\to\dots\stackrel{\sigma_{S,1}}{\to}S_0=S,
$$
where $\sigma_{S,i+1}\colon S_{i+1}\to S_i$ blows up a
non-singular subvariety $Z_{S,i}\subset S_i$. Set
$V_i=V\times_SS_i$ and $\pi_i\colon V_i\to S_i$; by the assumption
on the stability with respect to birational modifications of the
base $V_i/S_i$ is a Fano-Mori fibre space. Obviously,
$$
V_{i+1}=V_i\times_{S_i}S_{i+1}
$$
is the result of the blow up $\sigma_{i+1}\colon V_{i+1}\to V_i$
of the subvariety $Z_i=\pi^{-1}_i(Z_{S,i})\subset V_i$. Therefore,
we get the commutative diagram
$$
\begin{array}{ccccccccccccccc}
V^+ & = & V_N & \stackrel{\sigma_N}{\to} & \dots & \to & V_{i+1} &
\stackrel{\sigma_{i+1}}{\to} & V_i & \to & \dots &
\stackrel{\sigma_1}{\to} & V_0 & = & V \\
  &  & \downarrow &   &  \dots &  & \downarrow &   & \downarrow &
    &  \dots &  & \downarrow  &  & \\
S^+ & = & S_N & \stackrel{\sigma_{S,N}}{\to} & \dots & \to &
S_{i+1} & \stackrel{\sigma_{S,i+1}}{\to} & S_i & \to & \dots &
\stackrel{\sigma_{S,1}}{\to} & S_0 & = & S,
\end{array}
$$
where the vertical arrows $\pi\colon V_i\to S_i$ are Fano-Mori
fibre spaces. The symbol $\Sigma^i$ stands for the strict
transform of the system $\Sigma$ on $V_i$, $\Sigma^+=\Sigma^N$. In
these notations, let us consider a sequence of blow ups
$\sigma_{S,*}$ such that for every $i=0,1,\dots,N-1$
$$
Z_i\subset\mathop{\rm Bs}\Sigma^i,
$$
and the base set of the system $\Sigma^+$ contains entirely no
fibre $\pi^{-1}_+(s_+)$, where $s_+\in S^+$ and $\pi_+=\pi_N$. (If
this is true already for the original system $\Sigma$, then we set
$\sigma_S=\mathop{\rm id}_S$, $S^+=S$ and $V^+=V$, there is no
need to make any blow ups; but we will soon see that this case is
impossible.)

By the assumptions on the fibre space $V/S$ the fibre
$\pi^{-1}_+(s_+)$ is isomorphic to the fibre $F_s=\pi^{-1}(s)$ of
the original fibre space, where $s=\sigma_S(s_+)$. Let ${\cal T}$
be the set of all prime $\sigma_S$-exceptional divisors on $S^+$.
We get:
$$
\Sigma^+\subset\left|-n\sigma^*K_V+\pi^*_+\left(\sigma^*_SY-\sum_{T\in{\cal
T}}b_TT\right)\right|=
$$
$$
=\left|-nK^++\pi^*_+\left(\sigma^*_SY+\sum_{T\in{\cal T}}
(na_T-b_T)T\right)\right|,
$$
where $\sigma\colon V^+\to V$ is the composition of the morphisms
$\sigma_i$, $K^+=K_{V^+}$, $b_T\geqslant 1$ and $a_T\geqslant 1$
for all $T\in{\cal T}$, $a_T=a(T,S)$ is the discrepancy of $T$
with respect to $S$.

Let $\varphi\colon\widetilde{V}\to V^+$ be the resolution of
singularities of the composite map $\chi_+=\chi\circ\sigma\colon
V^+\dashrightarrow V'$, ${\cal E}$ the set of prime
$\varphi$-exceptional divisors on $\widetilde{V}$ and
$\psi=\chi\circ\sigma\circ\varphi\colon\widetilde{V}\to V'$ is a
birational morphism.

{\bf Proposition 2.1.} {\it For a general divisor $D^+\in\Sigma^+$
the pair $(V^+,\frac{1}{n}D^+)$ is canonical.}

{\bf Proof.} Assume that this is not the case. Then there is an
exceptional divisor $E\in{\cal E}$, satisfying the Noether-Fano
inequality
$$
\mathop{\rm ord}\nolimits_ED^+=\mathop{\rm
ord}\nolimits_E\Sigma^+>na(E,V^+)
$$
(we write $D^+,\Sigma^+$ instead of
$\varphi^*D^+,\varphi^*\Sigma^+$ for simplicity). Set
$B=\varphi(E)\subset V^+$.

There are two options:

(1) $\pi_+(B)=S^+$,

(2) $\pi_+(B)$ is a proper irreducible closed subset $S^+$.

If (1) is the case, then the fibre $F=F_s$ of general position
intersects $B$. The restriction
$$
\Sigma^+_F=\Sigma^+|_F\subset|-nK_F|
$$
is a mobile linear system, and moreover, the pair
$(F,\frac{1}{n}D^+_F)$ is not canonical for $D^+_F=D^+|_F$. This
contradicts the condition $\mathop{\rm mct}(F)\geqslant 1$.

Therefore, (2) is the case. Let $p\in B$ be a point of general
position and $F=\pi^{-1}_+(\pi_+(p))$, so that $p\in F$. Since
$F\not\subset\mathop{\rm Bs}\Sigma^+$, the restriction
$D^+_F=D^+|_F$ is well defined (although the linear system
$\Sigma^+_F$ may have fixed components). By inversion of
adjunction the pair $(F,\frac{1}{n}D^+_F)$ is not log canonical.
This contradicts the condition $\mathop{\rm lct}(F)\geqslant 1$.
Q.E.D. for the proposition.

Denote by the symbol $\widetilde{\Sigma}$ the strict transform of
the system $\Sigma^+$ on $\widetilde{V}$. Obviously,
\begin{equation}\label{05.11.22.1}
\widetilde{\Sigma}=\psi^*\Sigma'=|-m\psi^*K'+\psi^*(\pi')^*Y'|,
\end{equation}
where $K'=K_{V'}$, that is, $\widetilde{\Sigma}$ is a complete
linear system. We have another presentation for this linear
system:
$$
\widetilde{\Sigma}=\left|\varphi^*D^+-\sum_{E\in{\cal
E}}b_EE\right|=
$$
\begin{equation}\label{05.11.22.2}
=\left|-n\widetilde{K}+\varphi^*\pi^*_+\left(\sigma^*_SY+\sum_{T\in{\cal
T}}(na_T-b_T)T\right)+\sum_{E\in{\cal E}}(na_E-b_E)E\right|,
\end{equation}
where $\widetilde{K}=K_{\widetilde{V}}$, $D^+\in\Sigma^+$ is a
general divisor and $a_E=a(E,V^+)$ is the
discrepancy.\vspace{0.3cm}


{\bf 2.3. The mobile system of curves.} Take a family of
irreducible curves ${\cal C'}$ on $V'$, contracted by the
projection $\pi'$, sweeping out a Zariski dense subset of the
variety $V'$ and not meeting the set where the birational map
$\psi^{-1}$ is not determined. Assume that for a general pair of
points $p,q$ in a fibre of general position of the projection
$\pi'$ there is a curve $C'\in {\cal C'}$, containing the both
points. In the rationally connected case the curves of the family
${\cal C'}$ are rational (the existence of such family is shown in
\cite[Chapter II]{Kol96}), in the case of a Mori fibre space we do
not require this. For a curve $C'\in{\cal C'}$ set
$\widetilde{C}=\psi^{-1}(C')$ (at every point of the curve $C'$
the map $\psi^{-1}$ is an isomorphism), thus we get a family
$\widetilde{{\cal C}}$ of irreducible curves on $\widetilde{V}$.
Both in the rationally connected case and the case of a Mori fibre
space the inequality
$$
(C'\cdot K')<0
$$
holds, so that $(\widetilde{C}\cdot\widetilde{K})=(C'\cdot K')<0$.
Furthermore,
$$
(\widetilde{C}\cdot\widetilde{D})=(C'\cdot D')=-m(C'\cdot
K')\geqslant 0,
$$
and $(\widetilde{C}\cdot\widetilde{D})=0$ if and only if $m=0$
(since obviously $(C'\cdot(\pi')^*Y')=0$).

Let ${\cal C}^+=\varphi_*\widetilde{{\cal C}}$ be the image of the
family $\widetilde{{\cal C}}$ on $V^+$ and ${\cal C}
=\sigma_*{\cal C}^+$ its image on $V$.

{\bf Proposition 2.2.} {\it The curves $C\in{\cal C}$ are not
contracted by the projection} $\pi$.

{\bf Proof.} Assume the converse: $\pi(C)$ is a point on $S$. By
the construction of the family ${\cal C}'$ this means that the map
$\chi^{-1}$ is fibre-wise: there is a rational dominant map
$\beta'\colon S'\dashrightarrow S$, such that the diagram
$$
\begin{array}{ccc}
V & \stackrel{\phantom{xxx}\chi^{-1}}{\dashleftarrow} & V'\\
\downarrow &   &  \downarrow \\
S & \stackrel{\phantom{xx}\beta'}{\dashleftarrow} & S'
\end{array}
$$
is commutative, and moreover, $\dim S' > \dim S$ (otherwise
$\beta'$ is birational and then $\chi$ is fibre-wise, contrary to
our assumption). In that case for a point $s\in S$ of general
position the fibre $F_s=\pi^{-1}(s)$ is birational to
$(\pi')^{-1}(\beta')^{-1}(s)$. Here $\dim
(\beta')^{-1}(s)\geqslant 1$ and either the fibre
$(\pi')^{-1}(s')$ for a point $s'\in (\beta')^{-1}(s)$ of general
position is rationally connected, or the anti-canonical class of
the variety $(\pi')^{-1}(\beta')^{-1}(s)$ is $\pi'$-ample, and we
get a contradiction with the condition $\mathop{\rm mct}
(F_s)\geqslant 1$ (the fibre $F_s$ is a birationally superrigid
Fano variety). Q.E.D. for the proposition.

For a general curve $C\in{\cal C}$ set
$$
\pi_*C=d_C\overline{C},
$$
where $d_C\geqslant 1$. Replacing, if necessary, the family $\cal
{C}'$ by some open subfamily, we may assume that the integer $d_C$
does not depend on $C$. For the corresponding curve
$C^+\in\cal{C}^+$ we have $(\pi_+)_*C^+=d_C\overline{C}^+$, where
$\overline{C}^+$ is the strict transform of the curve
$\overline{C}$ on $S^+$.\vspace{0.3cm}


{\bf 2.4. Proof of Theorem 2.1.} Recall that we assume that $n>m$.
Using the two presentations (\ref{05.11.22.1}) and
(\ref{05.11.22.2}) for the class of a divisor
$\widetilde{D}\in\widetilde{\Sigma}$, we get
$$
d_C\left(\overline{C}^+\cdot\left(\sigma^*_SY+\sum_{T\in{\cal
T}}(na_T-b_T)T\right)\right)+\sum_{E\in{\cal
E}}(na_E-b_E)(\overline{C}\cdot E)=
(n-m)(\widetilde{C}\cdot\widetilde{K})<0,
$$
whence, taking into account the inequalities $b_E\leqslant na_E$
for all $E\in{\cal E}$ (Proposition 2.1), it follows that
$$
\left(\overline{C}^+\cdot \left(\sigma^*_SY+\sum_{T\in{\cal
T}}(na_T-b_T)T\right)\right)<0.
$$
However, the class $Y$ is pseudo-effective, so that
$$
(\overline{C}^+\cdot\sigma^*_S Y)=(\overline{C}\cdot Y)\geqslant
0,
$$
and $(\overline{C}^+\cdot T)\geqslant 0$ for all $T\in{\cal T}$,
so that ${\cal T}\neq\emptyset$ and for some $T\in{\cal T}$, such
that $(\overline{C}^+\cdot T)>0$, the inequality $b_T>na_T$ holds.
Since $a_T\geqslant 1$ for all $T\in{\cal T}$, we conclude that
$$
\left(\overline{C}^+\cdot\left(\sigma^*_SY-\sum_{T\in{\cal T}} b_T
T \right)\right)< -n\left(\overline{C}^+\cdot \sum_{T\in{\cal T}}
a_T T \right)\leqslant -n.
$$
For a general curve $\overline{C}^+$ consider the algebraic cycle
of the scheme-theoretic intersection
$$
(D^+\circ \pi_+^{-1}(\overline{C}^+))=\left(\left(\sigma^*
D-\pi_+^*\left(\sum_{T\in{\cal T}}b_T
T\right)\right)\circ\pi_+^{-1}(\overline{C}^+)\right).
$$
The numerical class of that effective cycle is
$$
n(\sigma^*(-K_V)\cdot\pi^{-1}_+(\overline{C}^+))+\left(\overline{C}^+
\cdot\left(\sigma^*_S Y-\sum_{T\in{\cal T}}b_TT\right)\right)F
$$
(where $F$ is the class of a fibre of the projection $\pi_+$), and
the class of the effective cycle
$\sigma_*(D^+\circ\pi^{-1}_+(\overline{C}^+ ))$ in the numerical
Chow group is
$$
-n(K_V\cdot\pi^{-1}(\overline{C}))+bF,
$$
where $b<-n$. This contradicts the condition (iii) of Theorem 0.2.
the proof of Theorem 2.1 is complete. Therefore, in both cases
(that of a rationally connected fibre space and of a Mori fibre
space) the map $\chi$ is fibre-wise. The first claim of Theorem
0.2 (in the rationally connected case) is shown. It remains to
prove the birational rigidity.\vspace{0.3cm}


{\bf 2.5. Proof of birational rigidity.} Starting from this
moment, we assume that $V'/S'$ is a Mori fibre space and the
birational map $\chi\colon V\dashrightarrow V'$ is fibre-wise,
however, the corresponding map of the bases $\beta\colon
S\dashrightarrow S'$ is not birational: $\mathop{\rm
dim}S>\mathop{\rm dim}S'$ and the fibres $\beta^{-1}(s')$ for
$s'\in S'$ are of positive dimension. We have to obtain a
contradiction, showing that this case is impossible.

First of all, let us consider the fibre-wise modification of the
fibre space $V/S$ (Subsection 2.2). Now we will need a composition
of blow ups $\sigma_S\colon S^+\to S$ with non-singular centres
such that as in Subsection 2.2, none of the fibres of the
Fano-Mori fibre space $V^+/S^+$ is contained in the base set
$\mathop{\rm Bs}\Sigma^+$ and, in addition, $\sigma_S$ resolves
the singularities of the rational dominant map $\beta\colon
S\dashrightarrow S'$, that is,
$$
\beta_+=\beta\circ\sigma_S\colon S^+\to S'
$$
is a morphism. (So that the inclusion
$Z_i=\pi^{-1}_i(Z_{S,i})\subset\mathop{\rm Bs}\Sigma^i$, see
Subsection 2.2, no longer takes place for all $i=0,\dots,N-1$.)

The fibre $\beta^{-1}_+(s')$ over a point $s'\in S'$ of general
position is an irreducible non-singular subvariety of positive
dimension. Set $G(s')=(\pi')^{-1}(s')$ and let $G^+(s')$ be the
strict transform of $G(s')$ on $V^+$. Obviously,
$$
G^+(s')=\pi^{-1}_+(\beta^{-1}_+(s'))
$$
is a union of fibres of the projection $\pi_+$ over the points of
the variety $\beta^{-1}_+(s')$.

Since $\pi'\colon V'\to S'$ is a Mori fibre space, we have the
equality $\rho(V')=\rho(S')+1$. Let ${\cal E}'$ be the set of all
$\psi$-exceptional divisors $E'\in{\cal E}'$, satisfying the
equality $\pi'(\psi'(E'))=S'$. Furthermore, let ${\cal
Z}\subset\mathop{\rm Pic}\widetilde{V}\otimes{\mathbb Q}$ be the
subspace, generated by the subspace $\psi^*(\pi')^*\mathop{\rm
Pic}S'\otimes{\mathbb Q}$ and the classes of all
$\psi$-exceptional divisors on $\widetilde{V}$, the images of
which on $V'$ do not cover $S'$. Then the equality
$$
\mathop{\rm Pic}\widetilde{V}\otimes{\mathbb Q}={}\mathbb
Q\widetilde{K}\oplus\left(\bigoplus_{E'\in{\cal E}'}{\mathbb
Q}E'\oplus {\cal Z}\right)
$$
holds, in particular, the subspace in brackets is a hyperplane in
$\mathop{\rm Pic}\widetilde{V}\otimes{\mathbb Q}$. Writing down
the class $\widetilde{K}$ with respect to the morphisms $\varphi$
and $\psi$, we get the equality
\begin{equation}\label{08.11.22.1}
\varphi^*K^++\sum_{E\in{\cal E}}a^+_EE=\psi^* K'+\sum_{E'\in{\cal
E}'}a'(E')E'+Z_1,
\end{equation}
where $Z_1\in{\cal Z}$ is some effective class, $a^+_E=a(E,V^+)$
for $\varphi$-exceptional divisors $E\in{\cal E}$ and
$a'(E')=a(E',V')$ for $\psi$-exceptional divisors $E'\in{\cal
E}'$, covering $S'$. Here all $a^+_E\geqslant 1$ and $a'(E')>0$.
Using the $\psi$-presentation (\ref{05.11.22.1}) and the
$\varphi$-presentation (\ref{08.11.22.1}) of the divisorial class
$\widetilde{D}$ and expressing $K'$ from the formula
(\ref{08.11.22.1}), we get the following equality in $\mathbb{\rm
Pic}\widetilde{V}\otimes{\mathbb Q}$:
\begin{equation}\label{05.12.22.1}
(m-n)\varphi^*K^++\varphi^*\pi^*_+Y_++\sum_{E\in{\cal
E}}(ma^+_E-b_E)E=m\sum_{E'\in{{\cal E}}'}a'(E')E'+Z_2,
\end{equation}
where $Y_+=\sigma^*_SY+\sum_{T\in{\cal T}}(na_T-b_T)T$ and
$Z_2=mZ_1+\psi^*(\pi')^*Y'\in{\cal Z}$ is an effective class.
Applying to both sides of (\ref{05.12.22.1}) $\varphi_*$ and
restricting onto a fibre of general position of the projection
$\pi_+$, we get that
$$
(m-n)K^+|_{\pi_+^{-1}(s_+)}
$$
is an effective class. Since $m\geqslant n$ and the fibre
$\pi_+^{-1}(s_+)$ is a Fano variety, we conclude that $m=n$ and
(\ref{05.12.22.1}) turns into
\begin{equation}\label{08.11.22.2}
\varphi^*\pi^*_+Y_++\sum_{E\in{\cal
E}}(na^+_E-b_E)E=n\sum_{E'\in{{\cal E}}'}a'(E')E'+Z_2.
\end{equation}
By Proposition 2.1, we have $b_E\leqslant na^+_E$ for all
$E\in{\cal E}$. Again we apply $\varphi_*$ and get that the class
$Y_+$ is effective on $S^+$.

Now let us consider the defined above fibre $G=G(s')$ of general
position of the morphism $\pi'$ and its strict transforms
$\widetilde{G}$ on $\widetilde{V}$ and $G^+$ on $V^+$ (the symbol
$s'$ for simplicity of notations is omitted). Obviously, for every
$Z\in{\cal Z}$ we have $Z|_{\widetilde{G}}=0$. Furthermore, for
any linear combination with non-negative coefficients
$$
\left.\left(\sum_{E'\in{\cal
E'}}b'_{E'}E'\right)\right|_{\widetilde{G}}
$$
is a fixed divisor on $\widetilde{G}$. Now let $\Delta$ be a very
ample divisor on $S^+$. Then the restriction
$\varphi^*\pi^*_+\Delta|_{\widetilde{G}}$ is mobile (recall that
$\beta^{-1}_+(s')$ is a variety of positive dimension, so that
$\Delta|_{\beta^{-1}_+(s')}$ is a mobile class). Therefore,
$$
\varphi^*\pi^*_+\Delta\not\in\bigoplus_{E'\in{\cal E'}}{\mathbb
Q}E'\oplus {\cal Z},
$$
whence we conclude that
$$
\mathop{\rm Pic}\widetilde{V}\otimes{\mathbb Q}={\mathbb
Q}[\varphi^*\pi^*_+\Delta]\oplus\left(\bigoplus_{E'\in{\cal
E'}}{\mathbb Q}E'\oplus{\cal Z}\right).
$$
However, this can not be the case. Let $F^+\subset G^+$ be a fibre
of general position of the morphism $\pi_+$ and
$\widetilde{F}\subset\widetilde{G}$ its strict transform on
$\widetilde{V}$. Restricting (\ref{08.11.22.2}) onto
$\widetilde{F}$, we obtain the equality
$$
\sum_{E\in{\cal E}}(na_E^+-b_E)E|_{\widetilde{F}}=
n\sum_{E'\in{\cal E'}}a'(E')E'|_{\widetilde{F}},
$$
where on the right hand side it is a linear combination of all
divisors $E'|_{\widetilde{F}}$, $E'\in{\cal E'}$, with {\it
positive} coefficients (it is here that we use the assumption that
the singularities of the variety $V'$ are terminal, see Subsection
0.2), and on the left hand side it is a linear combination of
$\varphi$-exceptional divisors $E|_{\widetilde{F}}$, $E\in{\cal
E}$, with non-negative coefficients. Since by construction
$\pi^*_+\Delta|_{F^+}=0$, we have
$\varphi^*\pi^*_+\Delta|_{\widetilde{F}}=0$, whence it follows
that the restriction of {\it every} divisorial class in
$\mathop{\rm Pic}\widetilde{V}\otimes{\mathbb Q}$ onto
$\widetilde{F}$ is fixed (is a linear combination of
$\varphi$-exceptional divisors $E|_F,E\in{\cal E}$), which is
impossible. This contradiction completes the proof of Theorem 0.2.


\section{Hyperplane sections}

This section is an immediate follow up of \S 1: we develop the
technique of working triples and consider its first
applications.\vspace{0.3cm}


{\bf 3.1. The working triple $(F,D_F,o)$.} Let us prove
Proposition 1.6. Proposition 1.2, shown in Subsection 1.5, implies
the condition (MC3). Theorem 1.1 gives the condition (MC1) for
$c_F=4k+2\varepsilon(k)$ (the inequality $c_F\geqslant l+4$ is
satisfied in the obvious way, since $l\leqslant k$). Finally, the
condition (MQ1) gives precisely (MC2). Therefore, $[F,o]$ is
indeed a marked complete intersection of level $(l,c_F)$.

Consider a general subspace $P^{\sharp}\ni o$ of dimension $k+c_F$
in ${\mathbb P}^{M+k}$. The pair
$$
\left(F\cap P^{\sharp},\frac{1}{n(D_F)}D_F|_{F\cap
P^{\sharp}}\right)
$$
is not canonical. By (MC1) the singularities of the variety
${F\cap P^{\sharp}}$ are zero-dimensional, and moreover,
$o\in\mathop{\rm Sing}{F\cap P^{\sharp}}$ and
$$
\mathop{\rm CS}\left(F\cap P^{\sharp},\frac{1}{n(D_F)}D_F|_{F\cap
P^{\sharp}}\right)\subset\mathop{\rm Sing} (F\cap P^{\sharp}),
$$
and the point $o$ is an (isolated) centre of some non canonical
singularity of that pair. For a general subspace $P\ni o$ of
dimension $k+c_F-1$ take a general hyperplane in $P^{\sharp}$,
containing the point $o$. By inversion of adjunction we have the
equalities
$$
\{o\}=\mathop{\rm LCS}\left(F\cap P,\frac{1}{n(D_F)}D_F|_{F\cap
P}\right)= \mathop{\rm CS}\left(F\cap
P,\frac{1}{n(D_F)}D_F|_{F\cap P}\right),
$$
and this is precisely (\ref{28.11.22.1}). Q.E.D. for Proposition
1.6.

As we explained in Subsection 1.7, from now our work is
constructing a certain special sequence of working triples. This
sequence starts with the working triple $(F,D_F,o)$. In order to
construct the sequence, we will need certain facts about working
triples.\vspace{0.3cm}


{\bf 3.2. Multiplicity at the marked point.} Let us prove
Proposition 1.5. We use the notations of Subsection 1.7, work with
a working triple $(X,D,o)$, where $[X,o]$ is a marked complete
intersection. Assume that $\nu(D)\leqslant 2n(D)$ (otherwise,
there is nothing to prove).

Since for a general subspace $P\ni o$ of dimension $m\in I_X$ the
inequality $a(E_{X\cap P})\geqslant 2$ holds (see Remark 1.2), the
pair
$$
\left((X\cap P)^+,\frac{1}{n(D)}D^+_P\right)
$$
is not log canonical, and moreover,
$$
\mathop{\rm LCS}\left((X\cap
P)^+,\frac{1}{n(D)}D^+_P\right)\subset E_{X\cap P}.
$$
Let $B(P)\subset E_{X\cap P}$ be the centre of some non log
canonical singularity of that pair. Then the inequality
$$
\mathop{\rm mult}\nolimits_{B(P)}D^+_P>n(D)
$$
holds, and the more so
$$\mathop{\rm mult}\nolimits_{B(P)}D^+_P|_{E_{X\cap
P}}>n(D).
$$
Considering a general subspace $P^*\ni o$ of the minimal
admissible dimension $k+l_X+3$ in $I_X$ as a general subspace of
codimension $\geqslant l_X$ in a general subspace $P\ni o$ of the
maximal admissible dimension $k+c_X-1$ in $I_X$ (recall that by
assumption $c_X\geqslant 2l_X+4$), we see that the centre $B(P)$
of some non log canonical singularity is of dimension $\geqslant
l_X$. However, $E_{X\cap P}$ is a non-singular complete
intersection of $l_X$ quadrics in the projective space of
dimension $l_X+c_X-2$, and the divisor $D^+_P|_{E_{X\cap P}}$ is
cut out on $E_{X\cap P}$ by a hypersurface of degree $\nu(D)$ in
that projective space. Therefore (for example, by
\cite[Proposition 3.6]{Pukh06b}), the inequality
$$
\nu(D)\geqslant\mathop{\rm mult}\nolimits_{B(P)} D^+_P|_{E_{X\cap
P}}
$$
holds. Therefore, $\nu(D)>n(D)$. Q.E.D. for Proposition
1.5.\vspace{0.3cm}


{\bf 3.3. Transversal hyperplane sections.} We still work with an
arbitrary working triple $(X,D,o)$, where $[X,o]$ is a marked
complete intersection of level $(l_X, c_X)$.

{\bf Proposition 3.1.} {\it Let $R\ni o$ be the section of the
variety $X$ by a hyperplane ${\mathbb P}(R)\subset{\mathbb P}(X)$,
which is not tangent to $X$ at the point $o$. Then $D\neq bR$ for
$b\geqslant 1$. Moreover, if $D$ contains $R$ as a component, that
is,
$$
D=D^*+bR,
$$
where $b\geqslant 1$, then $(X,D^*,o)$ is a working triple.}

{\bf Proof.} If $c_X\geqslant 2l_X+4$, then the first claim (that
$D$ is not a multiple of $R$) follows immediately from Proposition
1.5: indeed, the hyperplane ${\mathbb P}(R)$ is not tangent to $X$
at the point $o$, that is, for the strict transform $R^+$ on the
blow up of that point we have
$$
R^+\sim H_X-E_X,
$$
so that the equality $D=bR$ implies that $n(D)=b=\nu(D)$, which
contradicts Proposition 1.5. However we will show now that the
additional assumptions for the parameters $l_X$ and $c_X$ are not
needed.

By Remark 1.4 the condition (MC1) for $X$ implies the inequality
\begin{equation}\label{13.09.22.1}
\mathop{\rm codim}(\mathop{\rm Sing}R\subset R)\geqslant c_X-2.
\end{equation}
Since the hyperplane ${\mathbb P}(R)$ is not tangent to $X$ at the
point $o$, this point is a multi-quadratic singularity of the
variety $R$ of type $2^{l_X}$, the rank of which (by Remark 1.4
and the condition (MC2)) satisfies the inequality
$$
\mathop{\rm rk}(o\in R)\geqslant 2l_X+c_X-3.
$$

Consider a general linear subspace $P\ni o$ in ${\mathbb P}(X)$ of
dimension $k+c_X-2$. That dimension, generally speaking, does not
belong to $I_X$ and only the inequality
$$
a(E_{X\cap P})\geqslant 1
$$
holds. The variety $X\cap P$ has a unique singularity, the point
$o$, and its strict transform $(X\cap P)^+$ and the exceptional
divisor $E_{X\cap P}$ are non-singular.

The intersection $P\cap{\mathbb P}(R)$ is a general linear
subspace of dimension $k+c_X-3$ in ${\mathbb P}(R)$, containing
the point $o$. For that reason $R\cap P$ has a unique singularity,
the point $o$, and moreover, the exceptional divisor
$$
E_{R\cap P}=(R\cap P)^+\cap{\mathbb E}_X=R^+\cap E_{X\cap P}
$$
is non-singular, and the map $(R\cap P)^+\to R\cap P$ is the blow
up of the point $o$ on $R\cap P$, which resolves the singularities
of that variety. From here, taking into account that
$$
\nu(R)=1\leqslant a(E_{X\cap P}),
$$
it follows that the pair $(X\cap P,R\cap P)$ is canonical. By
inversion of adjunction we get that for every $m\in I_X$ and a
general subspace $P^{\sharp}\ni o$ of dimension $m$ the pair
$(X\cap P^{\sharp}, R\cap P^{\sharp})$ is canonical. Therefore,
$D\neq bR$, $b\geqslant 1$, and the first claim of the proposition
is shown.

Assume now that $D=D^*+bR$, where $b\geqslant 1$. Then for a
general subspace $P\ni o$ of dimension $k+c_X-1$ in ${\mathbb
P}(X)$ the pair $(X\cap P, \frac{1}{n(D)}D_P)$ is not log
canonical at the point $o$. As we saw above, the pair $(X\cap P,
R_P)$ is log canonical (and even canonical). The condition of
being log canonical is linear, so we conclude that the pair
$$
\left(X\cap P, \frac{1}{n(D^*)} D^*|_{X\cap P}\right)
$$
is not log canonical at the point $o$. Therefore, $(X,D^*,o)$ is a
working triple. Q.E.D. for the proposition.

{\bf Theorem 3.1 (on the transversal hyperplane section).} {\it
Let $[X,o]$ be a marked complete intersection of level
$(l_X=k,c_X)$, where $c_X\geqslant k+6$, and $(X,D,o)$ a working
triple. Let $R\ni o$ be a hyperplane section, which is not a
component of the divisor $D$. Assume that the inequality
$\mathop{\rm ct}(R\backslash\mathop{\rm Sing}R)\geqslant 1$ holds.
Then $(R,(D\circ R),o)$ is a working triple on the marked complete
intersection $[R,o]$ of level $(l_R=k,c_R)$, where} $c_R=c_X-2$.

{\bf Proof.} First of all, let us check that $[R,o]$ is a marked
complete intersection of level $(k,c_R)$. The inequality
$c_R\geqslant k+4$ holds by assumption. Furthermore,
$$
\mathop{\rm codim} (\mathop{\rm Sing}R\subset R) \geqslant
c_X-2=c_R,
$$
so that the condition (MC1) is satisfied. Furthermore, the point
$o\in R$ is a multi-quadratic singularity, the rank of which
satisfies the inequality
$$
\mathop{\rm rk}(o\in R)\geqslant\mathop{\rm rk}(o\in X)-2\geqslant
2k+c_R-1,
$$
so that the condition (MC2) holds. The condition (MC3) holds by
assumption. The bound for the codimension of the singular set
$\mathop{\rm Sing}R$ guarantees that the complete intersection
$R\subset{\mathbb P}(R)$ is irreducible, reduced and factorial.
Therefore, $[R,o]$ is a marked complete intersection of level
$(k,c_R)$. Set
$$
I_R=[2k+3,k+c_R-1].
$$
Obviously, $(D\circ R)\sim n(D\circ R)H_R=n(D)H_R$, where $H_R$ is
the class of a hyperplane section of $R$. It remains to check that
for a general subspace $P\ni o$ of dimension $k+c_R-1$ in
${\mathbb P}(R)$ the pair
\begin{equation}\label{14.09.22.1}
\left(R\cap P\,\frac{1}{n(D\circ R)}(D\circ R)|_{R\cap P}\right)
\end{equation}
is not log canonical at the point $o$. In order to do this, we
present $P$ as the intersection
$$
P=P^{\sharp}\cap{\mathbb P}(R),
$$
where $P^{\sharp}\ni o$ is a general subspace of dimension
$$
k+c_R=k+c_X-2
$$
in ${\mathbb P}(X)$. As $k+c_R\in I_X$, the point $o$ is the only
singularity of the variety $X\cap P^{\sharp}$ (and the only
singularity of the variety $R\cap P$), and
$$
\{o\}=\mathop{\rm LCS}\left(X\cap
P^{\sharp},\frac{1}{n(D)}D|_{X\cap P^{\sharp}}\right).
$$
The variety $R\cap P$ is the section of the variety $X\cap
P^{\sharp}$ by the hyperplane $P= P^{\sharp}\cap{\mathbb P}(R)$,
containing the point $o$, so that by inversion of adjunction the
pair (\ref{14.09.22.1}) is not log canonical. At the same time, it
is canonical outside the point $o$ since the subspace
$P\subset{\mathbb P}(R)$ is generic, the non-singular part
$R\backslash\mathop{\rm Sing}R$ is divisorially canonical and the
equality $\{o\}=\mathop{\rm Sing}(R\cap P)$ holds. Therefore, the
pair (\ref{14.09.22.1}) is not log canonical precisely at the
point $o$, which completes the proof of Theorem 3.1.\vspace{0.3cm}


{\bf 3.4. Tangent hyperplane sections.} Now let us consider a
marked complete intersection $[X,o]$ of level $(l_X,c_X)$, where
$l_X\leqslant k-1$. Let $R$ be the section of the variety $X$ by a
hyperplane ${\mathbb P}(R)\subset{\mathbb P}(X)$, which is tangent
to $X$ at the point $o$. By the symbol ${\mathbb P}(R)^+$ denote
the strict transform of the hyperplane ${\mathbb P} (R)$ on
${\mathbb P}(X)^+$ and set
$$
{\mathbb E}_R={\mathbb P}(R)^+\cap {\mathbb E}_X.
$$
Obviously, ${\mathbb E}_R\cong{\mathbb P}^{N(X)-2}$ is the
exceptional divisor of the blow up of the point $o$ on the
hyperplane ${\mathbb P}(R)$. Set also $E_R=R^+\cap{\mathbb E}_X$.
Obviously, $E_R\subset{\mathbb E}_R$, and
$$
\mathop{\rm codim}(E_R\subset{\mathbb E}_R)=k.
$$

{\bf Proposition 3.2.} {\it Assume that $c_X\geqslant l_X+5$ and
the point $o\in R$ is a multi-quadratic singularity of type
$2^{l_X+1}$, and moreover, the inequality
$$
\mathop{\rm rk} (o\in R)\geqslant 2l_X+c_X-2
$$
holds. Then $D\neq bR$ for $b\geqslant 1$. Moreover, if the
divisor $D$ contains $R$ as a component, that is,
$$
D=D^*+bR,
$$
where $b\geqslant 1$, then $(X,D^*,o)$ is a working triple.}

{\bf Proof} is completely similar to the proof Proposition 3.1,
but we give it in full details, because there some small points
where the two arguments are different. The inequality
(\ref{13.09.22.1}) holds in this case again. Let us use the
additional assumption about the singularity $o\in R$. Consider a
general linear subspace $P\ni o$ in ${\mathbb P}(X)$ of dimension
$k+l_X+3$ (it is the minimal admissible dimension). We get the
equality $a(E_{X\cap P})=2$. Obviously,
$$
R^+\sim H_X-2E_X
$$
and, respectively, on $(X\cap P)^+$ we have
$$
(R\cap P)^+\sim H_{X\cap P}-2E_{X\cap P}.
$$
Arguing as in the transversal case, we note that the intersection
$P\cap{\mathbb P}(R)$ is a general subspace of dimension $k+l_X+2$
in ${\mathbb P}(R)$. Taking into account that by the inequality
(\ref{13.09.22.1}) the inequality
$$
\mathop{\rm codim}(\mathop{\rm Sing}R\subset{\mathbb
P}(R))\geqslant k+c_X-2
$$
holds, and that by assumption $c_X\geqslant l_X+5$, we see that
$R\cap P$ has a unique singularity, the point $o$. Furthermore, by
the assumption about the rank of the singular point $o\in R$ we
get the inequality
$$
\mathop{\rm codim}(\mathop{\rm Sing}E_R\subset E_R)\geqslant
c_X-3\geqslant l_X+2,
$$
so that
$$
\mathop{\rm codim}(\mathop{\rm Sing}E_R\subset{\mathbb
E}_R)\geqslant k+l_X+2.
$$
The exceptional divisor $E_{R\cap P}$ is the section of the
subvariety $E_R\subset{\mathbb E}_R$ by a general linear subspace
of dimension $k+l_X+1$, whence we conclude that the variety
$E_{R\cap P}$ is non-singular. Thus we have shown that the
singularity $o\in R\cap P$ is resolved by one blow up. Therefore,
the pair
$$
((X\cap P)^+,(R\cap P)^+)
$$
is canonical, so that the pair
$$
(X\cap P),(R\cap P)
$$
is canonical, too. We have shown that $D\neq bR$ for $b\geqslant
1$.

By inversion of adjunction for every $m\in I_X$ and a general
subspace $P^{\sharp}\ni o$ of dimension $m$ the pair $(X\cap
P^{\sharp}, R\cap P^{\sharp})$ is canonical (recall that
$n(R)=1$). Repeating the arguments given in the transversal case
(the proof of Proposition 3.1) word for word, we complete the
proof of Proposition 3.2.

{\bf Remark 3.1.} If for $l_X\leqslant k-1$ the intersection
$X\cap T_oX$ has the point $o$ as a multi-quadratic singularity of
type $2^k$, the rank of which satisfies the inequality
$$
\mathop{\rm rk}(o\in X\cap T_oX)\geqslant 2l_X+c_X-2,
$$
the the assumption about the rank $\mathop{\rm rk}(o\in R)$ in the
statement of Proposition 3.2 holds automatically for every tangent
hyperplane at the point $o$.

{\bf Theorem 3.2 (on the tangent hyperplane section).} {\it Let
$[X,o]$ be a marked complete intersection of level $(l_X,c_X)$,
where $2\leqslant l_X\leqslant k-1$ and $c_X\geqslant l_X+7$, and
$(X,D,o)$ a working triple. Let $R$ be the section of $X$ by a
hyperplane which is tangent to $X$ at the point $o$, and assume
that $R$ is not a component of the divisor $D$. Assume that the
point $o\in R$ is a multi-quadratic singularity of type $2^{l_R}$,
where $l_R=l_X+1$, the rank of which satisfies the inequality
$$
\mathop{\rm rk}(o\in R)\geqslant 2l_R+c_R-1=2l_X+c_X-1,
$$
where $c_R=c_X-2$, and also that the inequality $\mathop{\rm
ct}(R\backslash\mathop{\rm Sing}R)\geqslant 1$ holds. Then
$(R,(D\circ R),o)$ is a working triple on the marked complete
intersection $[R,o]$ of level} $(l_R,c_R)$.

{\bf Proof} is similar to the transversal case (Theorem 3.1), and
we just emphasize the necessary modifications. The fact that
$[R,o]$ is a marked complete intersection of level $(l_R,c_R)$ is
checked in the tangent case even easier than in the transversal
one, because the assumption about the singularity $o\in R$ is
among the assumptions of the theorem.

A general subspace $P\ni o$ of dimension $k+c_R-1=k+c_X-3$ in
${\mathbb P}(R)$ is again presented as the intersection
$P=P^{\sharp}\cap{\mathbb P}(R)$, where $P^{\sharp}\ni o$ is a
general subspace of dimension $k+c_R\in I_X$ in ${\mathbb P}(X)$,
and now, repeating the arguments in the transversal case and using
inversion of adjunction, we get that $(R,(D\circ R),o)$ is a
working triple. Q.E.D. for the theorem.

{\bf Proof of Proposition 1.7.} We assume that $l\leqslant k-1$.
recall that the symbol $T$ stands for the intersection $F\cap
T_oF$; this is a subvariety of codimension $(k-l)$ in $F$. Let us
construct a sequence of subvarieties
$$
T_0=F\supset T_1\supset\dots\supset T_{k-l}=T,
$$
where $T_{i+1}$ is the section of $T_i\ni o$ by some hyperplane
${\mathbb P}(T_{i+1})=\langle T_{i+1}\rangle\ni o$, which is
tangent to $T_i$ at the point $o$. Theorem 1.1 implies that the
inequality
$$
c_F\geqslant l+3(k-l)+4
$$
holds (the inequality of Theorem 1.1 for the codimension $c_F$ is
much stronger, but for the clarity of exposition we give the
weakest estimate that is sufficient for the proof of Proposition
1.7; this remark also applies to the estimate of the rank of the
multi-quadratic singularity $o\in T$ below). Furthermore, the
condition (MQ2) implies that $o\in T$ is a multi-quadratic
singularity of type $2^k$, and moreover, the inequality
\begin{equation}\label{19.09.22.1}
\mathop{\rm rk}(o\in T)\geqslant 2k+c_F-1
\end{equation}
holds. Finally, by Theorem 1.2 for every hyperplane section $W$ of
every subvariety $T_i$, $i=0,1,\dots,k-l$, every non-singular
point $p\in W$ and every prime divisor $Y$ on $W$ the inequality
\begin{equation}\label{19.09.22.2}
\frac{\mathop{\rm mult}_p}{\mathop{\rm
deg}}Y\leqslant\frac{2}{\mathop{\rm deg} F}
\end{equation}
holds. Then for all $i=0,1,\dots,k-l$ the pair $[T_i,o]$ is a
marked complete intersection of level
$$
(l_i=l+i,c_i=c_F-2i).
$$
Indeed, the inequality $c_i\geqslant l_i+4$ is true by the
definition of the numbers $l_i$, $c_i$, the condition (MC1)
follows from Remark 1.4, the point $o\in T_i$ by construction is a
multi-quadratic singularity of type $2^{l+1}$, and moreover, by
(\ref{19.09.22.1}) we have
$$
\mathop{\rm rk}(o\in T_i)\geqslant 2l+c_F-1=2l_i+c_i-1,
$$
and, finally, repeating the proof of Proposition 1.1 and the
arguments of Subsection 1.5 word for word, we get that by
(\ref{19.09.22.2}) the condition (MC3) holds. Therefore, $[T,o]$
is a marked complete intersection of level $(k,c_{k-l})$, where
$c_{k-l}=c_F-2(k-l)\geqslant k+4$. Recall (Proposition 1.6), that
$c_F=4k+2\varepsilon(k)$. Since $l\geqslant 2$, the inequality
$$
c_{k-l}\geqslant c_T=2k+2\varepsilon(k)+4
$$
holds, so that $[T,o]$ is a marked complete intersection of level
$(k,c_T)$, as we claimed.

It remains to construct the working triple $(T,D_T,o)$. We will
construct a sequence of working triples $(T_i,D_i,o)$, where
$i=0,1,\dots,k-l$ and $D_0=D_F$. Assume that $(T_i,D_i,o)$ is
already constructed and $i\leqslant k-l-1$. Let us check that all
assumptions that allow us to apply Proposition 3.2 are satisfied.

Indeed, the fact that $i\leqslant k-l-1$ implies the inequality
$c_i\geqslant l_i+7$. The point $o\in T_{i+1}$ is a
multi-quadratic singularity of type $2^{l_i+1}$, the rank of which
satisfies the inequality
$$
\mathop{\rm rk}(o\in T_{i+1})\geqslant
2l_{i+1}+c_{i+1}-1=2l_i+c_i-1
$$
(see above). Applying Proposition 3.2, we remove $T_{i+1}$ from
the effective divisor $D_i$ (if it is necessary) and obtain the
working triple $(T_i,D^*_i,o)$, where the effective divisor
$D^*_i$ does not contain $T_{i+1}$ as a component.

it ie easy to see that we have all assumptions of Theorem 3.2. Set
$$
D_{i+1}=(D^*_i\circ T_{i+1}).
$$
Now $(T_{i+1}, D_{i+1}, o)$ is a working triple. Proof of
Proposition 1.7 is complete.\vspace{0.3cm}


{\bf 3.5. Plan of the proof of Proposition 1.8.} Recall that by
the condition (MQ2) the inequality
$$
\mathop{\rm rk} (o\in T)\geqslant 10k^2+8k+2\varepsilon(k)+5
$$
holds. The pair $[T,o]$ is a marked complete intersection of level
$(k,c_T)$, where $c_T=2k+2\varepsilon(k)+4$. Let
$$
R_0=T,R_1,\dots,R_a,
$$
where $a\leqslant\varepsilon(k)$, be an {\it arbitrary} sequence
of subvarieties in $T$, where $R_{i+1}$ is the section of $R_i$ by
the hyperplane ${\mathbb P}(R_{i+1})$ in ${\mathbb P}(R_i)$,
containing the point $o$. Set $c_i=c_T-2i$, where $i=0,1,\dots,a$.

{\bf Proposition 3.3.} {\it The pair $[R_i,o]$ is a marked
complete intersection of level} $(k,c_i)$.

{\bf Proof.} Since $a\leqslant\varepsilon(k)$, the inequality
$c_i\geqslant k+4$ holds in an obvious way (in fact, $c_i\geqslant
2k+4)$. The condition (MC1) holds by Remark 1.4. The condition
(MC3) is obtained by repeating the proof of Proposition 1.1 and
the arguments of Subsection 1.5 word for word, taking into account
Theorem 1.2. Finally, again by Remark 1.4 the inequality
$$
10k^2+8k+2\varepsilon(k)+4\geqslant 2k+c_i+2i-1
$$
implies the condition (MC2). Q.E.D. for the proposition.

Now let us construct for every $i=0,1,\dots,a$ an effective
divisor $D_i$ on $R_i$ in the same way as we did it in Subsection
3.4 in the proof of Proposition 1.7, applying instead of
Proposition 3.2 its ``transversal'' analog, Proposition 3.1, and
Theorem 3.1 instead of Theorem 3.2. More precisely, if the
effective divisor $D_i$, where $i\leqslant a-1$, is already
constructed, we remove from this divisor all components that are
hyperplane sections (if there are such components), and obtain an
effective divisor $D^*_i$ that does not contain hyperplane
sections as components, and such that $(R_i,D^*_i,o)$ is a working
triple (Proposition 3.1).

{\bf Proposition 3.4.} {\it The following inequality holds:}
$$
\frac{\nu(D^*_i)}{n(D^*_i)}\geqslant\frac{\nu(D_i)}{n(D_i)}.
$$

{\bf Proof.} It is sufficient to consider the case when $D^*_i$ is
obtained from $D_i$ by removing one hyperplane section $Z\ni o$.
Write down
$$
D_i=D^*_i+bZ,
$$
where $b\geqslant 1$. Since $c_i\geqslant 2k+4$, we can apply
Proposition 1.5: $\nu(D_i)>n(D_i)$. On the other hand,
$\nu(Z)=n(Z)=1$. Set $\nu(D_i)=\alpha n(D_i)$, where $\alpha>1$.
We get
$$
\frac{\nu(D^*_i)}{n(D^*_i)}= \frac{\alpha
n(D^*_i)+(\alpha-1)b}{n(D^*_i)}>\alpha,
$$
which proves the proposition. Q.E.D.

(If we remove from $D_i$ a hyperplane section that does not
contain the point $o$, the claim of Proposition 3.4 is trivial.)

Now we apply Theorem 3.1, setting $D_{i+1}=(D^*_i\circ R_{i+1})$,
this cycle of the scheme-theoretic intersection is well defined as
an effective divisor on $R_{i+1}$, and moreover, $(R_{i+1},
D_{i+1}, o)$ is a working triple and
$$
\frac{\nu(D_{i+1})}{n(D_{i+1})}\geqslant\frac{\nu(D^*_i)}{n(D^*_i)}
\geqslant\frac{\nu(D_i)}{n(D_i)}.
$$
We emphasize that $R_1,\dots, R_a$ is an arbitrary sequence of
consecutive hyperplane sections. By Remark 1.4, for all
$i=0,1,\dots,a$ the inequality $c_i\geqslant 2k+4$ holds, and the
rank of the multi-quadratic singularity $o\in R_i$ of type $2^k$
is at least $10k^2+8k+5$. By Theorem 1.4 and Proposition 1.5 w
have the inequalities
$$
n(D_i)<\nu(D_i)\leqslant\frac32n(D_i).
$$
Therefore, at every step of our construction the assumptions of
the following claim are satisfied.

{\bf Theorem 3.3 (on the special hyperplane section).} {\it Let
$[X,o]$ be a marked complete intersection of level $(k,c_X)$,
where $c_X\geqslant 2k+4$ and the inequality
$$
\mathop{\rm rk}(o\in X)\geqslant 10k^2+8k+5
$$
holds. Let $(X,D,o)$ be a working triple, where the effective
divisor $D$ does not contain hyperplane sections and satisfies the
inequalities
$$
n(D)<\nu(D)\leqslant\frac32n(D).
$$
Then there is a section $R\ni o$ of the variety $X$ by a
hyperplane ${\mathbb P}(R)=\langle R\rangle\subset {\mathbb P}(X)=
{\mathbb P}^{N(X)}$, such that the effective divisor $D_R=(R\circ
D)$ on $R$ satisfies the inequality}
$$
2-\frac{\nu(D_R)}{n(D_R)}<\frac{1}{1+\frac{1}{k}}\left(2-
\frac{\nu(D)}{n(D)}\right).
$$

Now by the definition of the integer $\varepsilon(k)$ and what was
said above, Theorem 3.3 immediately implies Proposition 1.8.

{\bf Proof} of Theorem 3.3 is given in \S 5.


\section{Multi-quadratic singularities}

In this section we consider the properties of multi-quadratic
singularities, the rank of which is bounded from below: they are
factorial, stable with respect to blow ups and terminal. In
Subsection 4.5 we study linear subspaces on complete intersections
of quadrics and the properties of projections from these
subspaces. \vspace{0.3cm}


{\bf 4.1. The definition and the first properties.} Let ${\cal X}$
be an (irreducible) algebraic variety, $o\in{\cal X}$ a point.

{\bf Definition 4.1.} The point $o$ s a {\it multi-quadratic
singularity} of the variety ${\cal X}$ of type $2^l$ and rank
$r\geqslant 1$, if in some neighborhood of this point ${\cal X}$
can be realized as a subvariety of a non-singular $N=(\mathop{\rm
dim} {\cal X}+ l)$-dimensional variety ${\cal Y}\ni o$, and for
some system $(u_1, \dots, u_N)$ of local parameters on ${\cal Y}$
at the point $o$ the subvariety ${\cal X}$ is the scheme of common
zeros of regular functions
$$
\alpha_1,\dots,\alpha_l\in{\cal O}_{o,{\cal Y}}\subset{\mathbb
C}[[u_1,\dots,u_N]],
$$
which are represented by the formal power series
$$
\alpha_i=\alpha_{i,2}+\alpha_{i,3}+\dots,
$$
where $\alpha_{i,j}(u_1,\dots,u_N)$ are homogeneous polynomials of
degree $j$ and
$$
\mathop{\rm rk}(\alpha_{1,2},\dots,\alpha_{l,2})=r.
$$
(Obviously, the order of the formal power series, representing
$\alpha_i$, and the rank of the tuple of quadratic forms
$\alpha_{i,2}$ do not depend on the choice of the local parameters
on ${\cal Y}$ at the point $o$.)

It is convenient to work in a more general context. Assume that in
a neighborhood of the point $o$ the variety ${\cal X}$ is realized
as a subvariety ${\cal X}\subset{\cal Z}$, where $\mathop{\rm
dim}{\cal Z}=\mathop{\rm dim}{\cal X}+e=N({\cal Z})$, and for a
certain system of local parameters $(v_1,\dots,v_{N({\cal Z})})$
on ${\cal Z}$ at the point $o$ the subvariety ${\cal X}$ is the
scheme of common zeros of regular functions
$$
\beta_1,\dots,\beta_e\in{\cal O}_{o,{{\cal Z}}}\subset{\mathbb
C}[[v_*]],
$$
which are represented by the formal power series
$$
\beta_i=\beta_{i,1}+\beta_{i,2}+\dots,
$$
where $\beta_{i,j}(v_*)$ are homogeneous polynomials of degree
$j$. Assume that for some $l\in\{0,1,\dots,e\}$
$$
\mathop{\rm dim}\langle\beta_{1,1},\dots,\beta_{e,1}\rangle=e-l,
$$
where we assume (for the convenience of notations), that the
linear forms $\beta_{j,1}$ for $l+1\leqslant j\leqslant e$ are
linearly independent, so that for $1\leqslant i\leqslant l$ and
$l+1\leqslant j\leqslant e$ there are uniquely determined numbers
$a_{i,j}$, such that
$$
\beta_{i,1}=\sum^e_{j=l+1}a_{i,j}\beta_{j,1}.
$$

Set ${\cal Y}=\{\beta_j=0\,|\, l+1\leqslant j\leqslant e\}$ and
$$
\beta_i^*=\beta_i-\sum_{j=l+1}^e a_{i,j}\beta_j.
$$
Then (in a neighborhood of the point $o$) the variety ${\cal Y}$
is non-singular, and ${\cal X}\subset {\cal Y}$ is realized as the
scheme of common zeros of the regular functions $\beta^*_i$,
$1\leqslant i\leqslant l$. Set
$$
T_o{\cal Y}=T_o{\cal X}=\{\beta_{j,1}=0\,|\, l+1\leqslant
j\leqslant e\}.
$$
If
$$
\mathop{\rm rk}\left(\beta^*_{i,2}|_{T_o{\cal X}}\, |\, 1\leqslant
i\leqslant l\right)=r,
$$
then obviously $o\in{\cal X}$ is a multi-quadratic singularity of
rank $r$.

The rank of the multi-quadratic point $o\in{\cal X}$ is denoted by
the symbol $\mathop{\rm rk} (o\in{\cal X})$ or just $\mathop{\rm
rk}(o)$, if it is clear which variety is meant. For uniformity of
notations we treat a non-singular point as a multi-quadratic one
of type $2^0$.

{\bf Proposition 4.1.} {\it Assume that $o\in{\cal X}$ is a
multi-quadratic singularity of type $2^l$, where $l\geqslant 1$,
and of rank $r\geqslant 2l$. Then in a neighborhood of the point
$o$ every point $p\in{\cal X}$ is either non-singular, or a
multi-quadratic of type $2^b$, where $b\in\{1,\dots,l\}$, of rank}
$\geqslant r-2(l-b)$.

{\bf Proof.} Using the notations for the embedding ${\cal
X}\subset{\cal Z}$, introduced above, with $e=l$ (so that
$\beta_{i,1}=0$ for all $i=1,\dots,l$) and setting $N({\cal
Z})=N$, consider an open set $U\subset{\cal Z}$, $U\ni o$, such
that for every point $p\in U$ the ``shifted'' functions
$$
v^{(p)}_i=v_i-v_i(p),\quad i=1,\dots,N,
$$
form a system of local parameters at the point $p$, and in the
formal expansion
$$
\beta_i=\beta^{(p)}_{i,0}+\beta^{(p)}_{i,1}+\beta^{(p)}_{i,2}+\dots
$$
with respect to the system of parameters $v^{(p)}_*$ the quadratic
components satisfy the inequality
$$
\mathop{\rm rk}(\beta^{(p)}_{i,2}\,|\, 1\leqslant i\leqslant
l)\geqslant r.
$$
If the point $p$ is a common zero of $\beta_1,\dots,\beta_l$, then
$\beta^{(p)}_{i,0}=0$ for $1\leqslant i\leqslant l$. Set
$$
T_p{\cal X}=\{\beta^{(p)}_{i,1}=0\,|\, 1\leqslant i\leqslant l\}
$$
and assume (for the convenience of notations) that the forms
$\beta^{(p)}_{i,1}$ for $b+1\leqslant i\leqslant l$ are linearly
independent, where
$$
\mathop{\rm dim}\langle\beta^{(p)}_{i,1}\,|\,\,1\leqslant
i\leqslant l\rangle=l-b.
$$
Since $\mathop{\rm codim}(T_p{\cal X}\subset T_p{\cal Z})=l-b$, by
Remark 1.4 the inequality
$$
\mathop{\rm rk}(\beta^{(p)}_{i,2}|_{T_p{\cal X}}\,|\, 1\leqslant
i\leqslant l) \geqslant r - 2(l-b)
$$
holds. It is easy to see from the construction of the quadratic
forms $\beta^{(p)*}_{i,2}$, $1\leqslant i\leqslant b$, that every
linear combination of these forms with coefficients
$(\lambda_1,\dots,\lambda_b)\neq(0,\dots,0)$ is a linear
combination of the original forms $\beta^{(p)}_{i,2}$, $1\leqslant
i\leqslant l$, not all coefficients in which are equal to zero.
Therefore, the point $p$ is a multi-quadratic singularity of rank
$\geqslant r-2(l-b)$, as we claimed. Q.E.D. for the
proposition.\vspace{0.3cm}


{\bf 4.2. Complete intersections of quadrics.} In the notations of
Definition 4.1 let ${\cal Y}^+\to{\cal Y}$ be the blow up of the
point $o$ with the exceptional divisor $E_{\cal Y}\cong{\mathbb
P}^{N-1}$ and ${\cal X}^+\subset{\cal Y}^+$ the strict transform
of ${\cal X}$ on ${\cal Y}^+$, so that ${\cal X}^+\to{\cal X}$ is
the blow up of the point $o$ on ${\cal X}$ with the exceptional
divisor $E_{\cal Y}|_{{\cal X}^+}=E_{\cal X}$. Therefore, $E_{\cal
X}$ is the scheme of common zeros of the quadratic forms
$\alpha_{i,2}$, $i=1,\dots,l$, on $E_{\cal Y}\cong{\mathbb
P}^{N-1}$.

Let $q_1,\dots,q_l$ be quadratic forms on ${\mathbb P}^{N-1}$,
where $N\geqslant l+4$. By the symbol $q_{[1,l]}$ we denote the
tuple $(q_1,\dots,q_l)$.

{\bf Proposition 4.2.} (i) {\it Assume that the inequality
$$
\mathop{\rm rk} q_{[1,l]}\geqslant 2l+3
$$
holds. Then the scheme of common zeros of the forms
$q_1,\dots,q_l$ is an irreducible non-degenerate factorial variety
$Q\subset{\mathbb P}^{N-1}$ of codimension $l$, that is, a
complete intersection of type $2^l$.}

(ii) {\it Assume that for some $e\geqslant 4$ the inequality
$$
\mathop{\rm rk} q_{[1,l]}\geqslant 2l+e-1
$$
holds. Then the following inequality is true:}
$$
\mathop{\rm codim}(\mathop{\rm Sing}Q\subset Q)\geqslant e.
$$

{\bf Proof} is given below in Subsection 4.4.

{\bf Corollary 4.1.} (i) {\it Assume that the rank of the tuple
$\alpha_{*,2}=(\alpha_{1,2},\dots,\alpha_{l,2})$ of quadratic
forms satisfies the inequality
$$
\mathop{\rm rk} (\alpha_{*,2})\geqslant 2l+3.
$$
Then in a neighborhood of the point $o$ the scheme of common zeros
of the regular functions $\alpha_1,\dots,\alpha_l$ is an
irreducible reduced factorial subvariety ${\cal X}$ of codimension
$l$ in ${\cal Y}$.}

(ii) {\it Assume that for some $e\geqslant 4$ the inequality
$$
\mathop{\rm rk} (\alpha_{*,2})\geqslant 2l+e-1
$$
holds. Then in a neighborhood of the point $o$ the following
inequality is true}
$$
\mathop{\rm codim}(\mathop{\rm Sing}{\cal X}\subset{\cal
X})\geqslant e.
$$

{\bf Proof.} Both claims obviously follow from Proposition 4.2,
taking into account Grothendieck's theorem \cite{CL} on the
factoriality of a complete intersection, the singular set of which
is of codimension $\geqslant 4$.

Therefore, for $r\geqslant 2l+3$ the assumption in Definition 4.1
that ${\cal X}$ is an irreducible variety, is unnecessary: in a
neighborhood of the point $o$ the scheme of common zeros of the
functions $\alpha_*$ is automatically irreducible and reduced, and
moreover, it is a factorial variety. This proves all claims of
Theorem 1.1, except for that the singularities of the variety $F$
are terminal.\vspace{0.3cm}


{\bf 4.3. Stability with respect to blow ups.} Let
$\underline{r}=(r_1,r_2,\dots,r_k)$ be a tuple of integers,
satisfying the inequalities $r_{i+1}\geqslant r_i+2$ for
$i=1,\dots,k-1$, where $r_1\geqslant 5$. Again, let ${\cal Y}$ be
a non-singular $N$-dimensional variety, where $N\geqslant k+3$,
and ${\cal X}\subset{\cal Y}$ an (irreducible) subvariety of
codimension $k$, every point $o\in{\cal X}$ of which is either
non-singular, or a multi-quadratic singularity of type $2^l$,
where $l\in\{1,\dots,k\}$, of rank $\geqslant r_l$. Somewhat
abusing the terminology, we say in this case that ${\cal X}$ has
{\it multi-quadratic singularities of type} $\underline{r}$.

{\bf Theorem 4.1.} {\it In the assumptions above let
$B\subset{\cal X}$ be an irreducible subvariety of codimension
$\geqslant 2$. Then there is an open subset $U\subset{\cal X}$,
such that $U\cap B\neq\emptyset$, $U\cap B$ is non-singular and
the blow up
$$
\sigma_B\colon U_B\to U
$$
along $B$ gives a quasi-projective variety $U_B$ with
multi-quadratic singularities of type} $\underline{r}$.

{\bf Proof.} If a point of general position $o\in B$ is
non-singular on ${\cal X}$, there is nothing to prove. If
$o\in{\cal X}$ is a multi-quadratic singularity of type $2^l$,
where $l>k$, then a certain Zariski open subset $U\subset{\cal
X}$, $U\ni o$, has multi-quadratic singularities of type
$(r_1,\dots,r_l)$ (see Subsection 4.1), so that it is sufficient
to consider the case when a point of general position $o\in B$ is
a multi-quadratic point of type $2^k$ on ${\cal X}$. Passing over
to an open subset, we may assume that the subvariety $B$ is
non-singular. Let $(u_1,\dots,u_N)$ be a system of local
parameters at the point $o$, such that $B=\{u_1=\dots=u_m=0\}$.
Since $B\subset{\cal X}$, the subvariety ${\cal X}\subset{\cal Y}$
is the scheme of common zeros of regular functions
$$
\beta_1,\dots,\beta_k\in{\cal O}_{o,{\cal Y}}\subset{\cal
O}_{o,B}[[u_1,\dots,u_m]],
$$
where for all $i=1,\dots,k$
$$
\beta_i=\beta_{i,2}+\beta_{i,3}+\dots,
$$
where $\beta_{i,j}$ are homogeneous polynomials of degree $j$ in
$u_1,\dots,u_m$ with coefficients from ${\cal O}_{o,B}$. Again
replacing ${\cal Y}$, if necessary, by an open subset, containing
the point $o$, we have
$$
\beta_i\in{\cal O}({\cal Y})\subset{\cal O}(B)[[u_1,\dots,u_m]],
$$
so that all coefficients of the forms $\beta_{i,j}$ are regular
functions on $B$; in particular,
$$
\beta_{i,2}=\sum_{1\leqslant j_1\leqslant j_2\leqslant
m}A_{j_1,j_2}u_{j_1}u_{j_2},
$$
where $A_{j_1,j_2}\in{\cal O}(B)$. In terms of the embedding
${\cal O}_{o,{\cal Y}}\subset{\mathbb C}[[u_1,\dots,u_N]]$ we get
the presentation
$$
\beta_i=\overline{\beta}_{i,2}+\overline{\beta}_{i,3}+\dots,
$$
where $\overline{\beta}_{i,j}$ is a homogeneous polynomial of
degree $j$ in $u_*$, and moreover, in the right hand side there
are no monomials that do not contain the variables
$u_1,\dots,u_m$, or that contain precisely one of them (in the
power 1): every monomial in the right hand side is divisible by
some quadratic monomial in $u_1,\dots,u_m$.

Let ${\cal Y}_B\to{\cal Y}$ be the blow up of the subvariety $B$
and ${\cal X}_B\subset{\cal Y}_B$ the strict transform of ${\cal
X}$. Obviously, the morphism ${\cal X}_B\to{\cal X}$ in the blow
up of $B$ on ${\cal X}$. The symbol $E_B$ denotes the exceptional
divisors of the blow up of $B$ on ${\cal Y}$. Since outside $E_B$
the varieties ${\cal X}_B$ and ${\cal X}$ are isomorphic, it is
sufficient to show that every point $p\in{\cal X}_B\cap E_B$ is
either non-singular, or a multi-quadratic singularity of the
variety $U_B$ of type $2^l$, where $l\geqslant 1$, and of rank
$\geqslant r_l$. We assume that the point $p$ lies over the point
$o\in U$ and is a singularity of the variety $U_B$.

By a linear change of local parameters $u_1,\dots,u_m$ we may
ensure that at the point $p\in{\cal Y}_B$ there is a system of
local parameters
$$
(v_1,\dots,v_m,u_{m+1},\dots,u_N),
$$
linked to the original system of parameters $u_*$ by the standard
relations
$$
u_1=v_1,\,\,u_2=v_1v_2,\dots,\,\,u_m=v_1v_m.
$$
The local equation of the exceptional divisor $E_B$ at the point
$p$ is $v_1=0$, and the subvariety ${\cal X}_B\subset{\cal Y}_B$
at that point is defined by the equations
$$
\widetilde{\beta}_1,\dots,\widetilde{\beta}_k\in{\cal O}_{p,{\cal
Y}_B}\subset{\mathbb C}[[v_1,\dots,v_m,u_{m+1},\dots,u_N]].
$$
Write down $\widetilde{\beta}_i=\widetilde{\beta}_{i,1}
+\widetilde{\beta}_{1,2}+\dots$ and assume that for some
$l\in\{1,\dots,k\}$ the linear forms $\widetilde{\beta}_{j,1}$,
$l+1\leqslant j\leqslant k$, are linearly independent, and
moreover,
$$
\mathop{\rm dim}\langle\widetilde{\beta}_{i,1}\,|\,1\leqslant
i\leqslant k\rangle=k-l,
$$
so that there are relations
$$
\widetilde{\beta}_{i,1}=\sum^k_{j=l+1}a_{i,j}\widetilde{\beta}_{i,1},
$$
$i=1,\dots,l$. Replacing the original system of local equations
$\beta_1,\dots,\beta_k$ by
$$
\beta_i-\sum^k_{j=l+1}a_{i,j}\beta_j,\,\,i=1,\dots,l,\quad
\beta_{l+1},\dots,\beta_k,
$$
we may assume that the linear forms $\widetilde{\beta}_{i,1}$,
$i=1,\dots,l$, are identically zero. In that case the following
claim is true.

{\bf Lemma 4.1.} {\it For $i=1,\dots,l$ the quadratic forms
$\overline{\beta}_{i,2}$ depend only on $u_2,\dots,u_m$ and
$$
\widetilde{\beta}_{i,2}=\overline{\beta}_{i,2}(v_2,\dots,v_m)+\beta^{\sharp}_{i,2},
$$
where every monomial in the quadratic form $\beta^{\sharp}_{i,2}$
is divisible either by $v_1$, or by} $u_i$, $m+1\leqslant
i\leqslant N$.

{\bf Proof.} This is obvious because every monomial in
$\overline{\beta}_{i,j}$ is divisible by some quadratic monomial
in $u_1,\dots,u_m$, and $\widetilde{\beta}_{i,1}\equiv 0$ for
$i=1,\dots,l$, and by the standard formulas, transforming regular
functions under a blow up. Q.E.D. for the lemma.

The lemma gives us the inequality
$$
\mathop{\rm rk}(\widetilde{\beta}_{i,2},\, 1\leqslant i\leqslant
l)\geqslant\mathop{\rm rk}(\overline{\beta}_{i,2},\, 1\leqslant
i\leqslant l)\geqslant r_k.
$$
Setting $T_p{\cal X}_B=\{\widetilde{\beta}_{i,1}=0\,|\,
l+1\leqslant j\leqslant k\}$ and using Remark 1.4, we get
$$
\mathop{\rm rk}(\widetilde{\beta}_{i,2}|_{T_p{\cal X}_B},\,
1\leqslant i\leqslant l)\geqslant r_k-2(k-l)\geqslant r_l.
$$
Therefore, $p\in{\cal X}_B$ is a multi-quadratic singularity of
type $2^l$ and rank $\geqslant r_l$. Q.E.D. for Theorem 4.1.

{\bf Corollary 4.2.} {\it Assume that $\cal X$ has multi-quadratic
singularities of type $\underline{r}$, where $r_l\geqslant 3l+1$
for all $l=1,\dots,k$. Then the singularities of ${\cal X}$ are
terminal.}

{\bf Proof.} In the notations of the proof of Theorem 4.1 it is
sufficient to show the inequality
$$
a({\cal X}_B\cap E_B,{\cal X})\geqslant 1.
$$
Assume that a point $o\in B$ of general position is a
multi-quadratic singularity of type $2^l$. From the claim (ii) of
Corollary 4.1 we get the inequality
$$
\mathop{\rm codim}(B\subset{\cal X})\geqslant l+2,
$$
so that $\mathop{\rm codim}(B\subset{\cal Y})\geqslant k+l+2$ and
for that reason
$$
a(E_B,{\cal Y})\geqslant k+l+1.
$$
By the adjunction formula
$$
a({\cal X}_B\cap E_B,{\cal X})=a(E_B,{\cal Y})-(k-l)-2l,
$$
which implies the required inequality. Q.E.D. for the corollary.

This completes the proof of Theorem 1.1.\vspace{0.3cm}


{\bf 4.4. Singularities of complete intersections.} Let us show
Proposition 4.2. We will prove the claims (i) and (ii)
simultaneously: by Grothendieck's theorem on parafactoriality
\cite{CL,Call1994} the claim (ii) for $e=4$ implies the
factoriality of the variety $Q$.

We argue by induction on $l\geqslant 1$. For one quadric ($l=1$)
the claims (i) and (ii) are obvious. Since
$$
\mathop{\rm rk} q_{[1,l-1]}\geqslant \mathop{\rm rk} q_{[1,l]},
$$
we may assume that the claims (i) and (ii) are true for the tuple
of quadratic forms $q_1$, \dots, $q_{l-1}$. In particular, the
scheme of their common zeros $Q_{l-1}$ is an irreducible reduced
factorial complete intersection of type $2^{l-1}$ in ${\mathbb
P}^{N-1}$, so that $\mathop{\rm Pic} Q_{l-1}={\mathbb Z} H_{l-1}$,
where $H_{l-1}$ is the class of a hyperplane section: every
effective divisor on $Q_{l-1}$ is cut out on $Q_{l-1}$ by a
hypersurface in ${\mathbb P}^{N-1}$.

The scheme of common zeros of the quadratic forms $q_1$, \dots,
$q_l$ is the divisor of zeros of the form $q_l$ on the variety
$Q_{l-1}$. This divisor is reducible or non-reduced if and only if
there is a form $q^*_l$ of rank $\leqslant 2$ such that
$$
q_l-q^*_l\in\langle q_1,\dots,q_{l-1}\rangle,
$$
and in that case $\mathop{\rm rk}q_{[1,l]}\leqslant 2$, which
contradicts the assumption. Therefore, $Q$ is an irreducible
reduced complete intersection. It is easy to see that
$Q\subset{\mathbb P}^{N-1}$ is non-degenerate. Since
$$
\mathop{\rm rk}q_{[1,l-1]}\geqslant 2(l-1)+(e+2)-1
$$
(for the claim (i) we set $e=4$), we have
$$
\mathop{\rm codim}(\mathop{\rm Sing}Q_{l-1}\subset
Q_{l-1})\geqslant e+2,
$$
so that
$$
\mathop{\rm codim}((Q\cap\mathop{\rm Sing}Q_{l-1})\subset
Q)\geqslant e+1.
$$
It is easy to see that a point $p\in Q$, which is non-singular on
$Q_{l-1}$, is singular on $Q$ if and only if for some
$\lambda_1,\dots,\lambda_{l-1}$ the quadric
$$
q_l-\lambda_1q_1-\dots-\lambda_{l-1}q_{l-1}=0
$$
is singular at that point. Since the singular set of a quadric of
rank $r$ in ${\mathbb P}^{N-1}$ has dimension $N-1-r$, we conclude
that the dimension of the set
$$
\mathop{\rm Sing}Q\cap(Q_{l-1}\setminus\mathop{\rm Sing}Q_{l-1})
$$
does not exceed $N-1-\mathop{\rm rk}q_{[1,l]}+(l-1)$, whence it
follows that the codimension of that set with respect to $Q$ is at
least $\mathop{\rm rk}q_{[1,l]}-2l+1\geqslant e$. Q.E.D. for
Proposition 4.2.\vspace{0.3cm}


{\bf 4.5. Linear subspaces and projections.} Now let us consider
the questions that are naturally close to Proposition 4.2 and its
proof. These questions are of key importance in the proof of
Theorem 3.3 (which will be given in \S 5). Since in Theorem 3.3
the multi-quadratic singularity is of type $2^k$, starting from
this moment we consider $k$ quadratic forms $q_1$, \dots, $q_k$ in
$N$ variables (that is, on ${\mathbb P}^{N-1}$), and the tuple of
them is denoted by the symbol $q_{[1,k]}$. The symbol $Q$, as
above, stands for the complete intersection of these $k$ quadrics
$\{q_i=0\}$ in ${\mathbb P}^{N-1}$.

{\bf Proposition 4.3.} {\it Assume that for some $b\geqslant 0$
the inequality
$$
\mathop{\rm rk} q_{[1,k]}\geqslant 2(1+b)k+3
$$
holds. Then for every point $p\in Q\setminus \mathop{\rm Sing} Q$
there is a linear space $\Pi\subset {\mathbb P}^{N-1}$ of
dimension $b$, such that $p\in\Pi\subset Q$, and moreover,}
$\Pi\cap\mathop{\rm Sing} Q=\emptyset$.

{\bf Proof} contains the (obvious) construction of such linear
subspaces. We argue by induction on $b$. If $b=0$, then $\Pi$ is
the point $p$ itself and there is nothing to prove. Assume that
$b\geqslant 1$ and for $b-1$ the claim of the Proposition is true.

Consider the linear subspace $T=T_pQ$ of codimension $k$ in
${\mathbb P}^{N-1}$. Obviously, every linear space in ${\mathbb
P}^{N-1}$ that contains the point $p$ and is contained in $Q$, is
contained in $T$, too. Furthermore, $Q\cap T$ is defined by the
quadratic forms $q_1|_T,\dots,q_k|_T$. Since $\mathop{\rm rk}
q_{[1,k]}|_T\geqslant\mathop{\rm rk} q_{[1,k]}-2k$, the inequality
$$
\mathop{\rm rk} q_{[1,k]}|_T\geqslant 2bk+3
$$
holds, where every quadric $\{q_i|_T=0\}$, $i=1,\dots,k$, is by
construction a cone with the vertex at $p$. Therefore, $Q\cap T$
is a cone with the vertex at the point $p$. Let $P\subset T$ be a
hyperplane in $T$ that does not contain the point $p$. Then the
cone $Q\cap T$ is a cone is the cone with the base $Q\cap P$,
where $Q\cap P$ is a complete intersection of the quadrics
$\{q_i|_P=0\}$, where, obviously,
$$
\mathop{\rm rk} q_{[1,k]}|_P=\mathop{\rm rk} q_{[1,k]}|_T\geqslant
2(1+(b-1))k+3.
$$
By the induction hypothesis, there is a linear subspace
$\Pi^{\sharp}\subset P$ of dimension $(b-1)$, such that
$\Pi^{\sharp}\subset Q\cap P$ and $\Pi^{\sharp}\cap\mathop{\rm
Sing}(Q\cap P)=\emptyset$.

Furthermore, the set of singular points $\mathop{\rm Sing}(Q\cap
T)$ is a cone with the vertex $p$, the base of which is
$\mathop{\rm Sing}(Q\cap P)$, so that for the subspace
$\Pi=\langle p,\Pi^{\sharp}\rangle$, which is a cone with the
vertex $p$ and the base $\Pi^{\sharp}$, we have
$\Pi\cap\mathop{\rm Sing}(Q\cap T)=\{p\}$. Since $T\cap\mathop{\rm
Sing}Q\subset\mathop{\rm Sing}(Q\cap T)$ and $p\not\in\mathop{\rm
Sing}Q$, we get $\Pi\cap\mathop{\rm Sing}Q=\emptyset$, which
completes the proof of the proposition.

{\bf Proposition 4.4.} {\it Let $b\geqslant \beta\geqslant 0$ be
some integers. Assume that the inequality
$$
\mathop{\rm rk} q_{[1,k]}\geqslant 2k(b+\beta+1)+2\beta+3
$$
holds. Then for every linear subspace $P\subset{\mathbb P}^{N-1}$
of codimension $\beta$ and a general linear subspace $\Pi\subset
Q$, $\Pi\cap\mathop{\rm Sing}=\emptyset$, of dimension $b$ the
intersection $P\cap\Pi$ has codimension $\beta$ in $\Pi$.}

{\bf Proof.} Again we argue by induction on $b,\beta$; the case
$\beta=0$ is trivial, only the equality $\Pi\cap\mathop{\rm Sing}
Q=\emptyset$ for a general subspace $\Pi\subset Q$ of dimension
$\beta\geqslant 0$ is needed, and it is true by Proposition 4.3.

Let us show our claim in the assumption that it is true for
$\beta-1$.

First of all, note that
$$
\mathop{\rm rk} q_{[1,k]}|_P\geqslant 2k(b+\beta+1)+3>2k+3,
$$
so that by Proposition 4.2 the intersection $Q\cap P$ is an
irreducible reduced complete intersection of type $2^k$ in $P$; in
particular, a point of general position $p\in Q\cap P$ is
non-singular. This means that
$$
T_p(Q\cap P)=T_pQ\cap P
$$
is of codimension $k$ in $P$, so that $T_pQ$ and $P$ are in
general position. The property to be in general position is an
open property, therefore for a point of general position $p\in Q$
(in particular, $p\not\in P$) the linear subspaces $T_pQ$ and $P$
are in general position and their intersection $T_pQ\cap P$ is of
codimension $k$ in $P$ and of codimension $\beta$ in $T_pQ$.

Consider a general hypersurface $Z$ in $T_pQ$, containing the
subspace $T_pQ\cap P$ and not containing the point $p$. We have
$$
\mathop{\rm rk} q_{[1,k]}|_Z\geqslant
2k(b+\beta)+2\beta+1=2k(b+(\beta-1)+1)+2(\beta-1))+3,
$$
so that by the induction hypothesis for a general linear subspace
$\Pi^{\sharp}\subset Q\cap Z$ of dimension $(b-1)$ that does not
meet the set $\mathop{\rm Sing}(Q\cap Z)$ the intersection
$$
(P\cap T_pQ)\cap\Pi^{\sharp}=P\cap\Pi^{\sharp}
$$
is of codimension $\beta-1=\mathop{\rm codim}((P\cap T_pQ)\subset
Z)$ with respect to $\Pi^{\sharp}$.

Then the linear space
$$
\Pi=\langle p,\Pi^{\sharp}\rangle\subset T_pQ
$$
of dimension $b$ is contained in $Q$, does not meet the set
$\mathop{\rm Sing}Q$ (see the proof of Proposition 4.3) and,
finally, the subspace
$$
P\cap\Pi=P\cap T_pQ\cap\Pi=P\cap T_pQ\cap
Z\cap\Pi=P\cap\Pi^{\sharp}
$$
is of codimension $\beta$ with respect to $\Pi$. Q.E.D. for the
proposition.

{\bf Corollary 4.3.} {\it In the assumptions of Proposition 4.4,
where $\beta\geqslant k$, let $Y\subset Q$ be an irreducible
subvariety of codimension $\beta-k$. Then the restriction onto $Y$
of the projection
$$
\mathop{\rm pr}\nolimits_{\Pi}\colon {\mathbb
P}^{N-1}\dashrightarrow {\mathbb P}^{N-b-2}
$$
from a general subspace $\Pi\subset Q$ of dimension $b$ is
dominant.}

{\bf Proof.} Let $p\in Y$ be a non-singular point. We apply
Proposition 4.4 to the subspace $P=T_pY\subset{\mathbb P^{N-1}}$
of codimension $\beta$. A general subspace $\Pi\subset Q$ of
dimension $b\geqslant\beta$ does not contain the point $p$ and is
in general position with $P$, so that $\mathop{\rm pr}_{\Pi}|_P$
is regular in a neighborhood of the point $p$ and its differential
at the point $p$ is an epimorphism. Therefore, $\mathop{\rm
pr}_{\Pi}|_Y$ is regular at the point $p$ and its differential at
that point is an epimorphism. Q.E.D. for the corollary.

Note an important particular case.

{\bf Corollary 4.4.} {\it Assume that $b\geqslant k$ and the
inequality
$$
\mathop{\rm rk}\nolimits_{[1,k]}\geqslant 2k(b+k+2)+3
$$
holds. Then the restriction of the projection $\mathop{\rm
pr}_{\Pi}$ from a general subspace $\Pi\subset Q$ of dimension $b$
onto $Q$ is dominant and its general fibre is a linear subspace of
dimension} $b+1-k$.

{\bf Proof.} That it is dominant, follows from the previous
corollary, so that the dimension of a general fibre is $b+1-k$.
Furthermore, $\mathop{\rm pr}_{\Pi}$ fibres ${\mathbb P}^{N-1}$
(more precisely, ${\mathbb P}^{N-1}\setminus\Pi$) into linear
subspaces $\Pi^{\sharp}\supset\Pi$ of dimension $b+1$. The centre
$\Pi$ of the projection is a hyperplane in $\Pi^{\sharp}$. Since
$\Pi\subset Q$, the quadric $\{q_i|_{\Pi^{\sharp}}=0\}$ is the
union of two hyperplanes, one of which is $\Pi$. Now the claim of
the corollary is obvious.

Let $\Pi\subset Q$ be a linear subspace of dimension $b\geqslant
k$, not meeting the set $\mathop{\rm Sing}Q$, and
$\sigma\colon\widetilde{Q}\to Q$ and $\sigma_{\mathbb
P}\colon\widetilde{{\mathbb P}^{N-1}}\to{\mathbb P}^{N-1}$ the
blow ups of $\Pi$ on $Q$ and ${\mathbb P}^{N-1}$, respectively, so
that we can identify $\widetilde{Q}$ with the strict transform of
$Q$ on $\widetilde{{\mathbb P}^{N-1}}$. By the symbols $E_Q$ and
$E_{\mathbb P}$ we denote the exceptional divisors of these blow
ups; we consider $E_Q$ as a subvariety in $E_{\mathbb P}$. Let
$\varphi\colon\widetilde{Q}\to{\mathbb P}^{N-b-2}$ and
$\varphi_{\mathbb P}\colon\widetilde{{\mathbb P}^{N-1}}
\to{\mathbb P}^{N-b-2}$ be the regularizations of the rational
maps $\mathop{\rm pr}_{\Pi}|_Q$ and $\mathop{\rm pr}_{\Pi}$,
respectively. We have the natural identification $E_{\mathbb
P}=\Pi\times{\mathbb P}^{N-b-2}$, where the map
$$
\varphi_{\mathbb P}|_{E_{\mathbb P}}\colon E_{\mathbb
P}\to{\mathbb P}^{N-b-2}
$$
is the projection onto the second factor. In the assumptions of
Corollary 4.4 the morphism $\varphi$ is surjective and for a point
of general position $p\in{\mathbb P}^{N-b-2}$ the fibre
$\varphi^{-1}(p)$ is a linear subspace of dimension $b+1-k$ in
$\varphi^{-1}_{\mathbb P}(p)\cong{\mathbb P}^{b+1}$, which is not
contained entirely in the hyperplane
$$
\varphi^{-1}_{\mathbb P}(p)\cap E_{\mathbb
P}=\left(\varphi_{\mathbb P}|_{E_{\mathbb P}}\right)^{-1}(p),
$$
which identifies naturally with $\Pi$, and for that reason
$\varphi^{-1}(p)\cap E_{\mathbb P}$ identifies naturally with a
subspace of dimension $b-k$ in $\Pi$ (and a hyperplane in
$\varphi^{-1}(p)$). However,
$$
\varphi^{-1}(p)\cap E_{\mathbb P}=\varphi^{-1}(p)\cap
E_Q=(\varphi|_{E_Q})^{-1}(p),
$$
so that arguing by dimensions, we conclude that the restriction
$\varphi|_{E_Q}$ is surjective.

{\bf Proposition 4.5.} {\it In the assumptions of Corollary 4.4,
let $Y\subset\Pi$ be an irreducible closed subset, and assume that
$$
b\geqslant k+\mathop{\rm codim} (Y\subset\Pi).
$$
Then the restriction $\varphi|_{\sigma^{-1}(Y)}$ is surjective, so
that for a point of general position $p\in{\mathbb P}^{N-b-2}$ the
intersection $\varphi^{-1}(p)\cap \sigma^{-1}(Y)$ is non-empty and
each of its components is of codimension $\mathop{\rm codim}
(Y\subset\Pi)$ in the projective space} $\varphi^{-1}(p)\cap
E_{\mathbb P}$.

{\bf Proof.} Obviously,
$$
\sigma^{-1}(Y)=\sigma^{-1}_{\mathbb P}(Y)\cap\widetilde{Q}=
\sigma^{-1}_{\mathbb P}(Y)\cap E_Q.
$$
Since $\varphi^{-1}(p)\subset\widetilde{Q}$, the equality
$$
\varphi^{-1}(p)\cap\sigma^{-1}(Y)=\varphi^{-1}(p)\cap\sigma^{-1}_{\mathbb
P}(Y)
$$
holds, but $\sigma^{-1}(Y)=Y\times{\mathbb P}^{N-b-2}$ in terms of
the direct decomposition of the exceptional divisor $E_{\mathbb
P}$. Therefore, identifying the fibre of the projection
$\varphi_{\mathbb P}|_{E_{\mathbb P}}$ with the projective space
$\Pi$, we get that the intersection
$\varphi^{-1}(p)\cap\sigma^{-1}(Y)$ identifies naturally with the
intersection of $Y$ and the linear subspace $\varphi^{-1}(p)\cap
E_{\mathbb P}$ of dimension $b-k$ in $\Pi$. By our assumption this
intersection is non-empty, so that the morphism
$\varphi|_{\sigma^{-1}(Y)}$ is surjective. Q.E.D. for the
proposition.


\section{The special hyperplane section}

In this section we prove Theorem 3.3.\vspace{0.3cm}

{\bf 5.1. Start of the proof.} We use the notations of Subsection
1.7 and the assumptions of Theorem 3.3. Recall that
$$
I_X=[2k+3,k+c_X-1]\cap {\mathbb Z}
$$
is the set of admissible dimensions for the working triple
$(X,D,o)$. Consider a general subspace $P\ni o$ in ${\mathbb
P}(X)$ of the minimal admissible dimension $2k+3$. Since
$a(E_{X\cap P})=2$ and $\nu(D)\leqslant\frac32n(D)<2n(D)$, we
conclude that the pair
$$
\left((X\cap P)^+,\frac{1}{n(D)} D|^+_{X\cap P}\right)
$$
is not log canonical, but canonical outside the exceptional
divisor $E_{X\cap P}$. By the inequality $\nu(D)<2n(D)$ we can
apply the connectedness principle to this pair:
\begin{equation}\label{16.01.23.1}
\mathop{\rm LCS}\left((X\cap P)^+,\frac{1}{n(D)}D|^+_{X\cap
P}\right)
\end{equation}
is a proper connected closed subset of the exceptional divisor
$E_{X\cap P}$. There are the following options:

$(1)_P$ this subset contains a divisor,

$(2)_P$ some irreducible component of maximal dimension
$B(P)\subset E_{X\cap P}$ in this set has a positive dimension and
codimension $\geqslant 2$ in $E_{X\cap P}$,

$(3)_P$ this subset is a point.

{\bf Remark 5.1.} In the case $(1)_P$ the divisor in the subset
(\ref{16.01.23.1}) is unique and is a hyperplane section of the
variety $E_{X\cap P}\subset{\mathbb E}_{X\cap P}$, since
$D|^+_{X\cap P}$ has along this subvariety the multiplicity
$>n(D)$ (since it is the centre of some non log canonical
singularity), whereas the restriction $D^+|_{E_{X\cap P}}$ is cut
out on $E_{X\cap P}$ by a hypersurface of degree $\nu(D)<2n(D)$.

Since $P\ni o$ is a subspace of general position, we go back to
the original variety $X$ and get that the pair
$(X^+,\frac{1}{n(D)}D|^+)$ is not log canonical, and moreover, for
the centre $B\subset E_X$ of some non log canonical singularity of
that pair one of the three option takes place:

(1) $B$ is a hyperplane section of $E_X\subset{\mathbb E}_X$,

(2) $\mathop{\rm codim}(B\subset E_X)\in\{2,\dots,k+1\}$,

(3) $B$ is a linear subspace of codimension $2k+2$ in ${\mathbb
E}_X$, which is contained in $E_X$.

{\bf Proposition 5.1.} {\it The option (1) does not take place.}

{\bf Proof.} Assume the converse: $B$ is a hyperplane section of
$E_X$. Let $R\subset X$, $R\ni o$ be the uniquely determined
hyperplane section, such that $R^+\cap E_X=B$ (in other words,
${\mathbb P}(R)^+\cap{\mathbb E}_X$ is the hyperplane in ${\mathbb
E}_X$ that cuts out $B$ on $E_X$). Since $\mathop{\rm
mult}_BD^+>n(D)$, we get that for the effective divisor
$D_R=(D\circ R)$ on $R$ the inequality
$$
\nu(D_R)\geqslant\nu(D)+\mathop{\rm
mult}\nolimits_BD^+>2n(D)=2n(D_R)
$$
holds, which is impossible by Theorem 1.4. Q.E.D. for the
proposition.

{\bf Proposition 5.2.} {\it The option (3) does not take place.}

{\bf Proof.} Since $\mathop{\rm codim}(B\subset E_X)=k+2$, this is
impossible by the Lefschetz theorem (in order to apply the
Lefschetz theorem, it is sufficient to have the inequality
$\mathop{\rm codim}(\mathop{\rm Sing}E_X\subset E_X)\geqslant
2k+6$, for which by Proposition 4.2 it is sufficient to have the
inequality $\mathop{\rm rk}(o\in X)\geqslant 4k+5$; we have a much
stronger condition for the rank of the singularity). Q.E.D. for
the proposition.

Therefore, the option (2) takes place. By construction (or arguing
by dimension), $B\not\subset\mathop{\rm Sing}E_X$. Recall that
there is a non log canonical singularity of the pair
$(X^+,\frac{1}{n(D)}D^+)$, the centre of which is $B$.

Let $p\in B$ be a point of general position; in particular,
$p\not\in\mathop{\rm Sing}E_X$ and the more so
$p\not\in\mathop{\rm Sing} X^+$. Applying inversion of adjunction
in the word for word the same way as in \cite[Chapter 7,
Proposition 2.3]{Pukh13a} (that is, restricting $D^+$ onto a
general non-singular surface, containing the point $p$), we get
the alternative: either $\mathop{\rm mult}_BD^+>2n(D)$, or on the
blow up
$$
\varphi_p\colon X^{(p)}\to X^+
$$
of the point $p$ with the exceptional divisor $E(p)\subset
X^{(p)}$, $E(p)\cong{\mathbb P}^{N(X)-1}$, there is a hyperplane
$\Theta(p)\subset E(p)$ in $E(p)$, satisfying the inequality
\begin{equation}\label{22.09.22.1}
\mathop{\rm mult}\nolimits_BD^++\mathop{\rm
mult}\nolimits_{\Theta(p)}D^{(p)}>2n(D),
\end{equation}
where $D^{(p)}$ is the strict transform of the divisor $D^+$ on
$X^{(p)}$, and moreover, the hyperplane $\Theta(p)$ is uniquely
determined by the pair $(X^+,\frac{1}{n(D)}D^+)$ and varies
algebraically with the point $p\in B$.

The case when the inequality $\mathop{\rm
mult}\nolimits_BD^+>2n(D)$ holds, is excluded (with
simplifications) by the arguments, excluding the option
$(2)_{\Theta}$, given below, see Subsection 5.3, Remark 5.2.

There are two options for the hyperplane $\Theta(p)$:

$(1)_{\Theta}$  $\Theta(p)\neq {\mathbb P}(T_pE_X)$ (where we
identify $E(p)$ with the projectivization of the tangent space
$T_pX^+$), so that $\Theta(p)$ intersects ${\mathbb P}(T_pE_X)$ by
some hyperplane $\Theta_E(p)$,

$(2)_{\Theta}$ the hyperplanes $\Theta(p)$ and ${\mathbb
P}(T_pE_X)$ in $E(p)$ are equal.

Below (see Subsection 5.3, Remark 5.2) we show that the option
$(2)_{\Theta}$ does not take place: it implies that $E_X\subset
D^+$, which is impossible; the same arguments exclude the
inequality $\mathop{\rm mult}\nolimits_BD^+>2n(D)$, too.

Therefore, we may assume that the option $(1)_{\Theta}$ takes
place.\vspace{0.3cm}


{\bf 5.2. The existence of the special hyperplane section.} Adding
the upper index $(p)$ means the strict transform on $X^{(p)}$: we
used this principle for the divisor $D$ above and will use it for
other subvarieties on $X^+$. Our aim is to prove the following
claim.

{\bf Theorem 5.1.} {\it There is a hyperplane section $\Lambda$ of
the exceptional divisor $E_X\subset{\mathbb E}_X$, containing $B$,
satisfying the inequality
$$
\mathop{\rm mult}\nolimits_{\Lambda}D^+>\frac{2n(D)-\nu(D)}{k+1}.
$$
Moreover, for a point of general position $p\in B$ the following
equality holds:}
$$
\Lambda^{(p)}\cap E(p)=\Theta_E(p).
$$

{\bf Proof.} Let $L\subset E_X,L\ni p$ be a line in the projective
space ${\mathbb E}_X$, such that $L\cap\mathop{\rm
Sing}E_X=\emptyset$ and
$$
L^{(p)}\cap E_(p)\in\Theta_E(p).
$$

{\bf Lemma 5.1.}{\it The line $L$ is contained in $D^+$.}

{\bf Proof.} Assume the converse. Then $D^+|_L$ is an effective
divisor on $L$ of degree $\nu(D)\leqslant\frac32n(D)<2n(D)$. At
the same time, the divisor $D^+|_L$ contains the point $p$ with
multiplicity $>2n(D)$ due to the inequality (\ref{22.09.22.1}).
The contradiction proves the lemma. Q.E.D.

{\bf Proposition 5.3.} {\it The following inequality holds:}
$$
\mathop{\rm mult}\nolimits_LD^+>\frac{2n(D)-\nu(D)}{k+1}.
$$

{\bf Proof} is given in \S 6.

Let us go back to the proof of Theorem 5.1.

We will construct the set $\Lambda\subset E_X$ explicitly, and
then prove that it is a hyperplane section. The exceptional
divisor $E_X$ is a complete intersection of $k$ quadrics in
${\mathbb E}_X$:
$$
E_X=\{q_1=\dots q_k=0\},
$$
using the notations of Subsection 4.5. Let $U_B\subset B$ be a
non-empty Zariski open subset, where
$$
U_B\cap\mathop{\rm Sing} E_X=\emptyset
$$
and for every point $p\in B$ the option $(1)_{\Theta}$ takes
place. By the assumption on the rank of the multi-quadratic point
$o\in X$ for $p\in U_{B}$ the set $E_X\cap T_pE_X$ (where
$T_pE_X\subset{\mathbb E}_X$ is the embedded tangent space, that
is, a linear subspace of codimension $k$ in ${\mathbb E}_X$) is
irreducible and reduced, and moreover, every hyperplane section of
that set is also irreducible and reduced. Indeed, by Proposition
4.2, in order to have these properties it is sufficient to have
the inequality $\mathop{\rm rk}q_{[1,k]}\geqslant 4k+5$, because
by Remark 1.4 it implies the inequality
$$
\mathop{\rm rk}q_{[1,k]}|_{T_pE_X}\geqslant 2k+5
$$
and we can apply Proposition 4.2. Obviously, $E_X\cap T_pE_X$ is a
cone with the vertex $p$, consisting of all lines $L\subset E_X$,
$L\ni p$. The singular set of that cone is of codimension
$\geqslant 6$ (Proposition 4.2), and so for a general line $L\ni
p$
$$
L\cap\mathop{\rm Sing}(E_X\cap T_p E_X)=\{p\},
$$
so that $L\cap\mathop{\rm Sing}E_X=\emptyset$ and the same is true
for every hyperplane section of the cone $E_X\cap T_pE_X$,
containing the point $p$, since its singular set is of codimension
$\geqslant 4$ (Remark 1.4).

Let ${\cal L}(p)$ be the union of all lines $L\subset E_X,L\ni p$,
such that
$$
L^{(p)}\cap E(p)\in\Theta_E(p).
$$
Obviously, ${\cal L}(p)$ is the section of the cone $E_X\cap
T_pE_X$ by some hyperplane, containing the point $p$ (this
hyperplane corresponds to the hyperplane $\Theta_E(p$)). As we
have shown above, ${\cal L}(p)$ is an irreducible closed subset of
codimension $k+1$ in $E_X$, and
$$
\mathop{\rm mult}\nolimits_{{\cal L}(p)}
D^+>\frac{2n(D)-\nu(D)}{k+1}.
$$
Set
$$
\Lambda=\overline{\mathop{\bigcup}\limits_{p\in U_B}{\cal L}(p)}
$$
(the overline means the closure). By what was said above, the
inequality
$$
\mathop{\rm mult}\nolimits_{\Lambda}D^+>\frac{2n(D)-\nu(D)}{k+1}
$$
holds.

{\bf Theorem 5.2.} {\it The subset $\Lambda\subset E_X$ is a
hyperplane section of the variety} $E_X\subset{\mathbb E}_X$.

We will prove Theorem 5.2 in two steps: first, we will show that
$\Lambda$ is a prime divisor on $E_X$, and then, that this divisor
is a hyperplane section. By construction, the set $\Lambda$ is
irreducible.\vspace{0.3cm}


{\bf 5.3. The set $\Lambda$ is a divisor.} By our assumption about
the rank of the point $o\in X$ for $b=k+1$ the inequality
\begin{equation}\label{01.10.22.1}
\mathop{\rm rk} q_{[1,k]}\geqslant 2k(b+2k+2)+2(2k+1)+3
\end{equation}
holds. By Corollary 4.3, for a general subspace $\Pi\subset E_X$
of dimension $b$ the restriction onto $B$ of the projection
$$
\mathop{\rm pr}\nolimits_{\Pi}\colon {\mathbb
P}^{N(X)-1}\dashrightarrow {\mathbb P}^{N(X)-b-2}
$$
from the subspace $\Pi$ is dominant. Let $s\in {\mathbb
P}^{N(X)-b-2}$ be a point of general position. By the symbol
$\langle\Pi,s\rangle$ denote the closure
$$
\overline{\mathop{\rm pr}\nolimits_{\Pi}^{-1}(s)}\subset {\mathbb
P}^{N(X)-1}
$$
(this is a $(\mathop{\rm dim}\Pi+1)$-dimensional subspace) and set
$$
E_X(\Pi,s)=E_X\cap \langle\Pi,s\rangle.
$$
For the blow ups $\sigma\colon\widetilde{E}_X\to E_X$ and
$\sigma_{\mathbb P}\colon\widetilde{{\mathbb
P}^{N(X)-1}}\to{\mathbb P}^{N(X)-1}$ of the subspace $\Pi$ on
$E_X$ and ${\mathbb P}^{N(X)-1}$, respectively, let
$\varphi\colon\widetilde{E}_X\to{\mathbb P}^{N(X)-b-2}$ and
$\varphi_{\mathbb P}\colon\widetilde{{\mathbb
P}^{N(X)-1}}\to{\mathbb P}^{N(X)-b-2}$ be the regularizations of
the projections $\mathop{\rm pr}_{\Pi}|_{E_X}$ and $\mathop{\rm
pr}_{\Pi}$, respectively. Obviously, the fibre
$\varphi^{-1}_{\mathbb P}(s)$ identifies naturally with
$\langle\Pi,s\rangle$, and the fibre $\varphi^{-1}(s)$ with
$E_X(\Pi,s)$. The fibre of the surjective morphism
$\varphi|_{\sigma^{-1}(B)}$ over the point $s$ we denote by the
symbol $B(s)$; this is a possibly reducible closed subset in
$\varphi^{-1}_{\mathbb P}(s)$, each irreducible component of which
is of codimension $c_B=\mathop{\rm codim} (B\subset E_X)$ and is
not contained entirely in the hyperplane $\Pi$ (with respect to
the identification $\varphi^{-1}_{\mathbb
P}(s)=\langle\Pi,s\rangle$). Write down $B(s)$ as a union of
irreducible components:
$$
B(s)=\mathop{\bigcup}\limits_{i\in I} B_i(s),
$$
and let $p\in B_i(s)$ be a point of general position on one of
them; in particular, $p\not\in\Pi$, so that the projection
$\mathop{\rm pr}\nolimits_{\Pi}$ is regular at that point and
$p\not\in B_j(s)$ for $j\neq i$. We will consider the point $p$ as
a point of general position on $B$, which was introduced in
Subsection 5.1, and use the notations for the blow up $\varphi_p$
of this point and for objects linked to this blow up. Note that
for $b=k+1$ we have the inequality
$$
\mathop{\rm dim}B(s)=\mathop{\rm dim}B_i(s) \geqslant 1.
$$
The set of lines $L\subset E_X(\Pi,s)$, $L\ni p$, such that
$L^{(p)}\cap E(p)\in\Theta(p)$, forms a hyperplane in
$E_X(\Pi,s)$, which we denote by the symbol $\Lambda(\Pi,s,p)$. By
construction, $\Lambda(\Pi,s,p)\subset\Lambda$.

Since any non-trivial algebraic family of hyperplanes in a
projective space sweeps out that space and for a general point $s$
we have $E_X(\Pi,s)\not\subset\Lambda$ (otherwise $\Lambda=E_X$,
which is impossible), we conclude that the hyperplane
$\Lambda(\Pi,s,p)$ does not depend on the choice of a point of
general position $p\in B_i(s)$, so that
$$
\Lambda(\Pi,s,p)=\Lambda(\Pi,s,B_i(s))
$$
is a hyperplane in $\varphi^{-1}(s)=E_X(\Pi,s)$, containing the
component $B_i(s)$. Therefore, for a general point $s$ the
intersection $\Lambda\cap E_X(\Pi,s)$ contains a divisor in
$E_X(\Pi,s)$, whence we get that $\Lambda\subset E_X$ is a (prime)
divisor on $E_X$, as we claimed. This divisor is cut out on $E_X$
by a hypersurface of degree $d_{\Lambda}$ in ${\mathbb E}_X$. It
remains to show that $d_{\Lambda}=1$.

{\bf Remark 5.2.} We promised above that the option $(2)_{\Theta}$
does not take place. Indeed, if it does, then every line $L\ni p$
in $E_X(\Pi,s)$ is contained in $\Lambda$, so that
$E_X(\Pi,s)\subset\Lambda$ and for that reason
$E_X\subset\Lambda$, which is absurd. In a similar way, if
$\mathop{\rm mult}_BD^+>\nu(D)$, then every line in $E_X(\Pi,s)$,
meeting $B$, is contained in $\Lambda$, so that
$E_X\subset\Lambda$, which is impossible. Therefore, the
inequality
$$
\mathop{\rm mult}\nolimits_BD^+\leqslant\nu(D)
$$
holds. \vspace{0.3cm}


{\bf 5.4. The divisor $\Lambda$ is a hyperplane section.} Let us
consider the intersection $\Lambda\cap E_X(\Pi,s)$ for a general
point $s$ in more details. This is a possibly reducible divisor,
each component of which has multiplicity 1, containing at least
one hyperplane. If in this divisor there are components of degree
$\geqslant 2$, then the union of hyperplanes in $\Lambda(\Pi,s)$
gives a proper closed subset of $\Lambda_1(\Pi,s)$, which is also
a divisor. Then
$$
\overline{\bigcup_s\Lambda_1(\Pi,s)}
$$
(the union is taken over a non-empty open subset in ${\mathbb
P}^{N(X)-b-2}$) is a proper closed subset in $\Lambda$, which is
of codimension 1 in $E_X$, which is impossible as $\Lambda$ is a
prime divisor. We conclude that $\Lambda(\Pi,s)$ is a union of
precisely $d_{\Lambda}$ distinct hyperplanes in $E_X(\Pi,s)$.

Assume that $d_{\Lambda}\geqslant 2$. By our assumptions about the
rank $\mathop{\rm rk}(o\in X)$ the inequality (\ref{01.10.22.1})
holds for $b=3k$:
$$
\mathop{\rm rk}q_{[1,k]}\geqslant 10k^2+8k+5.
$$
Again we apply Corollary 4.3, now to a general subspace
$\Pi^*=E_X(\Pi,s)$ of dimension $b^*=b+1-k\geqslant 2k+1$. The
subspace $\Pi^*$ does not meet the set $\mathop{\rm Sing} E_X$ and
the restriction of the projection from $\Pi^*$
$$
\mathop{\rm pr}\nolimits_{\Pi^*}\colon{\mathbb
P}^{N(X)-1}\dashrightarrow {\mathbb P}^{N(X)-b^*-2}
$$
onto $B$ is dominant. Let $s^*\in{\mathbb P}^{N(X)-b^*-2}$ be a
point of general position. We use the notations introduced above
and write $E_X(\Pi^*,s^*)$. For the blow ups of the subspace
$\Pi^*$ we use the symbols $\sigma_{\Pi^*}$ and $\sigma_{{\mathbb
P},\Pi^*}$, respectively, and for the regularized projections the
symbols $\varphi_{\Pi^*}$ and $\varphi_{{\mathbb P},\Pi^*}$. The
symbol $\langle\Pi^*,s^*\rangle$ has the same meaning as above.
Set
$$
E^*=\sigma_{\Pi^*}^{-1}(\Pi^*)\quad\mbox{and}\quad E^*_{\mathbb
P}= \sigma_{{\mathbb P},\Pi^*}^{-1}(\Pi^*)
$$
to be the exceptional divisors of the blow up of $\Pi^*$ on $E_X$
and ${\mathbb P}^{N(X)-1}$. By the arguments immediately before
the statement of Proposition 4.5, the map $\varphi_{\Pi^*}|_{E^*}$
is surjective, and by Proposition 4.5 (which applies since
$b^*\geqslant k+1$) the intersection
$$
\varphi_{\Pi^*}^{-1}(s^*)\cap \sigma_{\Pi^*}^{-1}(\Lambda\cap
\Pi^*)
$$
is non-empty and each of its irreducible components is of
codimension 1 in the projective space
$\varphi_{\Pi^*}^{-1}(s^*)\cap E^*_{\mathbb P}$.

By what was shown above, $\Lambda\cap{\Pi^*}$ is a union of
$d_{\Lambda}$ distinct hyperplanes $\Lambda^*_i$, $i\in I$. In a
similar way,
$$
\Lambda\cap
E_X(\Pi^*,s^*)=\sigma^{-1}_{\Pi^*}(\Lambda)\cap\varphi^{-1}_{\Pi^*}(s^*)
$$
is the union of $d_{\Lambda}$ distinct hyperplanes in
$\varphi^{-1}_{\Pi^*}(s^*)$, none of which coincides with the
hyperplane $\varphi^{-1}_{\Pi^*}(s^*)\cap E_{\mathbb P}^*$. Note
that the strict transform of the divisor $\Lambda$ with respect to
the blow up $\sigma_{\Pi^*}$ is just its full inverse image
$\sigma^{-1}_{\Pi^*}(\Lambda)$, since $\Lambda\not\subset\Pi^*$.
Furthermore,
$$
\sigma^{-1}_{\Pi^*}(\Lambda)\cap
E^*=\sigma^{-1}_{\Pi^*}(\Lambda\cap \Pi^*)=\bigcup_{i\in
I}\sigma^{-1}_{\Pi^*}(\Lambda^*_i),
$$
and every intersection
$\varphi^{-1}_{\Pi^*}(s^*)\cap\sigma^{-1}_{\Pi^*}(\Lambda^*_i)$ is
a hyperplane in $\varphi^{-1}_{\Pi^*}(s^*)\cap E^*_{\mathbb P}$.
It follows that each irreducible component of set $\Lambda\cap
E_X(\Pi^*,s^*)$ intersects the hyperplane
$\varphi^{-1}_{\Pi^*}(s^*)\cap E^*_{\mathbb P}$ by one of the
hyperplanes
$\sigma^{-1}_{\Pi^*}(\Lambda^*_i)\cap\varphi^{-1}_{\Pi^*}(s^*)$,
$i\in I$. Thus one can write down
$$
\Lambda\cap E_X(\Pi^*,s^*)=\bigcup_{i\in I}\Lambda_i(\Pi^*,s^*),
$$
where $\Lambda_i(\Pi^*,s^*)$ is a hyperplane in $E_X(\Pi^*,s^*)$,
satisfying the equality
$$
\Lambda_i(\Pi^*,s^*)\cap
E^*=\sigma^{-1}_{\Pi^*}(\Lambda^*_i)\cap\varphi^{-1}_{\Pi^*}(s^*).
$$
In other words, the choice of a component of the intersection
$\Lambda\cap\Pi^*$ determines uniquely the component of the
intersection of $\Lambda$ with
$\langle\Pi^*,s^*\rangle=\varphi^{-1}_{\Pi^*}(s^*)=E_X(\Pi^*,s^*)$
for a general point $s^*$. Now set
$$
\Lambda_i=\sigma_{\Pi^*}\left(\overline{\mathop{\bigcup}\limits_{s^*}
(\Pi^*,s^*)}\right),
$$
where the union is taken over a non-empty Zariski open subset of
the projective space ${\mathbb P}^{N(X)-b^*-2}$. This is a prime
divisor on $E_X$, and moreover, $\Lambda_i\subset\Lambda$ and for
that reason $\Lambda_i=\Lambda$, whence we conclude that all
hyperplanes $\Lambda_i(\Pi^*,s^*)$ are the same, which is a
contradiction with the assumption that $d_{\Lambda}\geqslant 2$.

Thus $d_{\Lambda}=1$ and $\Lambda$ is a hyperplane section of
$E_X\subset {\mathbb E}_X$. Q.E.D. for Theorem 5.2 and therefore,
for Theorem 5.1.\vspace{0.3cm}


{\bf 5.5. The construction of a new working triple.} Now we can
complete the proof Theorem 3.3 and construct the new working
triple $(R,D_R,o)$. Let $R\ni o$ the section of $X$ by the
hyperplane ${\mathbb P}(R)=\langle R\rangle$, such that
$$
R^+\cap {\mathbb E}_X=R^+\cap E_X=\Lambda
$$
(in other words, the hyperplane ${\mathbb P}(R)^+ \cap {\mathbb
E}_X$ cuts out $\Lambda$ on $E_X$). Since $R$ is not a component
of the effective divisor $D_X$, the scheme-theoretic intersection
$(R\circ D_X)$ is well defined, and we treat this intersection as
an effective divisor on $R$. Set $D_R=(R\circ D_X)$ in that sense.

On $R^+\subset X^+$ with the exceptional divisor
$$
E_R=(R^+\cap E_X)=\Lambda\subset {\mathbb E}_R= {\mathbb
P}(R)^+\cap{\mathbb E}_X
$$
we have the equivalence
$$
D^+_R\sim n(D_R) H_R-\nu(D_R) E_R,
$$
where $H_R$ is the class of a hyperplane section of $R$ and
$$
\nu(D_R)\geqslant\nu(D_X)+\mathop{\rm mult}\nolimits_\Lambda
D^+_X>\nu(D_X)+\frac{2n(D_X)-\nu(D_X)}{k+1}.
$$
Again $[R,o]$ is a marked complete intersection, of level
$(k,c_R)$, where $c_R=c_X-2$, $(R,D_R,o)$ is a working triple, and
the inequality
$$
2n(D_R)-\nu(D_R)<\left(1-\frac{1}{k+1}\right)(2n(D_X)-\nu(D_X))
$$
holds (since $n(D_R)=n(D_X)$).

The procedure of constructing the special hyperplane section is
complete. Q.E.D. for Theorem 3.3.


\section{Multiplicity of a line}

In this section we prove Proposition 5.3.\vspace{0.3cm}

{\bf 6.1. Blowing up a point and a curve.} Since we completed our
study of working triples, the symbol $X$ is now free and will mean
an arbitrary non-singular quasi-projective variety of dimension
$\geqslant 3$. Let $C\subset X$ be a non-singular projective
curve, $p\in C$ a point. Furthermore, let
$$
\sigma_C\colon X(C)\to X
$$
be the blow up of the curve $C$ with the exceptional divisor $E_C$
and $\sigma^{-1}_C(p)\cong {\mathbb P}^{\dim X-2}$ the fibre over
the point $p$. Let
$$
\sigma\colon X(C,\sigma^{-1}_C(p))\to X(C)
$$
be the blow up of that fibre with the exceptional divisor $E$ and
$E_C^{(p)}$ the strict transform of $E_C$ on that blow up.

On the other hand, consider the blow up
$$
\varphi_p\colon X(p)\to X
$$
of the point $p$ with the exceptional divisor $E_p$ and denote by
the symbol $C(p)$ the strict transform of the curve $C$ on $X(p)$.
Finally, let
$$
\varphi\colon X(p,C(p))\to X(p)
$$
be the blow up of the curve $C(p)$, $E_{C(p)}$ the exceptional
divisor of that blow up and $E^C_p$ the strict transform of $E_p$.

{\bf Proposition 6.1.} {\it The identity map $\mathop{\rm id}_X$
extends to an isomorphism
$$
X(C,\sigma^{-1}_C(p))\cong X(p,C(p)),
$$
identifying the subvarieties $E$ and $E^C_p$ and the subvarieties
$E^p_C$ and $E_{C(p)}$.}

{\bf Proof.} This is a well known fact, which can be checked by
elementary computations in local parameters. Q.E.D. for the
proposition.

Taking into account the identifications above, we will use the
notations $E^C_p$ and $E_{C(p)}$, and forget about $E$ and
$E^p_C$. The variety $X(C,\sigma^{-1}_C(p))$ will be denoted by
the symbol $\widetilde{X}$. Let $D$ be an effective divisor on
$X$. The symbols $D^C$ and $D^p$ stand for its strict transforms
on $X(C)$ and $X(p)$, respectively, and the symbol $\widetilde{D}$
for its strict transform on $\widetilde{X}$. Set
$$
\mu=\mathop{\rm
mult}\nolimits_CD\quad\mbox{and}\quad\mu_p=\mathop{\rm
mult}\nolimits_pD,
$$
where, of course, $\mu_p\geqslant\mu$.

{\bf Lemma 6.1.}  {\it The following equality holds:}
$$
\mathop{\rm mult}\nolimits_{\sigma^{-1}_C(p)}D^C=\mu_p-\mu.
$$

{\bf Proof.} (This is a well known fact, and we give a proof for
the convenience of the reader, and also because a similar argument
is used below.) We have the sequence of obvious equalities:
$$
\sigma^*_CD=D^C+\mu E_C,
$$
so that
$$
\sigma^*\sigma^*_CD=\widetilde{D}+\mu E_{C(p)}+(\mu+\mathop{\rm
mult}\nolimits_{\sigma^{-1}_{C(p)}}D^C)E^C_p.
$$
Considering the second sequence of blow ups, we get
$$
\varphi^*_pD=D^p+\mu_pE_p
$$
and, respectively,
$$
\varphi^*\varphi^*_pD=\widetilde{D}+\mu E_{C(p)}+\mu_pE^C_p.
$$
Comparing the two presentations of the same effective divisor, we
get the claim of the lemma.\vspace{0.3cm}


{\bf 6.2. Blowing up two points and a curve.} In the notations of
the previous subsection, let us consider the point
$$
q=C(p)\cap E_p.
$$
Set $\mu_q=\mathop{\rm mult}_qD^p$. Obviously,
$$
\mu_q\geqslant\mathop{\rm mult}\nolimits_{C(p)}D^p=\mathop{\rm
mult}\nolimits_CD=\mu.
$$
Let
$$
\varphi_q\colon X(p,q)\to X(p)
$$
the blow up of the point $q$ with the exceptional divisor $E_q$
and $C(p,q)$ the strict transform of the curve $C(p)$. Finally,
let
$$
\varphi_{\sharp}\colon X^{\sharp}\to X(p,q)
$$
be the blow up of the curve $C(p,q)$ with the exceptional divisor
$E_{C(p,q)}$ and $E^{\sharp}_q$ the strict transform of $E_q$.
Note that the curve $C(p)$ intersects $E_p$ transversally and
therefore $C(p,q)$ does not meet the strict transform $E^q_p$ of
the divisor $E_p\subset X(p)$ on $X(p,q)$.

{\bf Proposition 6.2.} {\it The restriction of the divisor $D^C$
onto the exceptional divisor $E_C$ contains the fibre
$\sigma^{-1}_C(p)$ with multiplicity at least} $\mu_p+\mu_q-2\mu$.

{\bf Proof.} Obviously, on $X^{\sharp}$ we have the equality
$$
\varphi^*_{\sharp}\varphi^*_q\varphi^*_p D=D^{\sharp}+
(\mu_q+\mu_p)E^{\sharp}_q+\mu_p E^q_p+\mu E_{C(p,q)},
$$
where $D^{\sharp}$ is the strict transform of $D$ on $X^{\sharp}$.
On the other hand, using the constructions of Subsection 6.1, we
see that $X^{\sharp}$ can be obtained as the blow up of the curve
$C(p)$ on $X(p)$ with the subsequent blowing up of the fibre of
the exceptional divisor $E_{C(p)}$ over the point $q$ or, applying
the construction of Subsection 6.1 twice, as the blow up of the
curve $C$ on $X$ with the subsequent blowing up of the fibre
$\sigma^{-1}_C(p)$ and then the blowing up of the subvariety
$$
E^C_p\cap E^p_C.
$$
In the last presentation the three prime exceptional divisors are
$$
E^{\sharp}_C=E_{C(p,q)},\quad E^{\sharp}_p\quad\mbox{and}\quad
E^{\sharp}_q.
$$
We denote the blow up $X^{\sharp}\to\widetilde{X}$ of the
subvariety $E^C_p\cap E^p_C$, mentioned above, by the symbol
$\sigma_{\sharp}$. Thus we obtain the following commutative
diagram of birational morphisms:
$$
\begin{array}{ccccc}
X(C) & \stackrel{\sigma}{\leftarrow} & \widetilde{X} &
\stackrel{\sigma_{\sharp}}{\leftarrow} &  X^{\sharp}\\
\downarrow & & \downarrow & & \downarrow\\
X &\stackrel{\varphi_p}{\leftarrow} & X(p)
&\stackrel{\varphi_q}{\leftarrow} & X(p,q),\\
\end{array}
$$
where the vertical arrows (from the left to the right) are
$\sigma_C$, $\varphi$ and $\varphi_{\sharp}$, respectively. We
have the equality
$$
\sigma^*_{\sharp}\sigma^*\sigma^*_CD=\varphi^*_{
\sharp}\varphi^*_q\varphi^*_pD.
$$
This pull back can be written down as
$$
D^{\sharp}+\mu
E^{\sharp}_C+\mu_pE^{\sharp}_p+(\mu_p+\mu_q)E^{\sharp}_q,
$$
since $E_{C(p,q)}=E^{\sharp}_C$ and $E^{\sharp}_q$ is the
exceptional divisor of the blow up $\sigma_{\sharp}$. On the other
hand, $D^C=\sigma^{*}_C D-\mu E_C$ and, besides,
$$
\sigma^*_{\sharp}\sigma^*E_C=E^{\sharp}_C+E^{\sharp}_p+2E^{\sharp}_q,
$$
so that the exceptional divisor $E^{\sharp}_q$ comes into the pull
back of the divisor $D^C$ on $X^{\sharp}$ with multiplicity
$(\mu_p+\mu_q-2\mu)$. However, the blow ups $\sigma$ and
$\sigma_{\sharp}$ do not change the divisor $E_C$, as they blow up
subvarieties of codimension 1 on the variety:
$$
\sigma\circ\sigma_{\sharp}|_{E^{\sharp}_C}\colon E^{\sharp}_C\to
E_C
$$
is an isomorphism. Since the restriction $D^{\sharp}$ onto
$E^{\sharp}_C$ is an effective divisor, it follows from here that
the restriction of the divisor $D^C$ onto $E_C$ contains the fibre
$\sigma^{-1}(p)$ (which is precisely the restriction of
$E^{\sharp}_q$ onto $E^{\sharp}_C$ in terms of the isomorphism
between $E^{\sharp}_C$ and $E_C$, discussed above) with
multiplicity $\geqslant\mu_p+\mu_q-2\mu$. Proof of Proposition 6.2
is complete.\vspace{0.3cm}


{\bf 6.3. The multiplicity of an infinitely near line.} Let us
come back to the proof of Proposition 5.3. We will obtain its
claim from a more general fact. Let $o\in{\cal X}$ be a germ of a
multi-quadratic singularity of type $2^k$, where ${\cal
X}\subset{\cal Y}$, ${\cal Y}$ is non-singular, $\mathop{\rm
codim}({\cal X}\subset{\cal Y})=k$ and the inequality
$$
\mathop{\rm rk}(o\in{\cal X})\geqslant 2k+3
$$
holds, so that $\mathop{\rm codim}(\mathop{\rm Sing{\cal
X}\subset{\cal X}})\geqslant 4$ and ${\cal X}$ is factorial. Let
$\sigma_{\cal Y}\colon{\cal Y}^+\to{\cal Y}$ be the blow up of the
point $o$ with the exceptional divisor $E_{\cal Y}$, ${\cal
X}^+\subset{\cal Y}^+$ the strict transform, so that
$$
\sigma=\sigma_{\cal Y}|_{{\cal X}^+}\colon{\cal X}^+\to{\cal X}
$$
is the blow up of the point $o$ on ${\cal X}$ with the exceptional
divisor $E_{\cal X}$, which is a complete intersection of $k$
quadrics in $E_{\cal Y}\cong{\mathbb P}^{\rm dim{\cal Y}-1}$. By
Proposition 4.2
$$
\mathop{\rm codim}(\mathop{\rm Sing} E_{\cal X}\subset E_{\cal
X})\geqslant 4.
$$
Let $L\subset E_{\cal X}$ be a line, where $L\cap \mathop{\rm
Sing} E_{\cal X}=\emptyset$ and $p\in L$ a point. Let us blow up
this point on ${\cal Y}^+$ and ${\cal X}^+$, respectively:
$$
\sigma_{p,{\cal Y}}\colon {\cal Y}_p\to {\cal
Y}^+\quad\mbox{and}\quad \sigma_p\colon {\cal X}_p\to {\cal X}^+
$$
are these blow ups with the exceptional divisors $E_{p,{\cal Y}}$
and $E_p$. Set
$$
q=L^{(p)}\cap E_p,
$$
where $L^{(p)}\subset {\cal X}_p$ is the strict transform.

Let $D_{\cal X}$ be an effective divisor on ${\cal X}$. For its
strict transform on ${\cal X}^+$ we have the equality
$$
D^+_{\cal X}=\sigma^* D_{\cal X}-\nu E_{\cal X}
$$
for some $\nu\in{\mathbb Z}_+$. Furthermore, we denote the strict
transform of $D^+_{\cal X}$ on ${\cal X}_p$ by the symbol
$D^{(p)}_{\cal X}$ and set
$$
\mu_p=\mathop{\rm mult}\nolimits_pD^+_{\cal X}\quad\mbox{and}
\quad\mu_q=\mathop{\rm mult}\nolimits_qD^{(p)}_{\cal X}.
$$
Set also $\mu=\mathop{\rm mult}_LD^+_{\cal X}$; obviously,
$\mu\leqslant\mu_p$.

{\bf Theorem 6.1.} {\it The following inequality holds:}
$$
\mu\geqslant\frac{1}{k+1}(\mu_p+\mu_q-\nu).
$$

{\bf Proof.} Let $P\subset E_{\cal Y}$ be a general linear
subspace of dimension $(k+2)$, containing the line $L$.

{\bf Lemma 6.2.} {\it The surface $S=P\cap{\cal X}^+=P\cap E_{\cal
X}$ is non-singular.}

{\bf Proof.} We argue by induction on $\mathop{\rm dim}E_{\cal
X}\geqslant 2$. If $\mathop{\rm dim}E_{\cal X}=2$, then there is
nothing to prove. Let $\mathop{\rm dim}E_{\cal X}\geqslant 3$. The
hyperplanes in $E_{\cal Y}$, tangent to $E_{\cal X}$ at at least
one point of the line $L$, form a $k$-dimensional family. The
hyperplanes, containing the line $L$, form a family (a linear
subspace) of codimension 2 in the dual projective space for
$E_{\cal Y}$, that is, of dimension $k+\mathop{\rm dim}E_{\cal
X}-2\geqslant k+1$, so that for a general hyperplane $R_{\cal
Y}\supset L$ in $E_{\cal Y}$ we have: $E_{\cal X}\cap R_{\cal Y}$
is non-singular along $L$ (and, of course, for the codimension of
the singular set we have the equality $\mathop{\rm
codim}(\mathop{\rm Sing}(E_{\cal X}\cap R_{\cal Y})\subset(E_{\cal
X}\cap R_{\cal Y}))=\mathop{\rm codim}(\mathop{\rm Sing}E_{\cal
X}\subset E_{\cal X}$)). Applying the induction hypothesis, we
complete the proof of the lemma. Q.E.D.

Let ${\cal Z}\subset{\cal Y}$, ${\cal Z}\ni o$, be a general
subvariety of dimension $(k+3)$, non-singular at the point $o$,
such that ${\cal Z}^+\cap E_{\cal Y}=P$, and
$$
{\cal X}_P={\cal X}\cap{\cal Z}
$$
(the notation ${\cal X}_P$ is chosen for convenience: ${\cal X}_P$
is determined by ${\cal Z}$, not by $P$). Then ${\cal X}_P$ is a
three-dimensional variety with the isolated multi-quadratic
singularity $o\in{\cal X}_P$, and the blow up of the point $o$
resolves this singularity: the exceptional divisor ${\cal
X}^+_P\cap E_{\cal Y}$ is the non-singular surface $S$.

The restriction of the divisor $D_{\cal X}$ onto ${\cal X}_P$ is
denoted by the symbol $D_P$, and its strict transform on ${\cal
X}^+_P$ by the symbol $D^+_P$.

{\bf Lemma 6.3.} {\it The normal sheaf ${\cal N}_{L/{\cal
X}^+_P}\cong{\cal O}_L(-\alpha)\oplus{\cal O}_L(-\beta)$, where
$\alpha+\beta=k$ and $\alpha\geqslant\beta\geqslant 1$.}

{\bf Proof.} Since $P\cong{\mathbb P}^{k+2}$, by the adjunction
formula $K_S=(k-3)H_{S}$, where $H_S$ is the class of a hyperplane
section of $S\subset P$, whence it follows that $(L^2)_{S}=1-k$.
Furthermore, the surface $S$ is the exceptional divisor of the
blow up of the point $o$ on ${\cal X}_P$ ($S={\cal X}^+_P\cap
E_{\cal Y}$), so that ${\cal O}_{{\cal X}^+_P}(S)|_L={\cal
O}_L(-1)$ and we have the exact sequence
$$
0\to{\cal N}_{L/S}\to{\cal N}_{L/{\cal X}^+_P}\to{\cal N}_{S/{\cal
X}^+_P}|_L\to 0
$$
or
$$
0\to{\cal O}_L(1-k)\to{\cal N}_{L/{\cal X}^+_P}\to{\cal O}(-1)\to
0.
$$
From this the claim of the lemma follows at once. Q.E.D.

Let $\sigma_L\colon{\cal X}_{P,L}\to{\cal X}^+_P$ be the blow up
of the line $L$, and $E_{P,L}\subset{\cal X}_{P,L}$ the
exceptional divisor. The lemma implies that $E_{P,L}$ is a ruled
surface of type ${\mathbb F}_{\alpha-\beta}$ and its Picard group
is ${\mathbb Z}s\oplus{\mathbb Z}f$, where $f$ is the class of a
fibre, $s$ the class of the exceptional section,
$s^2=-(\alpha-\beta)$. Again from the lemma shown above it follows
that
$$
(E^3_{P,L})_{{\cal X}_{P,L}}=(E_{P,L}|^2_{E_{P,L}})= -\mathop{\rm
deg}{\cal N}_{L/{\cal X}^+_P}=k,
$$
so that
$$
-E_{P,L}|_{P,L}=s+\frac12(k+\alpha-\beta)f.
$$
Obviously (since the subspace $P$ is general),
$$
\mathop{\rm mult}\nolimits_pD^+_P=\mu_p.
$$
Let ${\cal X}^{(p)}_P\subset{\cal X}_p$ be the strict transform of
${\cal X}^+_P$ on ${\cal X}_p$. By construction, $q\in{\cal
X}^{(p)}_P$. Setting $D^{(p)}_P=D^{(p)}_{\cal X}|_{{\cal
X}^{(p)}_P}$, we obtain
$$
\mu_q=\mathop{\rm mult}\nolimits_qD^{(p)}_P.
$$
Finally, let $D^{(L)}_P$ be the strict transform of the divisor
$D^+_P$ on ${\cal X}_{P,L}$. Obviously,
$$
D^{(L)}_P=\sigma^*_LD^+_P-\mu E_{P,L},
$$
so that, writing the pull back on ${\cal X}_{P,L}$ of the
restriction $E_{\cal X}|_{{\cal X}^+_P}$ for simplicity as the
restriction $E_{\cal X}|_{{\cal X}_{P,L}}$, we have
$$
(-\nu E_{\cal X}|_{{\cal X}_{P,L}}-\mu E_{P,L})|_{E_{P,L}}\sim
$$
$$
\sim\nu f+\mu(s+\frac12(k+\alpha-\beta)f)=\mu s+(\nu+\frac12\mu
(k+\alpha-\beta))f.
$$
By Proposition 6.2, this effective divisor contains the fibre
$\sigma^{-1}_L(p)$ with multiplicity at least $\mu_p+\mu_q-2\mu$,
whence we get the inequality
$$
\nu+\frac12\mu(k+\alpha-\beta)\geqslant\mu_p+\mu_q-2\mu,
$$
which after easy transformations gives us that
$$
\mu>\frac{2(\mu_p+\mu_q)-2\nu}{k+(\alpha-\beta)+4}.
$$
The denominator of the right hand side is maximal when
$\alpha=k-1$ and $\beta=1$ and so
$$
\mu>\frac{2(\mu_p+\mu_q)-2\nu}{2k+2}=\frac{(\mu_p+\mu_q)-\nu}{k+1}.
$$
Q.E.D. for the theorem.

Proposition 5.3 follows immediately from the theorem that we have
just shown, taking into account the construction of the line $L$
and the inequality (\ref{22.09.22.1}).


\section{Hypertangent divisors}

In this section we prove Theorems 1.2, 1.3 and 1.4.\vspace{0.3cm}

{\bf 7.1. Non-singular points. Tangent divisors.} Let us start the
proof of Theorem 1.2. Obviously, it is sufficient to consider the
case when the subspace $P$ is of maximal admissible codimension
$k+\varepsilon(k)-1$ in ${\mathbb P}^{M+k}$. Theorem 1.1 and
Remark 1.4 imply that the inequality
$$
\mathop{\rm codim}(\mathop{\rm Sing}(F\cap P)\subset(F\cap
P))\geqslant 2k+2
$$
holds. In particular, $F\cap P$ is a factorial complete
intersection of codimension $k$ in ${\mathbb P}\cong{\mathbb
P}^{M-\varepsilon(k)+1}$. Moreover, by the lefschetz theorem,
applied to the section of the variety $F\cap P$ by a general
linear subspace of dimension $3k+1$ in $P$ (this section is a
non-singular complete intersection of codimension $k$ in ${\mathbb
P}^{3k+1}$), we get that the section of $F\cap P$ by an arbitrary
linear subspace of codimension $a\leqslant k$ is irreducible and
reduced, since for the numerical Chow group we have
$$
A^aF\cap P={\mathbb Z}H^a_{F\cap P},
$$
where $H_{F\cap P}$ is the class of a hyperplane section.

Assume that Theorem 1.2 is not true and $\mathop{\rm
mult_o}Y>2n(Y)$. We will argue precisely as in \cite[\S
2]{Pukh01}, see also \cite[Chapter 3, Section 2.1]{Pukh13a}. Let
$T_1,\dots,T_k$ be the tangent hyperplane sections of $F\cap P$ at
the point $o$ (in the notations of Subsection 1.4 they are defined
by the linear forms $f_{i,1}|_P$, $i=1,\dots,k$). By what was said
above, for each $i=1,\dots,k$ the intersection $T_1\cap\dots\cap
T_i$ is of codimension $i$ in $F\cap P$, coincides with the
scheme-theoretic intersection $(T_1\circ\dots\circ T_i)$ and its
multiplicity at the point $o$ equals precisely $2^i$, since the
quadratic forms
$$
f_{1,2}|_{T_o(F\cap P)},\dots,f_{k,2}|_{T_o(F\cap P)}
$$
satisfy the regularity condition. Now we argue as in \cite[\S
2]{Pukh01}. We set $Y_1=Y$ and see that $Y_1\neq T_1$, because
$\mathop{\rm mult}_oT_1=2n(T_1)=2$. We consider the cycle
$(Y_1\circ T_1)$ of the scheme-theoretic intersection and take for
$Y_2$ the component of that cycle that has the maximal value of
the ratio $\mathop{\rm mult}_o/\mathop{\rm deg}$. Assume that the
subvariety $Y_i$ of codimension $i$ in $F\cap P$, satisfying the
inequality
$$
\mathop{\rm mult}\nolimits_oY_i>\frac{2^i}{\mathop{\rm
deg}F}\mathop{\rm deg}Y_i,
$$
is already constructed, and $i\leqslant k-1$. Then
$$
Y_i\neq T_1\cap\dots\cap T_i,
$$
however by construction $Y_i$ is contained in the divisors
$T_1,\dots,T_{i-1}$, so that $Y_i\not\subset T_i$ and the cycle of
scheme-theoretic intersection $(Y_i\circ T_i)$ of codimension
$i+1$ is well defined. For $Y_{i+1}$ we take the component of this
cycle with the maximal value of the ratio $\mathop{\rm
mult}_o/\mathop{\rm deg}$. Completing this process, we obtain an
irreducible subvariety $Y_{k+1}\subset(F\cap P)$ of codimension
$k+1$, satisfying the inequality
$$
\frac{\mathop{\rm mult}_o}{\mathop{\rm deg}}
Y_{k+1}>\frac{2^{k+1}}{\mathop{\rm deg}F}.
$$
\vspace{0.3cm}


{\bf 7.2. Non-singular points. Hypertangent divisors.} We continue
the proof of Theorem 1.2. In the notations of Subsection 1.4 for
each $j=2,\dots, d_k-1$ construct the hypertangent linear systems
$$
\Lambda_j=\left|\sum^k_{i=1}\sum^{\min\{j,d_i-1\}}_{\alpha=1}
f_{i,[1,\alpha]}s_{i,j-\alpha}\right|_{F\cap P},
$$
where $f_{i,[1,\alpha]}=f_{i,1}+\dots +f_{i,\alpha}$ is the left
segment of the polynomial $f_i$ of length $\alpha$, the
polynomials $s_{i,j-\alpha}$ are homogeneous polynomials of degree
$j-\alpha$, running through the spaces ${\cal P}_{j-\alpha,M+k}$
independently of each other and the restriction onto $F\cap P$
means the restriction onto the affine part of that variety in
${\mathbb A}^{M+k}_{z_*}$ followed by the closure.

Let $h_a$, where $a\geqslant k+1$, be the $a$-th polynomial in the
sequence ${\cal S}$. Then $h_a=f_{i,j}|_{{\mathbb P}(T_oF)}$ for
some $i$ and $j\geqslant 3$. Set ${\cal H}_a=\Lambda_{j-1}$. In
this way we obtain a sequence of linear systems ${\cal H}_{k+1}$,
${\cal H}_{k+2}$,\dots, ${\cal H}_{M}$, where the system
$\Lambda_j$ occurs, in the notations of \cite[\S 2]{Pukh01},
$$
w^+_j=\sharp\{i, 1\leqslant i\leqslant k\,|\, j\leqslant d_i-1\}
$$
times. By the symbol ${\cal H}[-m]$ we denote the space
$$
\prod^{M-m}_{a=k+1}{\cal H}_a
$$
of all tuples $(D_{k+1},\dots,D_{M-m})$ of divisors, where
$D_a\in{\cal H}_a$. For $a\in\{k+1,\dots,M\}$ set
$$
\beta_a=\frac{j+1}{j},
$$
if ${\cal H}_a=\Lambda_j$. The number $\beta_a$ is called the {\it
slope} of the divisor $D_a$. It is easy to see that
\begin{equation}\label{24.10.22.1}
\prod^M_{a=k+1}\beta_a=\frac{d_1\dots d_k}{2^k}=\frac{\mathop{\rm
deg} F}{2^k}.
\end{equation}
Set $m_*=k+\varepsilon(k)+3$. Let
$$
(D_{k+1},\dots,D_{M-m_*})\in{\cal H}[-m_*]
$$
be a general tuple. The technique of hypertangent divisors,
applied in precisely the same way as in \cite[\S 2]{Pukh01} or
\cite[Chapter 3, Subsection 2.2]{Pukh13a}, see also
\cite[Proposition 2.1]{Pukh2022a}, gives the following claim.

{\bf Proposition 7.1.} {\it There is a sequence of irreducible
subvarieties
$$
Y_{k+1},Y_{k+2},\dots,Y_{M-m_*},
$$
where $Y_{k+1}$ has been constructed above, such that $\mathop{\rm
codim}(Y_i\subset(F\cap P))=i$, the subvariety $Y_i$ is not
contained in the support of the divisor $D_{i+1}$ for $i\leqslant
M-m_*-1$, the subvariety $Y_{i+1}$ is an irreducible component of
the effective cycle $(Y_i\circ D_{i+1})$ and the following
inequality holds:}
$$
\frac{\mathop{\rm mult}_o}{\mathop{\rm deg}}Y_{i+1}\geqslant
\beta_{i+1}\frac{\mathop{\rm mult}_o}{\mathop{\rm deg}}Y_i.
$$

There is no need to give a {\bf proof} of that claim because it is
identical to the arguments mentioned above. Note only that the key
point in the construction of the sequence of subvarieties $Y_i$ is
the fact that $Y_i$ is not contained in the support of a general
divisor $D_{i+1}\in{\cal H}_{i+1}$, and this fact follows from the
regularity condition (R1). Since $\mathop{\rm dim}(F\cap
P)=M+1-k-\varepsilon(k)$, the subvariety $Y^*=Y_{M-m_*}$ is of
dimension 4 and satisfies the inequality
$$
\frac{\mathop{\rm mult}_o}{\mathop{\rm deg}}Y^*>
\frac{\displaystyle\frac{2^{k+1}}{\mathop{\rm
deg}F}\cdot\frac{\mathop{\rm deg}
F}{2^k}}{\displaystyle\frac32\cdot\prod\limits^M_{a=M-m_*+1}\beta_a}=
\frac43\frac{1}{\prod\limits^M_{a=M-m_*+1}\beta_a}.
$$
(The number $\frac32$ appears in the denominator because the
hypertangent divisor $D_{k+1}$ is skipped in the procedure of
intersection, in the same way and for the same reason as in
\cite[\S 2]{Pukh01}, and its slope is $\frac32$.) Now the
inequality
\begin{equation}\label{14.10.22.1}
\frac43\geqslant\prod^M_{a=M-m_*+1}\beta_a,
\end{equation}
shown below in Proposition 7.2, completes the proof of Theorem
1.2.

{\bf Proposition 7.2.} {\it Assume that for $k=3,4,5$ the
dimension $M$ is, respectively, at least $96$, $160$, $215$, and
for $k\geqslant 6$ the inequality $M\geqslant 8k^2+2k$ holds. then
the inequality (\ref{14.10.22.1}) is true.}

{\bf Proof.} Using the obvious fact that the function
$\frac{t+1}{t}$ is decreasing, it is easy to see that the right
hand side of the inequality (\ref{14.10.22.1}) with $k$ and $M$
fixed attains the maximum when the degrees $d_1,\dots,d_k$ are
equal or ``almost equal'' in the following sense: let $M\equiv
e\mathop{\rm mod}k$ with $e\in\{0,1,\dots,k-1\}$, then the
``almost equality'' means that
$$
d_1=\dots=d_{k-e}=\frac{M-e}{k}+1,\quad
d_{k-e+1}=\dots=d_k=\frac{M-k}{k}+2.
$$
For $k\in\{3,\dots,9\}$ the claim of the proposition can be
checked for each case of almost equal degrees, that is, for each
possible value of $e$, manually, computing $\varepsilon(k)$
explicitly. For $k\geqslant 10$ it is easy to see that
$\varepsilon(k)\leqslant k-3$, so that $m_*\leqslant 2k$.
Therefore (again considering the case of almost equal degrees) the
right hand side of (\ref{14.10.22.1}) does not exceed the number
$$
\left(\frac{\frac{M}{k}}{\frac{M}{k}-2}\right)^k=
\left(\frac{M}{M-2k}\right)^k,
$$
from which we get that (\ref{14.10.22.1}) is true if
$$
M\geqslant
2k\frac{(1+\frac13)^\frac{1}{k}}{(1+\frac13)^\frac{1}{k}-1}.
$$
If in the numerator and denominator we replace
$(1+\frac{1}{3})^{\frac{1}{k}}$ by the smaller number
$1+\frac{1}{4k}$, the right hand side of the last inequality gets
higher. This proves the proposition. Q.E.D.

Note that for $M\geqslant \rho(k)$ (see the inequality
(\ref{14.11.22.1}) in Subsection 0.1) the assumptions of the
previous proposition are satisfied. This completes the proof of
Theorem 1.2.\vspace{0.3cm}


{\bf 7.3. Quadratic points.} Let us show Theorem 1.3. Note first
of all that if $Y$ is a section of the variety $W$ by a hyperplane
that is tangent to $W$ at the point $o$ (that is, the equation of
the hyperplane is a linear combination of the forms
$f_{1,1},\dots,f_{k,1}$, restricted onto the hyperplane ${\mathbb
P}(W)\cong{\mathbb P}^{M+k-1}$), then
$$
\mathop{\rm mult}\nolimits_oY=4n(Y)=4,
$$
so that the claim of the theorem is optimal. Thus we assume the
converse: the inequality
$$
\mathop{\rm mult}\nolimits_oY>2n(Y)
$$
holds. We argue as in the non-singular case (Subsection 7.1): let
$T_1,\dots,T_{k-1}$ be the tangent hyperplane sections, given by
$(k-1)$ independent forms taken from the set
$\{f_{1,1},\dots,f_{k,1}\}$. Since $\mathop{\rm codim}(\mathop{\rm
Sing}F\subset F)\geqslant2k+2$, all scheme-theoretic intersections
$(T_1\circ\dots\circ T_i)$, $1\leqslant i\leqslant k-1$, are
irreducible, reduced and coincide with the set-theoretic
intersection $T_1\cap\dots\cap T_i$, and moreover, by the
condition (R2) the equality
$$
\mathop{\rm mult}\nolimits_oT_1\cap\dots\cap T_i=2^{i+1}
$$
holds. In particular, $\mathop{\rm mult}_oT_1=4n(T_1)=4$, so that
$Y\neq T_1$. Arguing as in Subsection 7.1, we construct a sequence
of irreducible subvarieties $Y_1=Y,Y_2,\dots,Y_k$, where
$\mathop{\rm codim}(Y_i\subset W)=i$ and the inequality
$$
\frac{\mathop{\rm mult}_o}{\mathop{\rm
deg}}Y_k>\frac{2^{k+1}}{\mathop{\rm deg}F}
$$
holds. Now the proof of Theorem 1.3 repeats the arguments of
Subsection 7.2, where $m_*$ is replaced by 4. Since $4<m_*$, the
inequality (\ref{14.10.22.1}) guarantees the inequality which is
obtained from (\ref{14.10.22.1}) when $m_*$ is replaced by 4. This
completes the proof of Theorem 1.3.\vspace{0.3cm}


{\bf 7.4. Multi-quadratic points. Tangent divisors.} We start the
proof of Theorem 1.4, the structure of which is similar to the
structure of the proof of Theorem 1.2. At first we argue as in
Subsection 7.1: it is sufficient to consider a linear subspace $P$
in $T_oF$ of maximal admissible codimension $\varepsilon(k)$.
Assume that the prime divisor $Y$ on $F\cap P$ satisfies the
inequality
$$
\mathop{\rm mult}\nolimits_oY>\frac32\cdot2^kn(Y),
$$
or the equivalent inequality
$$
\frac{\mathop{\rm mult}_o}{\mathop{\rm
deg}}Y>\frac32\cdot\frac{2^k}{\mathop{\rm deg}F},
$$
and consider the second hypertangent linear system (which in this
case plays the role of the tangent linear system)
$$
\Lambda_2=\left|\sum_{d_i\geqslant 3}s_{i,0}f_{i,2}\right|_{F\cap
P},
$$
where $s_{i,0}\in{\mathbb C}$ are constants, independent of each
other. Instead of the Lefschetz theorem we use the condition
(R3.1): the system of equations $f_{i,2}|_{F\cap P}=0$, where
$d_i\geqslant 3$, defines an irreducible reduced subvariety of
codimension $k+k_{\geqslant 3}$ in $P$, and by (R3.2) the
multiplicity of that subvariety at the point $o$ is precisely
$2^k\cdot(\frac32)^{k_{\geqslant 3}}$. More precisely, for a
general tuple $(D_{2,1},\dots,D_{2,k_{\geqslant 3}})$ of divisors
in the system $\Lambda_2$ the following claim is true: for each
$i=1,\dots,k_{\geqslant 3}$ the cycle $(D_{2,1}\circ\dots\circ
D_{2,i}$) of the scheme-theoretic intersection of the divisors
$D_{2,1},\dots,D_{2,i}$ is an irreducible reduced subvariety of
codimension $i$ in $F\cap P$, the multiplicity of which at the
point $o$ is $2^k\cdot(\frac32)^i$. Arguing as in Subsection 7.1,
we construct a sequence $Y_1=Y,Y_2,\dots,Y_{k_{\geqslant 3}}$ of
irreducible subvarieties of codimension $\mathop{\rm
codim}(Y_i\subset(F\cap P))=i$, where $Y_{i+1}$ is an irreducible
component of the cycle $(Y_i\circ D_{2,i})$ with the maximal value
of the ratio $\mathop{\rm mult}_o/\mathop{\rm deg}$. Therefore,
$$
\frac{\mathop{\rm mult}_o}{\mathop{\rm deg}}Y_{k_{\geqslant 3}}>
\left(\frac32\right)^{k_{\geqslant 3}}\cdot\frac{2^k}{\mathop{\rm
deg}F}.
$$
It follows from here that $Y_{k_{\geqslant 3}}\neq
D_{2,1}\cap\dots\cap D_{2,k_{\geqslant 3}}$, but since by
construction
$$
Y_{k_{\geqslant 3}}\subset D_{2,1}\cap\dots\cap D_{2,k_{\geqslant
3}-1},
$$
we conclude that $Y_{k_{\geqslant 3}}\not\subset D_{2,k_{\geqslant
3}}$, so that the effective cycle $(Y_{k_{\geqslant 3}}\circ
D_{2,k_{\geqslant 3}})$ of the scheme-theoretic intersection of
these varieties is well defined and one of its components
$Y_{k_{\geqslant 3}+1}$ satisfies the inequality
$$
\frac{\mathop{\rm mult}_o}{\mathop{\rm deg}}Y_{k_{\geqslant 3}+1}>
\left(\frac32\right)^{k_{\geqslant
3}+1}\cdot\frac{2^k}{\mathop{\rm deg}F}.
$$\vspace{0.3cm}


{\bf 7.5. Multi-quadratic points. Hypertangent divisors.} Now we
argue almost word for word as in Subsection 7.2: construct the
hypertangent systems
$$
\Lambda_j=\left|\sum^j_{\alpha=2}\sum_{d_i\geqslant
\alpha+1}f_{i,[2,\alpha]}s_{i,j-\alpha}\right|_{F\cap P},
$$
where $j=3,\dots,d_k-1$ and all symbols have the same meaning as
in Subsection 7.2. If $h_a$, where $a\geqslant k+k_{\geqslant
3}+1$, is the $a$-th polynomial in the sequence ${\cal
S},h_a=f_{i,j}|_{{\mathbb P}(T_oF)}$, for some $i$ and $j\geqslant
4$, then we set ${\cal H}_a=\Lambda_{j-1}$ and obtain the sequence
of linear systems
$$
{\cal H}_{k+k_{\geqslant 3}+1},\quad {\cal H}_{k+k_{\geqslant
3}+2},\quad\dots,\quad {\cal H}_M.
$$
By the symbol ${\cal H}[-m]$ we denote the space
$$
\prod^{M-m}_{a=k+k_{\geqslant 3}+1}{\cal H}_a.
$$
Instead of the equality (\ref{24.10.22.1}) we get the equality
$$
\prod^M_{a=k+k_{\geqslant 3}+1}\beta_a=\frac{\mathop{\rm
deg}F}{2^k\left(\frac32\right)^{k_{\geqslant 3}}}.
$$
Let $(D_{k+k_{\geqslant 3}+1},\dots,D_{M-m^*})\in{\cal H}[-m^*]$
be a general tuple. Now the technique of hypertangent divisors,
applied in the word for word the same way as in Subsection 7.2,
gives the following claim.

{\bf Proposition 7.3.} {\it There is a sequence of irreducible
subvarieties
$$
Y_{k_{\geqslant 3}+1},Y_{k_{\geqslant 3}+2},\dots,Y_{M-k-m^*},
$$
where $Y_{k_{\geqslant 3}+1}$ is constructed above, such that
$\mathop{\rm codim}(Y_i\subset(F\cap P))=i$, the subvariety $Y_i$
is not contained in the support of the divisor $D_{k+i+1}$ for
$i\leqslant M-m^*-1$, the subvariety $Y_{i+1}$ is an irreducible
component of the effective cycle $(Y_i\circ D_{k+i+1})$ and the
following inequality holds:}
$$
\frac{\mathop{\rm mult}_o}{\mathop{\rm
deg}}Y_{i+1}\geqslant\beta_{k+i+1}\frac{\mathop{\rm
mult}_o}{\mathop{\rm deg}}Y_i.
$$

Now since $\mathop{\rm dim}F\cap P=M-(k-l)-\varepsilon(k)$, by the
definition of the number $m^*$ the last subvariety
$Y^*=Y_{M-k-m^*}$ in that sequence is of dimension $\geqslant 4$
and satisfies the inequality
$$
\frac{\mathop{\rm mult}_o}{\mathop{\rm deg}}Y^*> \frac{
\displaystyle \left(\frac32\right)^{k_{\geqslant
3}+1}\cdot\frac{2^k}{\mathop{\rm deg}F}\cdot \frac{\mathop{\rm
deg}F}{2^k \left(\frac32\right)^{k_{\geqslant 3}}}}{\displaystyle
\frac43\prod\limits^M_{a=M-m_*+1}\beta_{k+a}}=
\frac98\frac{1}{\displaystyle\prod\limits^M_{a=M-m_*+1}\beta_{k+a}}.
$$
(The number $\frac43$ appears in the denominator of the right hand
side, because the hypertangent divisor $D_{k+k_{\geqslant 3}+1}$
is skipped in the process of constructing the sequence $Y_*$, see
the similar remark above, before the inequality
(\ref{14.10.22.1}).) If $m^*=0$, then the product in the
denominator is assumed to be equal to 1. Now the inequality
\begin{equation}\label{28.10.22.1}
\frac98\geqslant\prod^M_{a=M-m^*+1}\beta_{k+a},
\end{equation}
shown below in Proposition 7.4, completes the proof of Theorem
1.4.

{\bf Proposition 7.4.} {\it Assume that for $k\in\{3,\dots,7\}$
the number $M$ is at least the number shown in the corresponding
column of the table
\begin{center}
\begin{tabular}{|c|c|c|c|c|c|}
\hline $k$ & $3$ & $4$ & $5$ & $6$ & $7$\\
\hline $M\geqslant$ & $128$ & $204$ & $255$ & $357$ & $477$\\
\hline
\end{tabular},
\end{center}
and for $k\geqslant 8$ the inequality $M\geqslant 9k^2+k$ holds.
Then the inequality (\ref{28.10.22.1}) holds.}

{\bf Proof.} As in the non-singular case (the proof of Proposition
7.2), we see that the right hand side of the inequality
(\ref{28.10.22.1}) for $k$ and $M$ fixed, attains the maximum when
the degrees $d_i$ are equal or ``almost equal''. For
$k\in\{3,\dots,7\}$ the claim of the proposition is checked
manually. For $k\geqslant 8$ we have $\varepsilon(k)\leqslant
k-2$, so that $m^*\leqslant k$. Therefore (considering the case of
equal or almost equal degrees) the right hand side of the
inequality (\ref{28.10.22.1}) does not exceed the number
$$
\left(\frac{M}{M-k}\right)^k,
$$
which, in its turn, does not exceed $\frac98$ for $M\geqslant
9k^2+k$, which is easy to check by elementary computations,
similar to the proof Proposition 7.2. Q.E.D.


\section{The codimension of the complement}

In this section we show the estimate for the codimension of the
complement ${\cal P}\setminus{\cal F}$, given in Theorem
0.1.\vspace{0.3cm}

{\bf 8.1. Preliminary constructions.} Set
$$
\gamma=M-k+5+{M-\rho(k)+2\choose 2},
$$
see Subsection 0.1. We consider $\gamma$ as a function of $M$ with
$k\geqslant 3$ fixed, where $M\geqslant\rho(k)$. Let $o\in{\mathbb
P}^{M+k}$ be an arbitrary point. The symbol ${\cal P}(o)$ stands
for the linear subspace of the space ${\cal P}$, consisting of all
tuples $\underline{f}$, vanishing at the point $o$:
$\underline{f}(o)=(0,\dots,0)$. Obviously, $\mathop{\rm
codim}({\cal P}(o)\subset{\cal P})=k$. Fixing the point $o$, we
use the notations of Subsections 1.2-1.4, considering the
polynomials $f_i$ as non-homogeneous polynomials in the affine
coordinates $z_*$. By the symbols
$$
{\cal B}_{MQ1},\,{\cal B}_{MQ2},\,{\cal B}_{R1},\,{\cal
B}_{R2},\,{\cal B}_{R3.1},\,{\cal B}_{R3.2}
$$
we denote the subsets of the subspace ${\cal P}(o)$, consisting of
the tuples $\underline{f}$ that do not satisfy the conditions
$$
(MQ1),\,(MQ2),\,(R1),\,(R2),\,(R3.1),\,(R3.2)
$$
at the point $o$, respectively. Since the point $o$ varies in
${\mathbb P}^{M+k}$, it is sufficient to show that the codimension
of each of the six sets ${\cal B}_*$ in ${\cal P}(o)$ is at least
$\gamma+M$.

Furthermore, for an arbitrary tuple
$$
\underline{\xi}=(\xi_1,\dots,\xi_k)
$$
of linear forms in $z_*$ the symbol ${\cal P}(o,\underline{\xi})$
denotes the affine subspace, consisting of the tuples
$\underline{f}$, such that
$$
f_{1,1}=\xi_1,\quad\dots,\quad f_{k,1}=\xi_k.
$$
By the symbol $\mathop{\rm dim}\underline{\xi}$ denote the
dimension
$$
\mathop{\rm dim}\langle\xi_1,\dots,\xi_k\rangle,
$$
so that ${\cal P}(o)$ is fibred into disjoint subsets
$$
{\cal P}^{(i)}(o)=\bigcup_{\mathop{\rm dim}\underline{\xi}=i}{\cal
P}(o,\underline{\xi}),
$$
where $i=0,1,\dots,k$. Obviously, the equality
$$
\mathop{\rm codim}({\cal P}^{(i)}(o)\subset{\cal
P}(o))=(k-i)(M+k-i)
$$
holds. In particular, ${\cal P}^{(k)}(o)$ consists of the tuples
$\underline{f}$, such that the scheme of their common zeros is a
non-singular subvariety of codimension $k$ in a neighborhood of
the point $o$. Set
$$
{\cal B}_{R1}(\underline{\xi})={\cal B}_{R1}\cap{\cal
P}(o,\underline{\xi}).
$$
For the case of a non-singular point it is sufficient to prove the
inequality
$$
\mathop{\rm codim}({\cal B}_{R1}(\underline{\xi})\subset{\cal
P}(o,\underline{\xi}))\geqslant \gamma+M,
$$
where $\mathop{\rm dim}\underline{\xi}=k$.

Furthermore, let ${\cal B}_{MQ1}(\underline{\xi})={\cal
B}_{MQ1}\cap{\cal P} (o,{\underline{\xi}})$, where $\mathop{\rm
dim}\underline{\xi}=i\leqslant k-1$, be the set of the tuples
$\underline{f}$, such that the condition (MQ1) for $l=k-i$ is not
satisfied, and ${\cal B}_{MQ2}(\underline{\xi})={\cal
B}_{MQ2}\cap{\cal P }(o,\underline{\xi})$, where $\mathop{\rm
dim}\underline{\xi}=i\leqslant k-2$, be the set of the tuples
$\underline{f}$, such that the condition (MQ2) for $l=k-i$ is not
satisfied.

In a similar way, we define the sets ${\cal
B}_{R2}(\underline{\xi})$ for $\mathop{\rm
dim}\underline{\xi}=k-1$ and ${\cal B}_{R3.1}(\underline{\xi})$,
${\cal B}_{R3.2}(\underline{\xi})$ for $\mathop{\rm
dim}\underline{\xi}\leqslant k-2$.

Clearly, it is sufficient to prove that for $\mathop{\rm
dim}\underline{\xi}=i$ the codimension of the set ${\cal
B}_*(\underline{\xi})$ in ${\cal P}(o,\underline{\xi})$ is at
least
$$
\gamma+M-(k-i)(M+k-i).
$$

In the conditions (R1),(R2) and (R3.2) we have also an arbitrary
subspace $\Pi\subset{\mathbb P}(T_oF)$ of the corresponding
codimension, and in the condition (R3.1) an arbitrary subspace $P$
in the embedded tangent space $T_oF\subset{\mathbb P}^{M+k}$ of
codimension $\varepsilon(k)$, containing the point $o$. For an
arbitrary subspace $\Pi\subset{\mathbb P}(T_oF)$ of the
corresponding codimension let
$$
{\cal B}_{R1}(\underline{\xi},\Pi),\quad {\cal
B}_{R2}(\underline{\xi},\Pi),\quad {\cal
B}_{R3.2}(\underline{\xi},\Pi)
$$
be the set of tuples $\underline{f}\in{\cal
P}(o,\underline{\xi})$, such that the respective condition
(R1),(R2) and (R3.2) is violated precisely for that subspace
$\Pi$. In a similar way we define the subset ${\cal
B}_{R3.1}(\underline{\xi},P)$. These definitions are meaningful
because the tangent space $T_oF$ is given by the fixed linear
forms $\xi_i$ and for that reason is fixed.

Since the subspace $\Pi$ varies in a $(\mathop{\rm
dim}\Pi+1)\mathop{\rm codim}(\Pi\subset{\mathbb
P}(T_oF))$-dimensional Grassmanian, the estimate for the
codimension of the set ${\cal B}_*(\underline{\xi},\Pi)$ in ${\cal
P}(o,\underline{\xi})$ should be stronger than the estimate for
the codimension of the set ${\cal B}_*(\underline{\xi})$ by that
number. Similarly, $P$ varies in a
$$
\varepsilon(k)(\mathop{\rm dim}T_oF-\varepsilon(k))
$$
-dimensional family, so that the estimate for the codimension of
the set ${\cal B}_{R3.1}(\underline{\xi},P)$ should be stronger
than the estimate for ${\cal B}_{R3.1}(\underline{\xi})$ by that
number.

Now everything is ready to consider each of the six subsets ${\cal
B}_*$.\vspace{0.3cm}


{\bf 8.2. The conditions (MQ1) and (MQ2).} For a non-singular
point $o\in F$ these conditions contain no restrictions, so we
assume that $\dim\underline{\xi}\leqslant k-1$. It is well known
that the closed subset of quadratic forms of rank $\leqslant
r\leqslant N-1$ in the space ${\cal P}_{2,N}$ has the codimension
$$
{N-r+1\choose 2}.
$$
From here it is easy to see that the closed subset of tuples
$(q_1,\dots,q_e)=q_{[1,e]}$ of quadratic forms in $N$ variables,
defined by the condition
$$
\mathop{\rm rk}q_{[1,e]}\leqslant r,
$$
is of codimension
$$
\geqslant{N-r+1\choose 2}-(e-1)
$$
in the space ${\cal P}^{\times e}_{2,N}$. As we noted in
Subsection 1.2 (after stating the condition (MQ2)), for
$l\geqslant 2$ the condition (MQ2) is stronger than (MQ1),
therefore it is sufficient to estimate the codimension of the set
${\cal B}_{MQ2}$ (in the case of quadratic points, when $l=1$, it
is easy to check that the codimension of the set ${\cal B}_{MQ1}$
is higher than required). So we assume that $\mathop{\rm
dim}\underline{\xi}=k-l\leqslant k-2$. The condition (MQ2)
requires the rank of the tuple of quadratic forms $q_{[1,k]}$,
where $q_i=f_{i,2}|_{T_oF}$, to be at least $\rho(k)+2$, see
(\ref{14.11.22.1}) in Subsection 0.1. Taking into account the
variation of the tuple $\underline{\xi}$, from what was said above
it is easy to obtain that the codimension of the set ${\cal
B}_{MQ2}\cap{\cal P}^{k-l}(o)\geqslant$
$$
-k+1+l(M+l)+{M+l-\rho(k)\choose 2}.
$$
The minimum of this expression is attained for $l=2$ and it is
easy to check that this minimum is precisely $\gamma+M$.
Therefore, the codimension of the set ${\cal B}_{MQ2}$ is at least
$\gamma$, the codimension of the set ${\cal B}_{MQ1}$ for
$l\geqslant 2$ is higher. For $l=1$ the last codimension is also
higher. This completes our consideration of the conditions (MQ1)
and (MQ2).\vspace{0.3cm}


{\bf 8.3. Regularity at the non-singular and quadratic points.}
Let us estimate the codimension of the set ${\cal
B}_{R1}(\underline{\xi},\Pi)$ in the space ${\cal
P}(o,\underline{\xi})$. Here $\mathop{\rm dim}\underline{\xi}=k$
and $\Pi\subset{\mathbb P}(T_oF)$ is a subspace of codimension
$k+\varepsilon(k)-1=m_*-4$. Let
$$
{\cal G}(\underline{d})=\prod^{M-m_*}_{i=1}{\cal P}_{\mathop{\rm
deg}h_i,\mathop{\rm dim}\Pi+1}
$$
be the space, parameterizing all sequences ${\cal
S}[-m_*]|_{\Pi}$. Since the polynomials $h_i$ are distinct
homogeneous components of the polynomials of the tuple
$\underline{f}$, restricted onto the subspace $\Pi$, the
codimension of the subset ${\cal B}_{R1}(\underline{\xi},\Pi)$ in
${\cal P}(o,\underline{\xi})$ is equal to the codimension of the
subset ${\cal B}\subset{\cal G}(\underline{d})$, which consists of
the sequences that do not satisfy the condition (R1).

Using the approach that was applied in
\cite{Pukh98b,Pukh13a,EvansPukh2} and many other papers, let us
present ${\cal B}$ as a disjoint union
$$
{\cal B}=\bigsqcup^{M-m_*}_{i=1}{\cal B}_i,
$$
where ${\cal B}_i$ consists of sequences
$$
(h_1|_{\Pi},\dots,h_{M-m_*}|_{\Pi}),
$$
such that the first $i-1$ polynomials form a regular sequence but
$h_i$ vanishes on one of the components of the set of their common
zeros. The ``projection method'' estimates the codimension of
${\cal B}_i$ in ${\cal G}(\underline{d})$ from below by the
integer
\begin{equation}\label{17.10.22.1}
{\mathop{\rm dim}\Pi-i+1+\mathop{\rm deg}h_i\choose \mathop{\rm
deg}h_i} ={\mathop{\rm dim}\Pi-i+1+\mathop{\rm deg}h_i\choose
\mathop{\rm dim}\Pi-i+1}
\end{equation}
(we will use both presentations). It follows easily from here (see
\cite[\S3]{EvansPukh2}), that the worst estimate corresponds to
the case of equal or ``almost equal'' degrees $d_i$, described
above. We will consider this case.

Thus we need to estimate from below the minimum of $M-m_*$
integers (\ref{17.10.22.1}). Here are the first $(k+1)$ of them:
$$
{\mathop{\rm dim}\Pi+2\choose 2},\,{\mathop{\rm dim}\Pi+1\choose
2},\,\dots,\,{\mathop{\rm dim}\Pi+3-k\choose 2},\,{\mathop{\rm
dim}\Pi+3-k\choose 3}.
$$
We call the left hand side of the equality (\ref{17.10.22.1}) the
presentation of type (I), the right hand side is the presentation
of type (II). Let us write down each of the numbers
(\ref{17.10.22.1}) in the form
$$
{A(i)\choose B(i)},
$$
where $A(i)\geqslant 2B(i)$, using the presentation of type (I) or
of type (II).

At first (for the starting segment of the sequence) we use the
presentation of type (I). It is easy to see that when we change
$i$ by $i+1$, we have one of the two options:
\begin{itemize}
\item either $\mathop{\rm deg}h_{i+1}=\mathop{\rm deg}h_i$, and
then $A(i+1)=A(i)-1$ and $B(i+1)=B(i)$, so that
$$
{A(i+1)\choose B(i+1)}<{A(i)\choose B(i)},
$$
and moreover, $C(i)=A(i)-2B(i)$ decreases: $C(i+1)=C(i)-1$,

\item or $\mathop{\rm deg}h_{i+1}=\mathop{\rm deg}h_i+1$, and then
$A(i+1)=A(i)$ and $B(i+1)=B(i)+1$, so that $C(i+1)=C(i)-2$ and if
$C(i+1)\geqslant 0$, then
$$
{A(i+1)\choose B(i+1)}>{A(i)\choose B(i)}.
$$
\end{itemize}

This is how it goes on until the ``equilibrium'': $C(i_*)\geqslant
0$, but $C(i_*+1)<0$, and after that we use the presentation of
type (II).

Now when we change $i$ by $(i+1)$, we have one of the two options:

\begin{itemize}

\item either $\mathop{\rm deg}h_{i+1}=\mathop{\rm deg}h_i$, and
then $A(i+1)=A(i)-1$ and $B(i+1)=B(i)-1$, so that $C(i+1)=C(i)+1$
and
$$
{A(i+1)\choose B(i+1)}<{A(i)\choose B(i)},
$$

\item or $\mathop{\rm deg}h_{i+1}=\mathop{\rm deg}h_i+1$, and then
$A(i+1)=A(i)$ and $B(i+1)=B(i)-1$, so that $C(i+1)=C(i)+2$ and
$$
{A(i+1)\choose B(i+1)}<{A(i)\choose B(i)}.
$$
\end{itemize}

Therefore, after the ``equilibrium'' our sequence is strictly
decreasing. Moreover, if
$$
{A(i_1)\choose B(i_1)}\quad\mbox{and}\quad {A(i_2)\choose B(i_2)}
$$
are two numbers in our sequence, where $i_1\leqslant i_*$ and
$i_2>i_*$ and $B(i_1)\geqslant B(i_2)$, then, obviously,
$$
{A(i_1)\choose B(i_1)}>{A(i_2)\choose B(i_2)}.
$$
Recall that the degrees $d_i$ are equal or ``almost equal''.

{\bf Lemma 8.1.} {\it For $M\geqslant 3k^2$ the following
inequality holds:} $i_*<M-m_*$.

{\bf Proof.} Elementary computations, using the equality
$C(i+k)=C(i)-(k+1)$ if $C(i+k)\geqslant 0$. Q.E.D. for the lemma.

Therefore, the ``equilibrium'' is reached earlier than the
sequence $h_i,\dots,h_{M-m_*}$ comes to an end, so that there is a
non-empty segment after the ``equilibrium''. By construction,
$B(M-m_*)=4$. By what was said above, the minimum of the numbers
${A(i)\choose B(i)}$ for $i=1,\dots,M-m_*$ is the minimum of the
following three numbers:
$$
{\mathop{\rm dim}\Pi+3-k\choose 2},\quad {\mathop{\rm
dim}\Pi+4-2k\choose 3},\quad {\mathop{\rm deg}h_{M-m_*}+4\choose
4}.
$$

{\bf Lemma 8.2.} {\it For $\mathop{\rm dim}\Pi\geqslant 3k+1$ the
following inequality holds:}
$$
{\mathop{\rm dim}\Pi+4-2k\choose 3}>{\mathop{\rm
dim}\Pi+3-k\choose 2}.
$$

{\bf Proof.} Elementary computations. Q.E.D.

{\bf Lemma 8.3.} {\it For $M\geqslant 2\sqrt{3}k^2$ the following
inequality holds:}
$$
{\mathop{\rm deg}h_{M-m_*}+4\choose 4}>{\mathop{\rm
dim}\Pi+3-k\choose 2}.
$$

{\bf Proof.} It is easy to check the inequalities
$$
\frac{(M-2k)^2}{2}\geqslant{\mathop{\rm dim}\Pi+3-k\choose 2}
$$
and
$$
{\mathop{\rm deg}h_{M-m_*}+4\choose 4}\geqslant\frac{1}{24}
\left(\frac{M}{k}+1\right)\left(\frac{M}{k}\right)\left(\frac{M}{k}-1\right)
\left(\frac{M}{k}-2\right),
$$
so that it is sufficient to show that for $M\geqslant
2\sqrt{3}k^2$ the inequality
$$
\left(\frac{M^2}{k^2}-1\right)\left(\frac{M^2}{k^2}-2\frac{M}{k}\right)
>12(M-2k)^2
$$
holds or, equivalently, $M(M^2-k^2)>12k^4(M-2k)$. It is easy to
check the last inequality, considering the cubic polynomial
$$
t^3-(12k^4+k^2)t+24k^5
$$
in the real variable $t$. Q.E.D. for the lemma.

The work that was carried out above gives the inequality
$$
\mathop{\rm codim}({\cal B}\subset{\cal
G}(\underline{d}))\geqslant{M+3-2k-\varepsilon(k)\choose 2}.
$$
From here by elementary computations (taking into account the
variation of the subspace $\Pi$, see Subsection 8.1) it is easy to
obtain the required the inequality $\mathop{\rm codim}({\cal
B}_{R1}\subset{\cal P}(o))\geqslant\gamma+M$. This completes the
proof in the case of smooth points.

It is easy to see that the methods used above give a stronger
estimate for the codimension of the set ${\cal B}_{R2}$, because
the dimension of the subspace $\Pi$ is higher. The computations
are completely similar to the computations given above for the
case of a non-singular point, for that reason we do not consider
the case of a quadratic point and move on to estimating the
codimension of the sets ${\cal B}_{R3.1}$ and ${\cal
B}_{R3.2}$.\vspace{0.3cm}


{\bf 8.4. Regularity at the multi-quadratic points.} Let us
estimate the codimension of the set ${\cal
B}_{R3.2}(\underline{\xi},\Pi)$, where $\Pi\subset{\mathbb
P}(T_oF)$ is an arbitrary subspace of codimension
$\varepsilon(k)$. Our arguments are completely similar to the
arguments of Subsection 8.3 for a non-singular point and give a
stronger estimate for the codimension. We just point out the
necessary changes in the constructions of Subsection 8.3. Set
$$
{\cal G}(\underline{d})=\prod^{M-m^*}_{i=1}{\cal P}_{\mathop{\rm
deg} h_i,\mathop{\rm dim}\Pi+1}.
$$
Denote by the symbol ${\cal B}$ the subset in ${\cal
G}(\underline{d})$, consisting of the sequences that do not
satisfy the condition (R3.2). Again we break ${\cal B}$ into
subsets:
$$
{\cal B}=\bigsqcup^{M-m^*}_{i=1}{\cal B}_i,
$$
where ${\cal B}_i$ has the same meaning as in Subsection 8.3 (but
for the multi-quadratic point $o$). Again the codimension of
${\cal B}_i$ in ${\cal G}(\underline{d})$ is bounded from below by
the number (\ref{17.10.22.1}), and for $k$ and $M$ fixed the worst
estimate corresponds to the case of equal or ``almost equal''
degrees $d_i$.

Arguing precisely in the same way as in the non-singular case
(Subsection 8.3), we see, since the dimension of the subspace
$\Pi$ is higher than in the non-singular case, that the claim of
Lemma 8.1 is true. Note that if $m^*=0$, then in the notations of
Subsection 8.3 we have $B(M-m^*)\geqslant 4$. Thus, replacing
$B(M-m^*)$ by 4, we get that $\mathop{\rm codim}({\cal
B}\subset{\cal G}(\underline{d}))$ is bounded from below by the
least of the three numbers
$$
{\mathop{\rm dim}\Pi+3-k\choose 2},\quad {\mathop{\rm
dim}\Pi+4-2k\choose 3},\quad {\mathop{\rm deg}h_{M-m^*}+4\choose
4}.
$$
The claim of Lemma 8.2 is true since, as we noted above,
$\mathop{\rm dim} \Pi$ in the multi-quadratic case is higher than
in the non-singular case. Obviously, $m^*<m_*$, so that
$\mathop{\rm deg} h_{M-m^*}\geqslant\mathop{\rm deg}h_{M-m_*}$ and
the claim of Lemma 8.3 is also true. As a result, we get the
inequality
$$
\mathop{\rm codim}({\cal B}\subset{\cal
G}(\underline{d}))\geqslant{M+2+l-k-\varepsilon(k)\choose 2},
$$
where $\mathop{\rm dim}\underline{\xi}=k-l$. The minimum of the
right hand side is attained for $l=2$ and it is easy to see that
this minimum is significantly higher than in the non-singular
case. It is easy to check, taking into account the variation of
the subspace $\Pi$, that
$$
\mathop{\rm codim}(B_{R3.2}\subset{\cal P}(o))>\gamma+M.
$$
This completes our consideration of the condition (R3.2) in the
multi-quadratic case.\vspace{0.3cm}


{\bf 8.5. The condition (R3.1).} In order to estimate the
codimension of the set ${\cal B}_{R3.1}(\underline{\xi},P)$, we
need the following known general fact. Take $e\geqslant 1$ and let
$\underline{w}=(w_1,\dots,w_e)\in{\mathbb Z}^e$ be a tuple of
integers, where $2\leqslant w_1\leqslant\dots\leqslant w_e$.

Set
$$
{\cal P}(\underline{w})=\prod^e_{i=1}{\cal P}_{w_i,N+1}
$$
to be the space of tuples $\underline{g}=(g_1,\dots,g_e)$ of
homogeneous polynomials in $N+1$ variables, $\mathop{\rm
deg}g_i=w_i$, which we consider as homogeneous polynomials on
${\mathbb P}^N$. Let
$$
{\cal B}^*(\underline{w})\subset{\cal P}(\underline{w})
$$
be the set of tuples $\underline{g}$, such that the scheme of
their common zeros is not an irreducible reduced subvariety of
codimension $e$ in ${\mathbb P}^N$.

{\bf Theorem 8.1.} {\it The following inequality holds:}
$$
\mathop{\rm codim}({\cal B}^*(\underline{w})\subset{\cal
P}(\underline{w}))\geqslant\frac12(N-e-1)(N-e-4)+2.
$$

{\bf Proof:} this is Theorem 2.1 in \cite{Pukh2022a}.

Let us estimate the codimension of the set ${\cal
B}_{R3.1}(\underline{\xi},\Pi)$. In order to do this, consider in
the projective space $P$ a hypersurface $P^{\sharp}$ that does not
contain the point $o$, for instance, the intersection of the
hyperplane ``at infinity'' with respect to the system of affine
coordinates $(z_1,\dots,z_{M+k})$ with the subspace $P$. If the
scheme of common zeros of the tuple of polynomials, consisting of
$$
f_1|_P,\dots,f_k|_P
$$
and the polynomials $f_{i,2}|_P$ for $i$ such that $d_i\geqslant
3$, is not an irreducible reduced subvariety of codimension
$k+k_{\geqslant 3}$ in $P$ (that is, the condition (R3.1) is
violated, see Subsection 1.4), then the scheme of common zeros of
the set of polynomials
\begin{equation}\label{27.02.23.1}
f_1|_{P^{\sharp}},\dots,f_k|_{P^{\sharp}},f_{i,2}|_{P^{\sharp}}\quad
\mbox{for} \quad d_i\geqslant 3,
\end{equation}
respectively, is reducible, non-reduced or is of codimension
$<k+k_{\geqslant 3}$ in $P^{\sharp}$. However, for each $i$, such
that $d_i\geqslant 3$, the homogeneous polynomials
$$
f_i|_{P^{\sharp}}=f_{i,d_i}|_{P^{\sharp}}\quad\mbox{and}\quad
f_{i,2}|_{P^{\sharp}}
$$
on the projective space $P^{\sharp}$ are linear combinations of
disjoint sets of monomials in $f_i$, so that the coefficients of
those polynomials belong to disjoint subsets of coefficients of
the polynomial $f_i$. Therefore (re-ordering the polynomials of
the tuple (\ref{27.02.23.1}) so that their degrees do not
decrease), applying Theorem 8.1 to the tuple (\ref{27.02.23.1}),
we get that the codimension of the set ${\cal
B}_{R3.1}(\underline{\xi},P)$ is at least
$$
\frac12(\mathop{\rm dim}P^{\sharp}-k-k_{\geqslant 3}-1)
(\mathop{\rm dim}P^{\sharp}-k-k_{\geqslant 3}-4)+2,
$$
where $\mathop{\rm dim}P^{\sharp}=M+l-\varepsilon(k)-1$,
$\mathop{\rm dim}\underline{\xi}=k-l$. It is easy to check by
elementary computations that this estimate (with the correction
due to the variation of the subspace $P$ and the set of linear
forms $\underline{\xi}$) is stronger than we need.

This completes the proof of the estimate for the codimension of
the complement ${\cal P}\backslash{\cal F}$ in Theorem 0.1.

Note that (for the technique of estimating the codimension that we
used) the estimate of Theorem 0.1 is optimal for the condition
(MQ2), that requirement turns out to be the strongest.


\begin{flushleft}
Department of Mathematical Sciences,\\
The University of Liverpool
\end{flushleft}

\noindent{\it pukh@liverpool.ac.uk}

\end{document}